\documentclass[11pt,a4paper]{article}
\usepackage[utf8]{inputenc}
\usepackage[T1]{fontenc}
\usepackage{amsfonts}
\usepackage{amssymb}
\usepackage{amsthm}
\usepackage{amsmath}
\usepackage{graphicx}
\usepackage{mathrsfs}
\usepackage[left=2cm,right=2cm,top=2cm,bottom=2cm]{geometry}

\newtheorem{theorem}{Theorem}
\newtheorem{definition}[theorem]{Definition}
\newtheorem{example}[theorem]{Example}

\newtheorem{lemma}[theorem]{Lemma}

\newtheorem{proposition}[theorem]{Proposition}
\newtheorem{remark}[theorem]{Remark}

\begin{document}

\title{Singular Optimal Controls for Stochastic Recursive Systems under
Convex Control Constraint}
\author{Liangquan Zhang$^{1}$\thanks{%
L. Zhang acknowledges the financial support partly by the National Nature
Science Foundation of China(Grant No. 11701040, 61871058 \&11871010) and the
Fundamental Research Funds for the Central Universities (No.2019XD-A11).
E-mail: xiaoquan51011@163.com.} \\
%EndAName
{\small 1. School of Science, }\\
{\small Beijing University of Posts and Telecommunications, }\\
{\small Beijing 100876, China}}
\maketitle

\begin{abstract}
In this paper, we study two kinds of singular optimal controls (SOCs for
short) problems where the systems governed by forward-backward stochastic
differential equations (FBSDEs for short), in which the control has two
components: the regular control, and the singular one. Both drift and
diffusion terms may involve the regular control variable. The regular
control domain is postulated to be convex. Under certain assumptions, in the
framework of the Malliavin calculus, we derive the pointwise second-order
necessary conditions for stochastic SOC in the classical sense. This
condition is described by two adjoint processes, a maximum condition on the
Hamiltonian supported by an illustrative example. A new necessary condition
for optimal singular control is obtained as well. Besides, as a by-product,
a verification theorem for SOCs is derived via viscosity solutions without
involving any derivatives of the value functions. It is worth pointing out
that this theorem has wider applicability than the restrictive classical
verification theorems. Finally, we focus on the connection between the
maximum principle and the dynamic programming principle for such SOCs
problem without the assumption that the value function is smooth enough.
%, the
%set inclusions among the sub- and super-jets of the value function and the
%first-order and second-order adjoint processes as well as the generalized
%Hamiltonian function are established.
\end{abstract}



\section{Introduction}

\label{sect:1}Singular stochastic control problem is a fundamental topic in
fields of stochastic control. This problem was first introduced by Bather
and Chernoff \cite{BC} in 1967 by considering a simplified model for the
control of a spaceship. It was then found that there was a connection
between the singular control and optimal stopping problem. This link was
established through the derivative of the value function of this initial
singular control problem and the value function of the corresponding optimal
stopping problem. Subsequently, it was considered by Ben\v{e}s, Shepp,
Witzsenhausen (see \cite{BSW}) and Karatzas and Shreve (see \cite{K1, K2,
KS1, KS2, KS3}).
%In contrast to classical control problems, singular control
%problems admit both of the continuity of the cumulative displacement of the
%state caused by control and the jump of one in impulsive control problems,
%between which it is either constant or absolutely continuous.

The state process is described by a $n$-dimensional SDE of the following
type:%
\begin{equation}
\left\{
\begin{array}{lll}
\mathrm{d}X^{t,x;v,\xi } & = & b\left( s,X^{t,x;v,\xi }\left( s\right)
,v\left( s\right) \right) \mathrm{d}s+\sigma \left( s,X^{t,x;v,\xi }\left(
s\right) ,v\left( s\right) \right) \mathrm{d}W\left( s\right) +G\left(
s\right) \mathrm{d}\xi \left( s\right) , \\
X^{t,x;v,\xi }\left( t\right) & = & x,\qquad 0\leq t\leq s\leq T,%
\end{array}%
\right.  \label{SDE1}
\end{equation}%
on some filtered probability space $\left( \Omega ,\mathcal{F},P\right) $,
where $b\left( \cdot ,\cdot ,\cdot \right) :\left[ 0,T\right] \times \mathbb{%
R}^{n}\times \mathbb{R}^{k}\rightarrow \mathbb{R}^{n},$ $\sigma \left( \cdot
,\cdot ,\cdot \right) :\left[ 0,T\right] \times \mathbb{R}^{n}\times \mathbb{%
R}^{k}\rightarrow \mathbb{R}^{n\times d},$ $G\left( \cdot \right) :\left[ 0,T%
\right] \rightarrow \mathbb{R}^{n\times m}$ are given deterministic
functions, $\left( W_{s}\right) _{s\geq 0}$ is an $d$-dimensional Brownian
motion, $\left( x,t\right) $ are initial time and state, $v\left( \cdot
\right) :\left[ 0,T\right] \rightarrow \mathbb{R}^{k}$ is a \textit{regular}
control process, and $\xi \left( \cdot \right) :\left[ 0,T\right]
\rightarrow \mathbb{R}^{m}$, with nondecreasing left-continuous with right
limits stands for the \textit{singular control}\footnote{%
Because the measure $\mathrm{d}\xi _{s}$ may be singular with respect to the
Lebesgue measure $\mathrm{d}s$.} (SC for short). To avoid the risk of
confusion, we shall introduce the other definitions of singular control in
various senses. Indeed, they are just a coincidence of terminology usage.

The aim is to minimize the cost functional:%
\begin{equation}
J\left( t,x;v,\xi \right) =\mathbb{E}\left[ \int_{t}^{T}l\left(
s,X^{t,x;v,\xi }\left( s\right) ,v\left( s\right) \right) \mathrm{d}%
s+\int_{t}^{T}K\left( s\right) \mathrm{d}\xi \left( s\right) \right] ,
\label{cost1}
\end{equation}%
where
\begin{eqnarray*}
l\left( \cdot ,\cdot ,\cdot \right) &:&\left[ 0,T\right] \times \mathbb{R}%
^{n}\times \mathbb{R}^{k}\rightarrow \mathbb{R}, \\
K\left( \cdot \right) &:&\left[ 0,T\right] \rightarrow \mathbb{R}%
_{+}^{m}\triangleq \left\{ x\in \mathbb{R}^{m}:x_{i}\geq 0,\text{ }%
i=1,\ldots m\right\}
\end{eqnarray*}%
are given deterministic functions, where $l\left( \cdot \right) $ represents
the running cost tare of the problem and $K$ the cost rate of applying the
singular control.

We mention that there are four approaches to deal with singular control: The
first, partial differential equations (PDE for short) and on variational
arguments, can be found in the works of Alvarez \cite{A1, A2}, Chow,
Menaldi, and Robin \cite{CMR}, Karatzas \cite{K2}, Karatzas and Shreve \cite%
{KS3}, and Menaldi and Taksar \cite{MT}. The second one is related to
probabilistic methods; see Baldursson \cite{B1}, Boetius \cite{Bo1, Bo2},
Boetius and Kohlmann \cite{BK1}, El Karoui and Karatzas \cite{EK1, EK2},
Karatzas \cite{K1}, and Karatzas and Shreve \cite{KS1, KS2}. Third, the DPP,
has been studied in a general context, for example, by Boetius \cite{Bo2},
Haussmann and Suo \cite{HS2}, Fleming and Soner \cite{FS} and Zhang \cite%
{Zhangsingular}. At last the maximum principle for optimal singular controls
(see, for example, Cadenillas and Haussmann \cite{CH}, Dufour and Miller
\cite{DM2}, Dahl and \O ksendal \cite{DO} see references therein).

Singular controls are used in diverse fields such as mathematical finance
(see Baldursson and Karatzas \cite{BK}, Chiarolla, Haussmann \cite{CHF},
Kobila \cite{KOB}, Karatzas, Wang \cite{KW}, Davis, Norman \cite{DN} and Pag%
\`{e}s and Possama\"{\i} \cite{PPs}), manufacturing systems (see, Shreve,
Lehoczky, and Gaver \cite{SLG}), and queuing systems (see Martins and
Kushner \cite{MK}).

Completely different from the singular control introduced above, to the best
of our knowledge, there are two other types of singular optimal controls, in
which the first-order necessary conditions turn out to be trivial. We list
briefly as follows:

\begin{itemize}
\item \textit{Singular optimal control in the classical sense (SOCCS for
short)}, is the optimal control for which the gradient and the Hessian of
the corresponding Hamiltonian with respect to the control variable
vanish/degenerate.

\item \textit{Singular optimal control in the sense of Pontryagin-type
maximum principle (SOCSPMP for short)}, is the optimal control for which the
corresponding Hamiltonian is equal to a constant in the control region.
\end{itemize}

When an optimal control is singular in certain senses above (SOCCS and
SOCSPMP), usually the first-order necessary condition could not carry
sufficient information for the further theoretical analysis and numerical
computation, and consequently it is necessary to investigate the second
order necessary conditions. In the deterministic setting, reader can refer
many articles in this direction (see \cite{BJaco, FT, GK, Goh, Kno, Kre,
Kush} and references therein).

As for the second-order necessary conditions for stochastic singular optimal
controls (SOCCS and SOCSPMP), there are some work should be mentioned, for
instance \cite{ZZconvex, ZZgen} (note that singular control $\xi \left(
\cdot \right) $ in these articles does not appear in systems). Tang \cite%
{Tang} obtained a pointwise second order maximum principle for stochastic
singular optimal controls in the sense of the Pontryagin-type maximum
principle whenever the control variable $u$ does not enter into the
diffusion term. Meanwhile, Tang addressed an integral-type second-order
necessary condition for stochastic optimal controls with convex control
constraints. Zhang and Zhang \cite{ZZconvex} also establish certain
pointwise second-order necessary conditions for stochastic singular (SOCCS)
optimal controls, in which both drift and diffusion terms in may depend on
the control variable $u$ with convex control region $U$ by making use of
Malliavin calculus technique. Later, adopting the same idea but with large
complicated analysis, Zhang et al. \cite{ZZgen} deepen this research for the
general case when the control region is nonconvex.

The theory of \textit{backward stochastic differential equation} (BSDE for
short) can be traced back to Bismut \cite{BJ, BJ2} who studied linear BSDE
motivated by stochastic control problems. Pardoux and Peng 1990 \cite{PP1}
proved the well-posedness for nonlinear BSDE. Duffie and Epstein (1992)
introduced the notion of recursive utilities in continuous time, which is
actually a type of BSDE where the generator $f$ is independent of $z$. El
Karoui et al. (1997, 2001) extended the recursive utility to the case where $%
f$ contains $z$. The term $z$ can be interpreted as an ambiguity aversion
term in the market (see Chen and Epstein 2002 \cite{CZJ}). Particularly, the
celebrated Black-Scholes formula indeed provided an effective way of
representing the option price (which is the solution to a kind of linear
BSDE) through the solution to the Black-Scholes equation (parabolic partial
differential equation actually). Since then, BSDE has been extensively
studied and used in the areas of applied probability and optimal stochastic
controls, particularly in financial engineering (cf for instance \cite{KPQ}).

By means of BSDE, Peng (1990) \cite{Peng1990} considered the following type
of stochastic optimal control problem: Minimize a cost function
\begin{equation*}
J\left( v\left( \cdot \right) \right) =\mathbb{E}\int_{0}^{T}l\left( x\left(
t\right) ,v\left( t\right) \right) \text{d}t+\mathbb{E}\left( h\left(
T\right) \right) ,
\end{equation*}%
subject to
\begin{equation}
\left\{
\begin{array}{rcl}
\text{d}x\left( t\right) & = & g\left( t,x\left( t\right) ,v\left( t\right)
\right) \text{d}t+\sigma \left( t,x\left( t\right) ,v\left( t\right) \right)
\text{d}W\left( t\right) , \\
x\left( 0\right) & = & x_{0},%
\end{array}%
\right.  \label{1.1}
\end{equation}%
over an admissible control domain which need not be convex, and the
diffusion coefficients depends on the control variable. In his paper, by
spike variational method and the second order adjoint equations, Peng \cite%
{Peng1990} obtained a general stochastic maximum principle for the above
optimal control problem. It was just the adjoint equations in stochastic
optimal control problems that motivated the famous theory of BSDE (cf \cite%
{PP1}).

Later, Peng first \cite{Peng1993} studied a stochastic optimal control
problem where state variables are described by the system of FBSDEs:
\begin{equation}
\left\{
\begin{array}{rcl}
\text{d}x\left( t\right) & = & f\left( t,x\left( t\right) ,v\left( t\right)
\right) \text{d}t+\sigma \left( t,x\left( t\right) ,v\left( t\right) \right)
\text{d}W\left( t\right) , \\
\text{d}y\left( t\right) & = & g\left( t,x\left( t\right) ,v\left( t\right)
\right) \text{d}t+z\left( t\right) \text{d}W\left( t\right) , \\
x\left( 0\right) & = & x_{0},y\left( T\right) =y,%
\end{array}%
\right.  \label{1.2}
\end{equation}%
where $x$ and $y$ are given deterministic constants. The optimal control
problem is to minimize the cost function:
\begin{equation*}
J\left( v\left( \cdot \right) \right) =\mathbb{E}\left[ \int_{0}^{T}l\left(
t,x\left( t\right) ,y\left( t\right) ,v\left( t\right) \right) \text{d}%
t+h\left( x\left( T\right) \right) +\gamma \left( y\left( 0\right) \right) %
\right] ,
\end{equation*}%
over an admissible control domain which is convex. Later, Xu \cite{XWS}
studied the following non-fully coupled forward-backward stochastic control
system:
\begin{equation}
\left\{
\begin{array}{rcl}
\text{d}x\left( t\right) & = & f\left( t,x\left( t\right) ,v\left( t\right)
\right) \text{d}t+\sigma \left( t,x\left( t\right) \right) \text{d}W\left(
t\right) , \\
\text{d}y\left( t\right) & = & g\left( t,x\left( t\right) ,y\left( t\right)
,z\left( t\right) ,v\left( t\right) \right) \text{d}t+z\left( t\right) \text{%
d}W\left( t\right) , \\
x\left( 0\right) & = & x_{0},y\left( T\right) =h\left( x\left( T\right)
\right) .%
\end{array}%
\right.  \label{1.3}
\end{equation}%
The optimal control problem is to minimize the cost function $J\left(
v_{\left( \cdot \right) }\right) =\mathbb{E}\gamma \left( y_{0}\right) ,$
over $\mathcal{U}_{ad},$ but the control domain is non-convex. Wu \cite%
{Wu1998} firstly gave the maximum principle for optimal control problem of
fully coupled forward-backward stochastic system:
\begin{equation}
\left\{
\begin{array}{rcl}
\text{d}x\left( t\right) & = & f\left( t,x\left( t\right) ,y\left( t\right)
,z\left( t\right) ,v\left( t\right) \right) \text{d}t+\sigma \left(
t,x\left( t\right) ,y\left( t\right) ,z\left( t\right) ,v\left( t\right)
\right) \text{d}W\left( t\right) , \\
\text{d}y\left( t\right) & = & -g\left( t,x\left( t\right) ,y\left( t\right)
,z\left( t\right) ,v\left( t\right) \right) \text{d}t+z\left( t\right) \text{%
d}W\left( t\right) , \\
x\left( 0\right) & = & x_{0},\quad \quad y\left( T\right) =\xi ,%
\end{array}%
\right.  \label{1.4}
\end{equation}%
where $\xi $ is a random variable and the cost function:
\begin{equation*}
J\left( v\left( \cdot \right) \right) =\mathbb{E}\left[ \int_{0}^{T}L\left(
t,x\left( t\right) ,y\left( t\right) ,z\left( t\right) ,v\left( t\right)
\right) \text{d}t+\Phi \left( x\left( T\right) \right) +h\left( y\left(
0\right) \right) \right] .
\end{equation*}%
The optimal control problem is to minimize the cost function $J\left(
v_{\left( \cdot \right) }\right) $ over an admissible control domain which
is convex. Ji and Zhou \cite{JZ} obtained a maximum principle for stochastic
optimal control of non-fully coupled forward-backward stochastic system with
terminal state constraints. Shi and Wu \cite{SW} studied the maximum
principle for fully coupled forward-backward stochastic system:
\begin{equation}
\left\{
\begin{array}{rcl}
\text{d}x_{t} & = & b\left( t,x\left( t\right) ,y\left( t\right) ,z\left(
t\right) ,v\left( t\right) \right) \text{d}t+\sigma \left( t,x\left(
t\right) ,y\left( t\right) ,z\left( t\right) \right) \text{d}B_{t}, \\
\text{d}y_{t} & = & -f\left( t,x\left( t\right) ,y\left( t\right) ,z\left(
t\right) ,v\left( t\right) \right) \text{d}t+z\left( t\right) \text{d}B_{t},
\\
x\left( 0\right) & = & x_{0},\quad y\left( T\right) =h\left( x\left(
T\right) \right) .%
\end{array}%
\right.  \label{1.5}
\end{equation}%
and the cost function is
\begin{equation*}
J\left( v\left( \cdot \right) \right) =\mathbb{E}\left[ \int_{0}^{T}l\left(
t,x\left( t\right) ,y\left( t\right) ,z\left( t\right) ,v\left( t\right)
\right) \text{d}t+\Phi \left( x\left( T\right) \right) +\gamma \left(
y\left( 0\right) \right) \right] .
\end{equation*}%
The control domain is non-convex but the forward diffusion does not contain
the control variable.

Subsequently, in order to study the backward linear-quadratic optimal
control problem, Kohlmann and Zhou \cite{KZ}, Lim and Zhou \cite{LZ}
developed a new method for handling this problem. The term $z$ is regarded
as a control process and the terminal condition $y_{T}=h\left( x_{T}\right) $
as a constraint, and then it is possible to use the Ekeland variational
principle to obtain the maximum principle. Adopting this idea, Yong \cite%
{Yong2010} and Wu \cite{Wu2013} independently established the maximum
principle for the recursive stochastic optimal control problem (noting the
diffusion term containing control variable with non-convex control region).
Nonetheless, the maximum principle derived by these method involves two
unknown parameters. Therefore, the hard questions raise as follows: What is
the second-order variational equation for the BSDE? How to obtain the
second-order adjoint equation since the quadratic form with respect to the
variation of $z$. All of which seem to be extremely complicated.

Hu \cite{Hms} overcomes the above difficulties by introducing two new
adjoint equations. Then, the second-order variational equation for the BSDE
and the maximum principle are obtained. The main difference of his
variational equations with those in Peng \cite{Peng1990} consists in the
term $\left\langle p\left( t\right) ,\delta \sigma \left( t\right)
\right\rangle I_{E_{\varepsilon }}\left( t\right) $ in the variation of $z$.
Due to the term $\left\langle p\left( t\right) ,\delta \sigma \left(
t\right) \right\rangle I_{E_{\varepsilon }}\left( t\right) $ in the
variation of $z$, Hu obtained a global maximum principle which is novel and
different from that in Wu \cite{Wu2013}, Yong \cite{Yong2010} and previous
work, which solves completely Peng's open problem. Furthermore, Hu's maximum
principle is stronger than the one in Wu \cite{Wu2013}, Yong \cite{Yong2010}%
. For a general case, reader can refer \cite{HJX}.

Motivated by above work, in this paper, we consider singular controls
problem of the following type:%
\begin{equation}
\left\{
\begin{array}{rcl}
\mathrm{d}X^{t,x;v,\xi }\left( s\right)  & = & b\left( s,X^{t,x;v,\xi
}\left( s\right) ,v\left( s\right) \right) \mathrm{d}s+\sigma \left(
s,X^{t,x;v,\xi }\left( s\right) ,v\left( s\right) \right) \mathrm{d}W\left(
s\right) +G\left( s\right) \mathrm{d}\xi \left( s\right) , \\
\mathrm{d}Y^{t,x;v,\xi }\left( s\right)  & = & -f\left( t,X^{t,x;v,\xi
}\left( s\right) ,Y^{t,x;v,\xi }\left( s\right) ,Z^{t,x;v,\xi }\left(
s\right) ,v\left( s\right) \right) \mathrm{d}s \\
&  & +Z^{t,x;v,\xi }\left( s\right) \mathrm{d}W\left( s\right) -K\mathrm{d}%
\xi \left( s\right) , \\
X^{t,x;v,\xi }\left( t\right)  & = & x,\text{ }Y^{t,x;v,\xi }\left( T\right)
=\Phi \left( X^{t,x;v,\xi }\left( T\right) \right) ,\qquad 0\leq t\leq s\leq
T,%
\end{array}%
\right.   \label{FBSDE1}
\end{equation}%
with the similar cost functional%
\begin{equation}
J\left( t,x;v,\xi \right) =\left. Y^{t,x;v,\xi }\left( s\right) \right\vert
_{s=t}.  \label{aim}
\end{equation}%
%
%
%
%
%
%
%
%
%
%
%
%
%
%
%
%
%
%
%
%
%
%
%
%
%
%
%
%
%
%
%
%
%
%
%
%
%
%
%
%
%
%
%
%
%
%
%
%
%
%
%
%
%
%
%
%
%
%
%
%
%
%
%
%
%
%
%
%
%
%
%
%
%
%
%
%
%
%
%
%
%We postulate that $K$ is a deterministic matrix in Eq. (\ref{FBSDE1}). The
%justification will be given in Remark \ref{GK} below.
Wang \cite{WB} firstly introduced and studied a class of singular control
problems with recursive utility, where the cost function is determined by
BSDE. Under certain assumptions, the author proved that the value function
is a nonnegative, convex solution of the H-J-B equation. However, FBSDEs in
Wang \cite{WB} do not contain the regular control and the generator is not
general case. In our work, using some properties of the BSDE and analysis
technique, we expand the extension of the MP for SOC to the recursive
control problem in Zhang and Zhang \cite{ZZconvex}. To the best of our
knowledge, such singular optimal controls problems of FBSDEs (\ref{FBSDE1})
via two kinds of singular controls have not been explored before. We shall
establish some pointwise second-order necessary conditions for stochastic
optimal controls of FBSDEs. Both drift and diffusion terms may contain the
control variable $u$, and we assume that the control region $U$ is convex.
We also consider the pointwise second-order necessary condition, which is
easier to verify in practical applications.

As claimed in \cite{ZZconvex}, quite different from the deterministic
setting, there exist some essential difficulties in deriving the pointwise
second-order necessary condition from an integral-type one whenever the
diffusion term depends on the control variable, even for the case of convex
control domain. We overcome these difficulties by means of some technique
from the Malliavin calculus. For general case, namely, the control region is
non-convex can be found in \cite{ZZgen}.

In this paper, we are interested in studying singular optimal controls for
FBSDEs (\ref{FBSDE1}). Compared with above literature, our paper has several
new features. The novelty of the formulation and the contribution in this
paper may be stated as follows:

\begin{itemize}
\item Our control systems in this paper are governed by FBSDEs which exactly
extends the work of Zhang and Zhang \cite{ZZconvex} to utilities. Our work
is the first time to establish the pointwise second order necessary
condition for stochastic singular optimal control in the classical sense for
FBSDEs, a new necessary condition for singular control is involved as well.
In this sense, our paper actually considers two kinds of singular controls
problems simultaneously, which is interesting to deepen this research.

\item We derive a new verification theorem for optimal singular controls via
viscosity solution, which responses to the question raised in Zhang \cite%
{Zhangsingular}; Meanwhile, we study the relationship between the adjoint
equations derived and value function, which extends the smooth case
considered by Cadenillas and Haussmann \cite{CH} to the framework of
viscosity solution for stochastic recursive systems.
\end{itemize}

The rest of this paper is organized as follows: after some preliminaries in
the second section, we are devoted the third section to the MP for two kinds
of singular optimal controls. A concrete example is concluded with as well.
Then, in Section \ref{secverif}, we study the verification theorem for
singular optimal controls via viscosity solutions. Finally, we establish the
relationship between the DPP and MP for viscosity solution. Some proofs of
lemmas are displayed in Appendix \ref{APP}.

\section{Preliminaries and Notations}

\label{sect2}

Throughout this paper, we denote by $\mathbb{R}^{n}$ the space of $n$%
-dimensional Euclidean space, by $\mathbb{R}^{n\times d}$ the space the
matrices with order $n\times d$. Let $(\Omega ,\mathcal{F},\{\mathcal{F}%
_{t}\}_{t\geq 0},P)$ be a complete filtered probability space on which a
\textit{one}-dimensional standard Brownian motion $W(\cdot )$ is defined,
with $\{\mathcal{F}_{t}\}_{t\geq 0}$ being its natural filtration, augmented
by all the $P$-null sets. Given a subset $U$ (nonempty, bounded, and convex)
of $\mathbb{R}^{k},$ we will denote $\mathcal{U}\left[ 0,T\right] \mathcal{=U%
}_{1}\times \mathcal{U}_{2}$, separately, the class of measurable, adapted
processes $\left( v,\xi \right) :\left[ 0,T\right] \times \Omega \rightarrow
U\times \left[ 0,\infty \right) ^{m},$ with $\xi $ nondecreasing
left-continuous with right limits and $\xi _{0}=0$, moreover, $\mathbb{E}%
\left[ \sup\limits_{0\leq t\leq T}\left\vert v\left( t\right) \right\vert
^{2}+\left\vert \xi \left( T\right) \right\vert ^{2}\right] <\infty .$ $\xi $
is called \textit{singular control.} %\footnote{%
For each $t>0$, we denote by $\left\{ \mathcal{F}_{s}^{t},t\leq s\leq
T\right\} $ the natural filtration of the Brownian motion $\{{W}\left(
s\right) {-W}\left( t\right) {\}}_{{t\leq s\leq T}}$, augmented by the $P$%
-null sets of $\mathcal{F}$. $\top $ appearing as superscript denotes the
transpose of a matrix. In what follows, $C$ represents a generic constant,
which can be different from line to line.

We now introduce the following spaces of processes:
\begin{align*}
\mathcal{S}^{2}(0,T;\mathbb{R})\triangleq & \left\{ \mathbb{R}^{n}\text{%
-valued }\mathcal{F}_{t}\text{-adapted process }\phi (t)\text{; }\mathbb{E}%
\left[ \sup\limits_{0\leq t\leq T}\left\vert \phi \left( t\right)
\right\vert ^{2}\right] <\infty \right\} , \\
\mathcal{M}^{2}(0,T;\mathbb{R})\triangleq & \left\{ \mathbb{R}^{n}\text{%
-valued }\mathcal{F}_{t}\text{-adapted process }\varphi (t)\text{; }\mathbb{E%
}\left[ \int_{0}^{T}\left\vert \varphi \left( t\right) \right\vert ^{2}%
\mbox{\rm d}t\right] <\infty \right\} ,
\end{align*}%
and denote $\mathcal{N}^{2}\left[ 0,T\right] =\mathcal{S}^{2}(0,T;\mathbb{R}%
^{n})\times \mathcal{S}^{2}(0,T;\mathbb{R})\times \mathcal{M}^{2}(0,T;%
\mathbb{R}^{n}).$ Clearly, $\mathcal{N}^{2}\left[ 0,T\right] $ forms a
Banach space.

For any $v\left( \cdot \right) \times \xi \left( \cdot \right) \in \mathcal{U%
}_{1}\times \mathcal{U}_{2},$ we study the stochastic control systems
governed by FBSDEs (\ref{FBSDE1}).%\begin{equation}
%\left\{
%\begin{array}{rcl}
%\mathrm{d}X_{s}^{t,x;v,\xi } & = & b\left( s,X_{s}^{t,x;v,\xi },v_{s}\right)
%\mathrm{d}s+\sigma \left( s,X_{s}^{t,x;v,\xi },v_{s}\right) \mathrm{d}%
%W_{s}+G_{s}\mathrm{d}\xi _{s}, \\
%\mathrm{d}Y_{s}^{t,x;v,\xi } & = & -f\left( s,X_{s}^{t,x;v,\xi
%},Y_{s}^{t,x;v,\xi },Z_{s}^{t,x;v,\xi },v_{s}\right) \mathrm{d}%
%s+Z_{s}^{t,x;v,\xi }\mathrm{d}W_{s}-K_{s}\mathrm{d}\xi _{s}, \\
%X_{t}^{t,x;v,\xi } & = & x,\text{ }Y_{T}^{t,x;v,\xi }=\Phi \left(
%X_{T}^{t,x;v,\xi }\right) ,\qquad 0\leq t\leq s\leq T,%
%\end{array}%
%\right.  \label{eq1}
%\end{equation}%
\newline

We assume that the following conditions hold:

\begin{description}
\item[(A1)] The coefficients $b:[0,T]\times \mathbb{R}^{n}\times \mathbb{R}%
^{k}\rightarrow \mathbb{R}^{n},$ $\sigma :[0,T]\times \mathbb{R}^{n}\times
\mathbb{R}^{k}\rightarrow \mathbb{R}^{n},$ are twice continuously
differentiable with respect to $x;$ $b,$ $b_{x},$ $b_{xx},$ $\sigma ,$ $%
\sigma _{x},$ $\sigma _{xx}$ are continues in $\left( x,u\right) ;$ $b_{x},$
$b_{xx},$ $\sigma _{x},$ $\sigma _{xx}$ are bounded $b$, $\sigma $ are
bounded by $C\left( 1+\left\vert x\right\vert +\left\vert u\right\vert
\right) $ for some positive constant $C.$ Moreover, for any $\left(
t,x_{1},u_{1}\right) ,$ $\left( t,x_{2},u_{2}\right) \in \left[ 0,T\right]
\times \mathbb{R}^{n}\times \mathbb{R}^{k},$%
\begin{eqnarray*}
\left\vert b\left( t,0,x\right) \right\vert +\left\vert \sigma \left(
t,0,u\right) \right\vert &\leq &C, \\
\left\vert b_{\left( x,u\right) ^{2}}\left( t,x_{1},u_{1}\right) -b_{\left(
x,u\right) ^{2}}\left( t,x_{2},u_{2}\right) \right\vert &\leq &C\left(
\left\vert x_{1}-x_{2}\right\vert +\left\vert u_{1}-u_{2}\right\vert \right)
, \\
\left\vert \sigma _{\left( x,u\right) ^{2}}\left( t,x_{1},u_{1}\right)
-\sigma _{\left( x,u\right) ^{2}}\left( t,x_{2},u_{2}\right) \right\vert
&\leq &C\left( \left\vert x_{1}-x_{2}\right\vert +\left\vert
u_{1}-u_{2}\right\vert \right) .
\end{eqnarray*}

\item[(A2)] The coefficients $f:[0,T]\times \mathbb{R}^{n}\times \mathbb{R}%
\times \mathbb{\mathbb{R}}\times \mathbb{R}^{k}\rightarrow \mathbb{R},$ $%
\Phi :\mathbb{R}^{n}\rightarrow \mathbb{R},$ are twice continuously
differentiable with respect to $\left( x,y,z\right) .$ $K$ is a given
deterministic matrix. $f,$ $Df,$ $D^{2}f$ are continuous in $\left(
x,y,z,u\right) $. There exists constant $C>0$ such that for any $\left(
t,x_{1},y_{1},z_{1},u_{1}\right) ,$ $\left( t,x_{2},y_{2},z_{2},u_{2}\right)
\in \left[ 0,T\right] \times \mathbb{R}^{n}\times \mathbb{R}\times \mathbb{R}%
\times \mathbb{R}^{k},$
\begin{equation*}
\left\vert f\left( t,x,y,z,u\right) \right\vert \leq C\left( 1+\left\vert
x\right\vert +\left\vert y\right\vert +\left\vert z\right\vert \right) ,
\end{equation*}%
\begin{equation*}
\left\vert f_{x}\left( t,x,y,z,u\right) \right\vert +\left\vert f_{y}\left(
t,x,y,z,u\right) \right\vert +\left\vert f_{z}\left( t,x,y,z,u\right)
\right\vert +\left\vert f_{u}\left( t,x,y,z,u\right) \right\vert \leq C,
\end{equation*}%
\begin{equation*}
\begin{array}{l}
\left\vert f_{xx}\left( t,x,y,z,u\right) \right\vert +\left\vert
f_{xu}\left( t,x,y,z,u\right) \right\vert +\left\vert f_{yu}\left(
t,x,y,z,u\right) \right\vert \\
\qquad +\left\vert f_{yy}\left( t,x,y,z,u\right) \right\vert +\left\vert
f_{zz}\left( t,x,y,z,u\right) \right\vert +\left\vert f_{zu}\left(
t,x,y,z,u\right) \right\vert \\
\qquad +\left\vert f_{uu}\left( t,x,y,z,u\right) \right\vert \leq C,%
\end{array}%
\end{equation*}%
\begin{equation*}
\begin{array}{l}
\left\vert f_{\left( x,y,z,u\right) ^{2}}\left(
t,x_{1},y_{1},z_{1},u_{1}\right) -f_{\left( x,y,z,u\right) ^{2}}\left(
t,x_{2},y_{2},z_{2},u_{2}\right) \right\vert \\
\qquad \leq C\left( \left\vert x_{1}-x_{2}\right\vert +\left\vert
y_{1}-y_{2}\right\vert +\left\vert z_{1}-z_{2}\right\vert +\left\vert
u_{1}-u_{2}\right\vert \right) ,%
\end{array}%
\end{equation*}%
and%
\begin{equation*}
\begin{array}{l}
\Phi \left( x\right) \leq C\left( 1+\left\vert x\right\vert ^{2}\right) ,%
\text{ }\Phi _{x}\left( x\right) \leq C\left( 1+\left\vert x\right\vert
\right) , \\
\Phi _{xx}\left( x\right) \leq C,\text{ }\left\vert \Phi _{xx}\left(
x_{1}\right) -\Phi _{xx}\left( x_{2}\right) \right\vert \leq C\left\vert
x_{1}-x_{2}\right\vert .%
\end{array}%
\end{equation*}
\end{description}

Under above assumptions (A1)-(A2), for any $v\left( \cdot \right) \times \xi
\left( \cdot \right) \in \mathcal{U}_{1}\times \mathcal{U}_{2}$, it is easy
to check that FBSDEs (\ref{FBSDE1}) admit a unique $\mathcal{F}_{t}$-adapted
solution denoted by the triple $(X_{\cdot }^{t,x;v,\xi },Y_{\cdot
}^{t,x;v,\xi },Z_{\cdot }^{t,x;v,\xi })\in \mathcal{N}^{2}\left[ 0,T\right] $
(See Pardoux and Peng \cite{PP1}).

Like Peng \cite{P3}, given any control processes $v\left( \cdot \right)
\times \xi \left( \cdot \right) \in \mathcal{U}_{1}\times \mathcal{U}_{2}$,
we introduce the following cost functional:
\begin{equation}
J(t,x;v\left( \cdot \right) ,\xi \left( \cdot \right) )=\left.
Y_{s}^{t,x;v,\xi }\right\vert _{s=t},\qquad \left( t,x\right) \in \left[ 0,T%
\right] \times \mathbb{R}^{n}.  \label{c1}
\end{equation}%
We are interested in the \emph{value function} of the stochastic optimal
control problem:
\begin{eqnarray}
u\left( t,x\right) &=&J(t,x;\hat{v}\left( \cdot \right) ,\hat{\xi}\left(
\cdot \right) )  \notag \\
&=&ess\inf_{v\left( \cdot \right) \times \xi \left( \cdot \right) \in
\mathcal{U}_{1}\times \mathcal{U}_{2}}J\left( t,x;v\left( \cdot \right) ,\xi
\left( \cdot \right) \right) ,\text{ }\left( t,x\right) \in \left[ 0,T\right]
\times \mathbb{R}^{n}.  \label{value}
\end{eqnarray}%
Since the value function (\ref{value}) is defined by the solution of
controlled FBSDEs (\ref{FBSDE1}), so from the existence and uniqueness, $u$
is \emph{well-defined}.

The following estimate is very useful whose proof can be found in Briand et
al. 2003 \cite{BDH}.

\begin{lemma}
\label{estlemma} Let $\left( y^{i},z^{i}\right) ,$ $i=1,2,$ be the solution
to the following
\begin{equation}
y^{i}\left( t\right) =\xi ^{i}+\int_{t}^{T}f^{i}\left( s,y^{i}\left(
s\right) ,z^{i}\left( s\right) \right) \mathrm{d}s-\int_{t}^{T}z^{i}\left(
s\right) \mathrm{d}W\left( s\right) ,  \label{estbdsde}
\end{equation}%
where $\xi ^{i}\in \mathcal{F}_{T}$ and $\mathbb{E}\left[ \left\vert \xi
^{i}\right\vert ^{\beta }\right] <\infty ,$ whilst $f^{i}\left(
s,y^{i},z^{i}\right) $ satisfies the conditions \emph{(A2)}\textsl{,} and
\begin{equation*}
\mathbb{E}\left[ \left( \int_{t}^{T}\left\vert f^{i}\left( s,y^{i}\left(
s\right) ,z^{i}\left( s\right) \right) \right\vert \mathrm{d}s\right)
^{\beta }\right] <\infty ;
\end{equation*}%
Then, for some $\beta \geq 2,$ there exists a positive constant $C_{\beta }$
such that
\begin{eqnarray*}
&&\mathbb{E}\left[ \sup_{0\leq t\leq T}\left\vert y^{1}\left( t\right)
-y^{2}\left( t\right) \right\vert ^{\beta }+\left( \int_{0}^{T}\left\vert
z^{1}\left( s\right) -z^{2}\left( s\right) \right\vert ^{2}\mathrm{d}%
s\right) ^{\frac{\beta }{2}}\right] \\
&\leq &C_{\beta }\mathbb{E}\Bigg [\left\vert \xi ^{1}-\xi ^{2}\right\vert
^{\beta }+\left( \int_{t}^{T}\left\vert f^{1}\left( s,y^{1}\left( s\right)
,z^{1}\left( s\right) \right) -f^{2}\left( s,y^{2}\left( s\right)
,z^{2}\left( s\right) \right) \right\vert \mathrm{d}s\right) ^{\beta }\Bigg ]%
.
\end{eqnarray*}%
Particularly, whenever putting $\xi ^{2}=0,$ and $f^{2}=0,$ one has%
\begin{equation*}
\mathbb{E}\left[ \sup_{0\leq t\leq T}\left\vert y^{1}\left( t\right)
\right\vert ^{\beta }+\left( \int_{0}^{T}\left\vert z^{1}\left( s\right)
\right\vert ^{2}\mathrm{d}s\right) ^{\frac{\beta }{2}}\right] \leq C_{\beta }%
\mathbb{E}\Bigg [\left\vert \xi ^{1}\right\vert ^{\beta }+\left(
\int_{t}^{T}\left\vert f^{1}\left( s,0,0\right) \right\vert \mathrm{d}%
s\right) ^{\beta }\Bigg ].
\end{equation*}
\end{lemma}

Now let us recall briefly the notion of differentiation on Wiener space (see
the expository papers by Nualart 1995 \cite{Nua1995}, Nualart and Pardoux
\cite{NuaPar} and Ocone 1988 \cite{Oco}).

\begin{itemize}
\item $C_{b}^{k}\left( \mathbb{R}^{k},\mathbb{R}^{q}\right) $ will denote
the set of functions of class $C^{k}$ from $\mathbb{R}^{k}$ into $\mathbb{R}%
^{q}$ whose partial derivatives of order less than or equal to $k$ are
bounded.

\item Let $\mathcal{S}$ denote the set of random variables $\xi $ of the
form $\xi =\varphi (W\left( h^{1}\right) ,W\left( h^{2}\right) ,\cdots
,W\left( h^{k}\right) ),$ where $\varphi \in C_{b}^{\infty }\left( \mathbb{R}%
^{k},\mathbb{R}\right) $, $h^{1},h^{2},\cdots h^{k}\in L^{2}\left( \left[ 0,T%
\right] ;\mathbb{R}^{n}\right) ,$ and $W\left( h^{i}\right)
=\int_{0}^{T}\left\langle h_{s}^{i},\mathrm{d}W\left( s\right) \right\rangle
$.

\item If $\xi \in \mathcal{S}$ is of the above form, we define its
derivative as being the $n$-dimensional process%
\begin{equation*}
\mathcal{D}_{\theta }\xi =\sum_{j=1}^{k}\frac{\partial }{\partial x_{j}}%
\varphi \left( W\left( h^{1}\right) ,W\left( h^{2}\right) ,\cdots ,W\left(
h^{k}\right) \right) h_{\theta }^{j},\text{ }0\leq \theta \leq T.
\end{equation*}%
For $\xi \in \mathcal{S},$ $p>1,$ we define the norm
\begin{equation*}
\left\Vert \xi \right\Vert _{1,p}=\left[ E\left\{ \left\vert \xi \right\vert
^{p}+\left( \int_{0}^{T}\left\vert \mathcal{D}_{\theta }\xi \right\vert ^{2}%
\mathrm{d}\theta \right) ^{\frac{p}{2}}\right\} \right] ^{\frac{1}{p}}.
\end{equation*}
\end{itemize}

It can be shown (Nualart 1995) that the operator $\mathcal{D}$ has a closed
extension to the space $\mathbb{D}^{1,p}$, the closure of $\mathcal{S}$ with
respect to the norm $\left\Vert \cdot \right\Vert _{1,p}$. Observe that if $%
\xi $ is $\mathcal{F}_{t}$-measurable, then $\mathcal{D}_{\theta }\xi =0$
for $\theta \in \left( t,T\right] $. We denote by $\mathcal{D}_{\theta
}^{i}\xi $, the i\emph{th} component of $\mathcal{D}_{\theta }\xi $.

Let $\mathbb{L}^{1,p}\left( \mathbb{R}^{d}\right) $ denote the set of $%
\mathbb{R}^{d}$-valued progressively measurable processes $\{u\left(
t,\omega \right) ,0\leq t\leq T;\omega \in \Omega \}$ such that

\begin{itemize}
\item For a.e. $t\in \left[ 0,T\right] ,$ $u\left( t,\cdot \right) \in
\mathbb{D}^{1,p}\left( \mathbb{R}^{n}\right) ;$

\item $\left( t,\omega \right) \rightarrow \mathcal{D}u\left( t,\omega
\right) \in \left( L^{2}\left( \left[ 0,T\right] \right) \right) ^{n\times
d} $ admits a progressively measurable version;

\item We have%
\begin{equation*}
\left\Vert u\right\Vert _{1,p}=\mathbb{E}\left[ \left(
\int_{0}^{T}\left\vert u\left( t\right) \right\vert ^{2}\mathrm{d}t\right) ^{%
\frac{p}{2}}+\left( \int_{0}^{T}\int_{0}^{T}\left\vert \mathcal{D}_{\theta
}u\left( t\right) \right\vert ^{2}\mathrm{d}\theta \mathrm{d}t\right) ^{%
\frac{p}{2}}\right] <+\infty .
\end{equation*}
\end{itemize}

Note that for each $\left( \theta ,t,\omega \right) ,$ $\mathcal{D}_{\theta
}u\left( t,\omega \right) $ is an $n\times d$ matrix. Hence, $\left\vert
\mathcal{D}_{\theta }u\left( t\right) \right\vert ^{2}=\sum_{i,j}\left\vert
\mathcal{D}_{\theta }^{i}u_{j}\left( t\right) \right\vert ^{2}.$ Obviously, $%
\mathcal{D}_{\theta }u\left( t,\omega \right) $ is defined uniquely up to
sets of $\mathrm{d}\theta \otimes \mathrm{d}t\otimes \mathrm{d}P$ measure
zero. Moreover, denote by $\mathbb{L}_{\mathbb{F}}^{1,p}\left( \mathbb{R}%
^{d}\right) $ the set of all adapted processes in $\mathbb{L}^{1,p}\left(
\mathbb{R}^{d}\right) .$

We define the following notations from Zhang and Zhang \cite{ZZconvex}:%
\begin{equation*}
\begin{array}{lll}
\mathbb{L}_{2+}^{1,p}\left( \mathbb{R}^{d}\right) & := & \Big \{\left.
\varphi \left( \cdot \right) \in \mathbb{L}^{1,p}\left( \mathbb{R}%
^{d}\right) \right\vert \exists \mathcal{D}^{+}\varphi \left( \cdot \right)
\in L^{2}\left( \left[ 0,T\right] \times \Omega ;\mathbb{R}^{n}\right) \text{
such that} \\
&  & f_{\varepsilon }\left( s\right) :=\sup_{s<t<\left( s+t\right) \wedge T}%
\mathbb{E}\left\vert \mathcal{D}_{s}\varphi \left( t\right) -\mathcal{D}%
^{+}\varphi \left( s\right) \right\vert ^{2}<\infty ,\text{ }a.e.\text{ }%
s\in \left[ 0,T\right] , \\
&  & f_{\varepsilon }\left( s\right) \text{ is measurable on }\left[ 0,T%
\right] \text{ for any }\varepsilon >0,\text{ and }\lim_{\varepsilon
\rightarrow 0^{+}}\int_{0}^{T}f_{\varepsilon }\left( s\right) \mathrm{d}s=0%
\Big \};%
\end{array}%
\end{equation*}%
and%
\begin{equation*}
\begin{array}{lll}
\mathbb{L}_{2-}^{1,p}\left( \mathbb{R}^{d}\right) & := & \Big \{\left.
\varphi \left( \cdot \right) \in \mathbb{L}^{1,p}\left( \mathbb{R}%
^{d}\right) \right\vert \exists \mathcal{D}^{-}\varphi \left( \cdot \right)
\in L^{2}\left( \left[ 0,T\right] \times \Omega ;\mathbb{R}^{n}\right) \text{
such that} \\
&  & g_{\varepsilon }\left( s\right) :=\sup_{\left( s-\varepsilon \right)
\vee 0<t<s}\mathbb{E}\left\vert \mathcal{D}_{s}\varphi \left( t\right) -%
\mathcal{D}^{-}\varphi \left( s\right) \right\vert ^{2}<\infty ,\text{ }a.e.%
\text{ }s\in \left[ 0,T\right] , \\
&  & g_{\varepsilon }\left( s\right) \text{ is measurable on }\left[ 0,T%
\right] \text{ for any }\varepsilon >0,\text{ and }\lim_{\varepsilon
\rightarrow 0^{+}}\int_{0}^{T}g_{\varepsilon }\left( s\right) \mathrm{d}s=0%
\Big \}.%
\end{array}%
\end{equation*}%
Denote $\mathbb{L}_{2}^{1,p}\left( \mathbb{R}^{d}\right) =\mathbb{L}%
_{2+}^{1,p}\left( \mathbb{R}^{d}\right) \cap \mathbb{L}_{2-}^{1,p}\left(
\mathbb{R}^{d}\right) .$ For any $\varphi \left( \cdot \right) \in \mathbb{L}%
_{2}^{1,p}\left( \mathbb{R}^{d}\right) ,$ denote $\nabla \varphi \left(
\cdot \right) =\mathcal{D}^{+}\varphi \left( \cdot \right) +\mathcal{D}%
^{-}\varphi \left( \cdot \right) .$ Whenever $\varphi $ is adapted, it
follows that $\mathcal{D}_{s}\varphi \left( t\right) =0$ for $t<s.$
Furthermore, $\nabla \varphi \left( \cdot \right) =\mathcal{D}^{+}\varphi
\left( \cdot \right) $ since $\mathcal{D}^{-}\varphi \left( \cdot \right)
=0. $ Put $\mathbb{L}_{2,\mathbb{F}}^{1,p}\left( \mathbb{R}^{d}\right) $ as
the set of all adapted processes in $\mathbb{L}_{2}^{1,p}\left( \mathbb{R}%
^{d}\right) .$

%We need the following lemma:
%
%\begin{lemma}
%\label{l2} Let $\varphi \in \mathbb{L}_{2,\mathbb{F}}^{1,p}\left( \mathbb{R}%
%^{d}\right) .$ Then, there exists a positive sequence $\left\{ \theta
%_{n}\right\} _{n=1}^{\infty }$ such that $\theta _{n}\rightarrow 0^{+}$ as $%
%n\rightarrow +\infty $ and
%\begin{equation*}
%\lim_{n\rightarrow \infty }\frac{1}{\theta _{n}^{2}}\int_{\tau }^{\tau
%+\theta }\int_{\tau }^{t}\mathbb{E}\left\vert \mathcal{D}_{s}\varphi \left(
%t\right) -\mathcal{D}\varphi \left( s\right) \right\vert ^{2}\mathrm{d}s%
%\mathrm{d}t=0,\text{ }a.e.\text{ }\tau \in \left[ 0,T\right] .
%\end{equation*}
%\end{lemma}

\section{Maximum Principle of Singular Optimal Controls}

\label{sec3}

This section will study the optimal controls separately. Due to the special
structure of control systems, we shall first consider the singular control
part, deriving the necessary condition, subsequently, regular part. The
initial condition will fixed to be $\left( 0,x\right) ,$ $x\in \mathbb{R}%
^{n}.$ At the beginning let us suppose that $\left( \bar{u}\left( \cdot
\right) ,\bar{\xi}\left( \cdot \right) \right) \in \mathcal{U}_{1}\times
\mathcal{U}_{2}$ is an optimal control and denote by $\left( X^{0,x;\bar{u},%
\bar{\eta}}\left( \cdot \right) ,Y^{0,x;\bar{u},\bar{\eta}}\left( \cdot
\right) ,Z^{0,x;\bar{u},\bar{\eta}}\left( \cdot \right) \right) $ the
optimal solution of (\ref{FBSDE1}). Our maximum principle will be proved in
two steps. The first variational inequality is derived from the fact
\begin{equation}
J\left( 0,x,u^{\varepsilon }\left( \cdot \right) ,\bar{\xi}\left( \cdot
\right) \right) -J\left( 0,x,\bar{u}\left( \cdot \right) ,\bar{\xi}\left(
\cdot \right) \right) \geq 0,  \label{3.1.1}
\end{equation}%
where $u^{\varepsilon }\left( \cdot \right) $ is a convex perturbation of
optimal control. The second variational inequity is attained from the
inequity
\begin{equation}
J\left( 0,x,\bar{u}\left( \cdot \right) ,\xi ^{\varepsilon }\left( \cdot
\right) \right) -J\left( 0,x,\bar{u}\left( \cdot \right) ,\bar{\xi}\left(
\cdot \right) \right) \geq 0,  \label{3.1.2}
\end{equation}%
where $\xi ^{\varepsilon }\left( \cdot \right) $ is a convex perturbation of
$\xi .$

\subsection{Optimal Singular Control}

\label{sec3.1}

For $l=b\left( \cdot \right) ,$ $\sigma \left( \cdot \right) ,$ $f\left(
\cdot \right) ,$ we denote%
\begin{eqnarray*}
\bar{l}_{x}\left( r,\cdot \right) &=&l_{x}\left( r,X^{0,x;\bar{u},\bar{\eta}%
}\left( r\right) ,Y^{0,x;\bar{u},\bar{\eta}}\left( r\right) ,Z^{0,x;\bar{u},%
\bar{\eta}}\left( r\right) ,\bar{u}\left( r\right) \right) , \\
\bar{l}_{y}\left( r,\cdot \right) &=&l_{y}\left( r,X^{0,x;\bar{u},\bar{\eta}%
}\left( r\right) ,Y^{0,x;\bar{u},\bar{\eta}}\left( r\right) ,Z^{0,x;\bar{u},%
\bar{\eta}}\left( r\right) ,\bar{u}\left( r\right) \right) .
\end{eqnarray*}%
Let us introduce the following

\begin{proposition}
Let \emph{(A1)-(A2)} hold, and let $\left( X^{0,x;\bar{u},\bar{\eta}}\left(
\cdot \right) ,Y^{0,x;\bar{u},\bar{\eta}}\left( \cdot \right) ,Z^{0,x;\bar{u}%
,\bar{\eta}}\left( \cdot \right) \right) \in \mathcal{N}^{2}(0,T;\mathbb{R}%
^{n})$ be an optimal solution. Then, the following FBSDEs:%
\begin{equation}
\left\{
\begin{array}{lll}
\mathrm{d}\mathfrak{p}\left( r\right) & = & -\left[ \bar{b}_{x}\left(
r,\cdot \right) ^{\top }\mathfrak{p}\left( r\right) +\bar{\sigma}_{x}\left(
r,\cdot \right) ^{\top }\mathfrak{k}\left( r\right) -\bar{f}_{x}\left(
r,\cdot \right) ^{\top }\mathfrak{q}\left( r\right) \right] \mathrm{d}r+%
\mathfrak{k}\left( r\right) \mathrm{d}W\left( r\right) , \\
\mathrm{d}\mathfrak{q}\left( r\right) & = & \bar{f}_{y}\left( r,\cdot
\right) ^{\top }\mathfrak{q}\left( r\right) \mathrm{d}r+\bar{f}_{z}\left(
r,\cdot \right) ^{\top }\mathfrak{q}\left( r\right) \mathrm{d}W\left(
r\right) , \\
\mathfrak{p}\left( T\right) & = & -\Phi _{x}\left( X^{t,x;u^{\varepsilon
}}\left( T\right) \right) \mathfrak{q}\left( T\right) ,\qquad \mathfrak{q}%
\left( 0\right) =1,%
\end{array}%
\right.  \label{firstadjoint}
\end{equation}%
admit an adapted solution $\left( \mathfrak{p}\left( \cdot \right) ,%
\mathfrak{q}\left( \cdot \right) ,\mathfrak{k}\left( \cdot \right) \right)
\in \mathcal{N}^{2}(0,T;\mathbb{R}^{n}).$
\end{proposition}

\begin{theorem}
\label{mpsing}Let \emph{(A1)-(A2)} hold. If $\left( X^{\bar{u},\bar{\eta}%
}\left( \cdot \right) ,Y^{\bar{u},\bar{\eta}}\left( \cdot \right) ,Z^{\bar{u}%
,\bar{\eta}}\left( \cdot \right) ,\bar{u}\left( \cdot \right) ,\bar{\xi}%
\left( \cdot \right) \right) $ is an optimal solution of (\ref{FBSDE1}),
then there exists a unique pair of adapted processes $\left( \mathfrak{p}%
\left( \cdot \right) ,\mathfrak{q}\left( \cdot \right) \right) $ satisfying (%
\ref{firstadjoint}) such that%
\begin{equation}
P\left\{ \mathfrak{q}\left( t\right) K_{\left( i\right) }-\mathfrak{p}^{\top
}\left( t\right) G_{\left( i\right) }\left( t\right) \geq 0,\text{ }t\in %
\left[ 0,T\right] ,\text{ }\forall i\right\} =1,  \label{3.1.7}
\end{equation}%
and%
\begin{equation}
P\left\{ \sum_{i=1}^{m}\int_{0}^{T}I_{\left\{ \mathfrak{q}\left( r\right)
K_{\left( i\right) }-\mathfrak{p}^{T}\left( r\right) G_{\left( i\right)
}\left( t\right) >0\right\} }\mathrm{d}\bar{\xi}_{r}^{\left( i\right)
}=0\right\} =1.  \label{3.1.8}
\end{equation}
\end{theorem}

Before the proof, we need some lemmas. At the beginning, we introduce the
convex perturbation%
\begin{equation*}
\left( \bar{u}\left( t\right) ,\xi ^{\alpha }\left( t\right) \right) =\left(
\bar{u}\left( t\right) ,\bar{\xi}\left( t\right) +\alpha \left( \xi \left(
t\right) -\bar{\xi}\left( t\right) \right) \right) ,
\end{equation*}%
where $\alpha \in \left[ 0,1\right] $ and $\xi \left( \cdot \right) $ is an
arbitrary element of $\mathcal{U}_{2}.$ We now introduce the following
variational equations of (\ref{FBSDE1}):%
\begin{equation}
\left\{
\begin{array}{lll}
\mathrm{d}x^{1}\left( t\right) & = & \bar{b}_{x}\left( t\right) x^{1}\left(
t\right) \text{d}t+\bar{\sigma}_{x}\left( t\right) x^{1}\left( t\right)
\text{d}W\left( t\right) +G\left( t\right) \mathrm{d}\left( \xi \left(
t\right) -\bar{\xi}\left( t\right) \right) , \\
\mathrm{d}y^{1}\left( t\right) & = & -\bar{f}_{x}\left( t\right) x^{1}\left(
t\right) -\bar{f}_{y}\left( t\right) y^{1}\left( t\right) -\bar{f}_{z}\left(
t\right) z^{1}\left( t\right) \text{d}t+z^{1}\left( t\right) \text{d}W\left(
t\right) \\
&  & -K\mathrm{d}\left( \xi \left( t\right) -\bar{\xi}\left( t\right)
\right) , \\
y^{1}\left( T\right) & = & \Phi _{x}\left( X^{0,x;\bar{u},\bar{\xi}}\left(
T\right) \right) x^{1}\left( T\right) ,\text{ }x^{1}\left( 0\right) =0.%
\end{array}%
\right.  \label{3.1.9}
\end{equation}%
From (A1)-(A2) it is easy to check that (\ref{3.1.9}) has a unique strong
solution. Moreover, we have

\begin{lemma}
\label{estvs}Under the Assumptions \emph{(A1)-(A2)}, we have%
\begin{eqnarray}
\lim\limits_{\alpha \rightarrow 0}\mathbb{E}\left[ \left\vert \frac{X^{0,x;%
\bar{u},\xi ^{\alpha }}\left( t\right) -X^{0,x;\bar{u},\bar{\xi}}\left(
t\right) }{\alpha }-x^{1}\left( t\right) \right\vert ^{2}\right] &=&0,\text{
}t\in \left[ 0,T\right] ,  \label{e1} \\
\lim\limits_{\alpha \rightarrow 0}\mathbb{E}\left[ \left\vert \frac{Y^{0,x;%
\bar{u},\xi ^{\alpha }}\left( t\right) -Y^{0,x;\bar{u},\bar{\xi}}\left(
t\right) }{\alpha }-y^{1}\left( t\right) \right\vert ^{2}\right] &=&0,\text{
}t\in \left[ 0,T\right] ,  \label{e2} \\
\lim\limits_{\alpha \rightarrow 0}\mathbb{E}\left[ \int_{0}^{T}\left\vert
\frac{Z^{0,x;\bar{u},\xi ^{\alpha }}\left( t\right) -Z^{0,x;\bar{u},\bar{\xi}%
}\left( t\right) }{\alpha }-z^{1}\left( t\right) \right\vert ^{2}\mathrm{d}t%
\right] &=&0.  \label{e3}
\end{eqnarray}
\end{lemma}

The proof can be seen in the Appendix.

\paragraph{Proof of Theorem \protect\ref{mpsing}.}

Applying It\^{o}'s formula to $\left\langle \mathfrak{p}\left( \cdot \right)
,x^{1}\left( \cdot \right) \right\rangle +\mathfrak{q}\left( \cdot \right)
y^{1}\left( \cdot \right) $ on $\left[ 0,T\right] $ yields%
\begin{equation}
0\leq y^{1}\left( 0\right) =\mathbb{E}\left[ \int_{0}^{T}\left( \mathfrak{q}%
\left( r\right) K-\mathfrak{p}^{\top }\left( r\right) G\left( t\right)
\right) \mathrm{d}\left( \xi \left( t\right) -\bar{\xi}\left( t\right)
\right) \text{d}t\right] .  \label{singularinq}
\end{equation}%
In particular, let $\xi \in \mathcal{U}_{2}$ be a process satisfying $%
P\left\{ \sum_{i}\int_{0}^{T}G\left( s\right) \mathrm{d}\xi _{s}^{\left(
i\right) }<\infty \right\} $ and such that (\ref{singularinq}) and
\begin{equation*}
\mathrm{d}\xi _{s}^{\left( i\right) }=\left\{
\begin{array}{ll}
0 & \text{if }\mathfrak{q}\left( r\right) K_{\left( i\right) }-\mathfrak{p}%
^{\top }\left( r\right) G_{\left( i\right) }\left( t\right) >0, \\
\hat{\xi}_{s}^{\left( i\right) } & \text{otherwise,}%
\end{array}%
\right.
\end{equation*}%
holds where $\xi _{s}^{\left( i\right) }$ denotes the $i$th component. Then,
\begin{eqnarray*}
&&\mathbb{E}\left[ \sum_{i=1}^{m}\int_{0}^{T}\left( \mathfrak{q}\left(
r\right) K_{\left( i\right) }-\mathfrak{p}^{\top }\left( r\right) G_{\left(
i\right) }\left( t\right) \right) I_{\left\{ \mathfrak{q}\left( r\right)
K_{\left( i\right) }-\mathfrak{p}^{T}\left( r\right) G_{\left( i\right)
}\left( t\right) >0\right\} }\mathrm{d}\left( -\bar{\xi}^{\left( i\right)
}\left( t\right) \right) \right] \\
&=&\mathbb{E}\left[ \int_{0}^{T}\left( \mathfrak{q}\left( r\right) K-%
\mathfrak{p}^{\top }\left( r\right) G\left( t\right) \right) \mathrm{d}%
\left( \xi \left( t\right) -\bar{\xi}\left( t\right) \right) \text{d}t\right]
\\
&\geq &0.
\end{eqnarray*}%
Thus%
\begin{equation*}
\mathbb{E}\left[ \sum_{i=1}^{m}\int_{0}^{T}\left( \mathfrak{q}\left(
r\right) K_{\left( i\right) }-\mathfrak{p}^{\top }\left( r\right) G_{\left(
i\right) }\left( t\right) \right) I_{\left\{ \mathfrak{q}\left( r\right)
K_{\left( i\right) }-\mathfrak{p}^{T}\left( r\right) G_{\left( i\right)
}\left( t\right) >0\right\} }\mathrm{d}\left( -\bar{\xi}^{\left( i\right)
}\left( t\right) \right) \text{d}t\right] =0
\end{equation*}%
which proves (\ref{3.1.8}). Next we show that (\ref{3.1.7}) is valid. For
that, let us define the events:
\begin{equation*}
\mathcal{A}^{\left( i\right) }\triangleq \{\left( t,\omega \right) \in \left[
0,T\right] \times \Omega :\mathfrak{q}\left( r\right) K_{\left( i\right) }-%
\mathfrak{p}^{\top }\left( r\right) G_{\left( i\right) }\left( t\right) <0\},
\end{equation*}%
where $t\in \left[ 0,T\right] ,$ $1\leq i\leq m.$

Define the stochastic
process $\breve{\xi}^{\left( i\right) }:\left[ 0,T\right] \times \Omega
\rightarrow \left[ 0,\infty \right) $ by
\begin{equation*}
\breve{\xi}^{\left( i\right) }\left( t\right) =\bar{\xi}^{\left( i\right)
}\left( t\right) +\int_{0}^{t}I_{\mathcal{A}^{\left( i\right) }}\left(
r,\omega \right) \mathrm{d}r.
\end{equation*}%
Then one can easily check that $\breve{\xi}=\left( \breve{\xi}^{\left(
1\right) },\breve{\xi}^{\left( 2\right) },\ldots ,\breve{\xi}^{\left(
m\right) }\right) $ is a measurable, adapted process which is nondecreasing
left-continuous with right limits and $\breve{\xi}\left( 0\right) =0$, and
which satisfies
\begin{equation*}
P\left\{ \sum_{i}\int_{0}^{T}G\left( s\right) \mathrm{d}\breve{\xi}%
_{s}^{\left( i\right) }<\infty \right\} .
\end{equation*}%
Further, we have
\begin{eqnarray*}
&&\mathbb{E}\left[ \int_{0}^{T}\left( \mathfrak{q}\left( r\right) K-%
\mathfrak{p}^{\top }\left( r\right) G\left( t\right) \right) \mathrm{d}%
\left( \breve{\xi}\left( t\right) -\bar{\xi}\left( t\right) \right) \text{d}t%
\right] \\
&=&\mathbb{E}\left[ \sum_{i=1}^{m}\int_{0}^{T}\left( \mathfrak{q}\left(
r\right) K_{\left( i\right) }-\mathfrak{p}^{\top }\left( r\right) G_{\left(
i\right) }\left( t\right) \right) I_{\mathcal{A}^{\left( i\right) }}\text{d}t%
\right] \\
&<&0,
\end{eqnarray*}%
which obviously contradicts to (\ref{singularinq}), unless for any $i,$ we
have $\left( Leb\otimes P\right) \left\{ \mathcal{A}^{\left( i\right)
}\right\} =0.$ We thus complete the proof. \hfill $\Box $

\begin{remark}
One can easily check that%
\begin{equation*}
\mathfrak{q}\left( s\right) =\exp \Bigg \{\int_{t}^{s}\bigg [\bar{f}%
_{y}\left( r\right) -\frac{1}{2}\left\vert \bar{f}_{z}\left( r\right)
\right\vert ^{2}\bigg ]\mathrm{d}r+\int_{0}^{s}\bar{f}_{z}\left( r\right) %
\mbox{\rm d}W\left( r\right) \Bigg \},
\end{equation*}%
which implies that\textsl{\ }$\mathfrak{q}\left( r\right) >0,$ $r\in \left[
0,T\right] $, $P$-a.s. So $-\mathfrak{p}\left( \cdot \right) /\mathfrak{q}%
\left( \cdot \right) $ makes sense. Clearly, our Theorem \ref{mpsing} for
optimal singular control is completely different from \cite{BDM1}. Ours
contains two variables $\left( \mathfrak{p}\left( \cdot \right) ,\mathfrak{q}%
\left( \cdot \right) \right) $. As a matter of fact, we have
\begin{equation*}
P\left\{ K_{\left( i\right) }+\left( -\frac{\mathfrak{p}^{\top }\left(
t\right) }{\mathfrak{q}\left( t\right) }\right) G_{\left( i\right) }\left(
t\right) \geq 0,\text{ }t\in \left[ 0,T\right] ,\text{ }\forall i\right\} =1.
\end{equation*}%
We claim that $-\mathfrak{p}\left( \cdot \right) /\mathfrak{q}\left( \cdot
\right) $ is the partial derivative of value function, which will be studied
in Section \ref{secconne}.
\end{remark}

\subsection{Optimal Regular Control}

\label{sec3.2}

In this subsection, we study the optimal regular controls for systems driven
by FBSDEs (\ref{FBSDE1}) under the types of Pontryagin, namely, necessary
maximum principles for optimal control. To this end, we fix $\bar{\xi}\in
\mathcal{U}_{2}$ and introduce the following convex perturbation control.
Taking $u\left( \cdot \right) \in \mathcal{U}_{1},$ we define $v\left( \cdot
\right) =u\left( \cdot \right) -\bar{u}\left( \cdot \right) ,$ $%
u^{\varepsilon }\left( \cdot \right) =\bar{u}\left( \cdot \right)
+\varepsilon v\left( \cdot \right) ,$ where $\varepsilon >0$ is sufficiently
small$.$ Since $U$ is convex, $u^{\varepsilon }\left( \cdot \right) \in
\mathcal{U}\left( 0,T\right) .$ Let $\left( x^{\varepsilon },y^{\varepsilon
},z^{\varepsilon },u^{\varepsilon }\right) $ be the trajectory of the
control system (\ref{FBSDE1}) corresponding to the control $u^{\varepsilon
}. $ Put $\delta x\left( \cdot \right) =x^{\varepsilon }\left( \cdot \right)
-\bar{x}\left( \cdot \right) .$ When $l=b,$ $\sigma $ and $\Phi ,$ we denote
\begin{eqnarray*}
l_{x}\left( t\right) &=&l_{x}\left( t,\bar{x}\left( t\right) ,\bar{u}\left(
t\right) \right) ,\text{ }l_{u}\left( t\right) =l_{u}\left( t,\bar{x}\left(
t\right) ,\bar{u}\left( t\right) \right) , \\
l_{xx}\left( t\right) &=&l_{xx}\left( t,\bar{x}\left( t\right) ,\bar{u}%
\left( t\right) \right) ,\text{ }l_{xu}\left( t\right) =l_{xu}\left( t,\bar{x%
}\left( t\right) ,\bar{u}\left( t\right) \right) , \\
l_{uu}\left( t\right) &=&l_{uu}\left( t,\bar{x}\left( t\right) ,\bar{u}%
\left( t\right) \right) ,\text{ }f_{x}\left( t\right) =f_{x}\left( t,\bar{x}%
\left( t\right) ,\bar{y}\left( t\right) ,\bar{z}\left( t\right) ,\bar{u}%
\left( t\right) \right) , \\
f_{y}\left( t\right) &=&f_{y}\left( t,\bar{x}\left( t\right) ,\bar{y}\left(
t\right) ,\bar{z}\left( t\right) ,\bar{u}\left( t\right) \right) ,\text{ }%
f_{z}\left( t\right) =f_{z}\left( t,\bar{x}\left( t\right) ,\bar{y}\left(
t\right) ,\bar{z}\left( t\right) ,\bar{u}\left( t\right) \right) , \\
f_{xx}\left( t\right) &=&f_{xx}\left( t,\bar{x}\left( t\right) ,\bar{y}%
\left( t\right) ,\bar{z}\left( t\right) ,\bar{u}\left( t\right) \right) ,%
\text{ }f_{xu}\left( t\right) =f_{xu}\left( t,\bar{x}\left( t\right) ,\bar{y}%
\left( t\right) ,\bar{z}\left( t\right) ,\bar{u}\left( t\right) \right) , \\
f_{yy}\left( t\right) &=&f_{yy}\left( t,\bar{x}\left( t\right) ,\bar{y}%
\left( t\right) ,\bar{z}\left( t\right) ,\bar{u}\left( t\right) \right)
,\cdots \text{\textit{etc.}}
\end{eqnarray*}%
Let us introduce the following two kinds of variational equations, mainly
taken from \cite{BS}. For simplicity, we omit the superscript.
\begin{equation}
\left\{
\begin{array}{lll}
\mathrm{d}x_{1}\left( t\right) & = & \left[ b_{x}\left( t\right) x_{1}\left(
t\right) +b_{u}\left( t\right) v\left( t\right) \right] \mathrm{d}t+\left[
\sigma _{x}\left( t\right) x_{1}\left( t\right) +\sigma _{u}\left( t\right)
v\left( t\right) \right] \mathrm{d}W\left( t\right) , \\
x_{1}\left( 0\right) & = & 0,\qquad t\in \left[ 0,T\right] ,%
\end{array}%
\right.  \label{varsde1}
\end{equation}%
and
\begin{equation}
\left\{
\begin{array}{lll}
\mathrm{d}x_{2}\left( t\right) & = & \Big [b_{x}\left( t\right) x_{2}\left(
t\right) +x_{1}^{\top }\left( t\right) b_{xx}\left( t\right) x_{1}\left(
t\right) \\
&  & +2v\left( t\right) ^{\top }b_{xu}\left( t\right) x_{1}\left( t\right)
+v^{\top }\left( t\right) b_{uu}\left( t\right) v\left( t\right) \Big ]%
\mathrm{d}t, \\
&  & \Big [\sigma _{x}\left( t\right) x_{2}\left( t\right) +x_{1}^{\top
}\left( t\right) \sigma _{xx}\left( t\right) x_{1}\left( t\right) \\
&  & +2v\left( t\right) ^{\top }\sigma _{xu}\left( t\right) x_{1}\left(
t\right) +v^{\top }\left( t\right) \sigma _{uu}\left( t\right) v\left(
t\right) \Big ]\mathrm{d}W\left( t\right) , \\
x_{2}\left( 0\right) & = & 0,\qquad t\in \left[ 0,T\right] .%
\end{array}%
\right.  \label{varsde2}
\end{equation}

From Lemma 3.5 and Lemma 3.11 in \cite{BS}, we have following result.

\begin{lemma}
\label{l3}Assume that \emph{(A1)-(A2)} is in force. Then, we have, for any $%
\beta \geq 2,$%
\begin{eqnarray*}
\left\Vert x_{1}\right\Vert _{\infty }^{\beta } &\leq &C,\text{ }\left\Vert
x_{2}\right\Vert _{\infty }^{\beta }\leq C,\text{ }\left\Vert \delta
x\right\Vert _{\infty }^{\beta }\leq C\varepsilon ^{\beta },\text{ } \\
\left\Vert \delta x-\varepsilon x_{1}\right\Vert _{\infty }^{\beta } &\leq
&C\varepsilon ^{2\beta },\text{ }\left\Vert \delta x-\varepsilon x_{1}-\frac{%
\varepsilon ^{2}}{2}x_{2}\right\Vert _{\infty }^{\beta }\leq C\varepsilon
^{3\beta },
\end{eqnarray*}%
where
\begin{equation*}
\left\Vert x_{1}\right\Vert _{\infty }^{\beta }=\left[ \mathbb{E}\left(
\sup_{t\in \left[ 0,T\right] }\left\vert x_{1}\left( t\right) \right\vert
^{\beta }\right) \right] ^{\frac{1}{\beta }}.
\end{equation*}
\end{lemma}

We shall introduce the so called variational equations for FBSDEs (\ref%
{FBSDE1}) beginning from the following two adjoint equations:%
\begin{equation}
\left\{
\begin{array}{rcl}
-\mathrm{d}p\left( t\right) & = & \Gamma \left( t\right) \mathrm{d}t-q\left(
t\right) \mathrm{d}W\left( t\right) , \\
p\left( T\right) & = & \Phi _{x}\left( \bar{x}\left( T\right) \right) ,%
\end{array}%
\right.  \label{adj1}
\end{equation}%
and%
\begin{equation}
\left\{
\begin{array}{rcl}
-\mathrm{d}P\left( t\right) & = & \Pi \left( t\right) \mathrm{d}t-Q\left(
t\right) \mathrm{d}W\left( t\right) , \\
P\left( T\right) & = & \Phi _{xx}\left( \bar{x}\left( T\right) \right) ,%
\end{array}%
\right.  \label{adj2}
\end{equation}%
where $\Gamma \left( \cdot \right) ,$ $\Pi \left( \cdot \right) $ are
unknown two processes to be determined. Next we will derive two kinds of
adjoint equations. The main idea is borrowed from \cite{Hms}. First of all,
we observe that
\begin{eqnarray*}
\Phi \left( x^{\varepsilon }\left( T\right) \right) -\Phi \left( \bar{x}%
\left( T\right) \right) &=&\left\langle \delta x\left( T\right) ,\Phi
_{x}\left( \bar{x}\left( T\right) \right) \right\rangle \\
&&+\frac{1}{2}\left\langle \Phi _{xx}\left( \bar{x}\left( T\right) \right)
\delta x\left( T\right) ,\delta x\left( T\right) \right\rangle +o\left(
\varepsilon ^{2}\right) \\
&=&\left\langle \varepsilon x_{1}\left( T\right) +\frac{1}{2}\varepsilon
^{2}x_{2}\left( T\right) ,\Phi _{x}\left( \bar{x}\left( T\right) \right)
\right\rangle \\
&&+\frac{1}{2}\varepsilon ^{2}\left\langle \Phi _{xx}\left( \bar{x}\left(
T\right) \right) x_{1}\left( T\right) ,x_{1}\left( T\right) \right\rangle
+o\left( \varepsilon ^{2}\right) ,
\end{eqnarray*}%
which inspires us to use the adjoint equations to expand the following:
\begin{equation}
\left\langle \varepsilon x_{1}\left( t\right) +\frac{1}{2}\varepsilon
^{2}x_{2}\left( t\right) ,p\left( t\right) \right\rangle +\frac{1}{2}%
\varepsilon ^{2}\text{tr}\left[ P\left( t\right) x_{1}\left( t\right)
x_{1}^{\top }\left( t\right) \right] \text{ on }\left[ 0,T\right] .
\label{ito}
\end{equation}%
It\^{o}'s formula applied to (\ref{ito}) yields for $t\in \left[ 0,T\right]
, $%
\begin{eqnarray*}
&&\left\langle \varepsilon x_{1}\left( T\right) +\frac{1}{2}\varepsilon
^{2}x_{2}\left( T\right) ,p\left( T\right) \right\rangle +\frac{1}{2}%
\varepsilon ^{2}\text{tr}\left[ P\left( T\right) x_{1}\left( T\right)
x_{1}^{\top }\left( T\right) \right] \\
&=&\left\langle \varepsilon x_{1}\left( t\right) +\frac{1}{2}\varepsilon
^{2}x_{2}\left( t\right) ,p\left( t\right) \right\rangle +\frac{1}{2}%
\varepsilon ^{2}\text{tr}\left[ P\left( t\right) x_{1}\left( t\right)
x_{1}^{\top }\left( t\right) \right] \\
&&+\int_{t}^{T}\Bigg [\left\langle \left( \varepsilon x_{1}\left( s\right) +%
\frac{1}{2}\varepsilon ^{2}x_{2}\left( s\right) \right) ,\Lambda _{1}\left(
s\right) \right\rangle +\varepsilon ^{2}\left\langle x_{1}\left( s\right)
,\Lambda _{2}\left( s\right) \right\rangle \\
&&+\frac{1}{2}\varepsilon ^{2}\left\langle \Lambda _{3}\left( s\right)
x_{1}\left( s\right) ,x_{1}\left( s\right) \right\rangle +\Lambda _{4}\left(
s\right) \Bigg ]\mathrm{d}s \\
&&+\int_{t}^{T}\Bigg [\left\langle \varepsilon x_{1}\left( s\right) +\frac{1%
}{2}\varepsilon ^{2}x_{2}\left( s\right) ,\Lambda _{5}\left( s\right)
\right\rangle +\varepsilon ^{2}\left\langle x_{1}\left( s\right) ,\Lambda
_{6}\left( s\right) \right\rangle \\
&&+\frac{1}{2}\varepsilon ^{2}\left\langle \Lambda _{7}\left( s\right)
x_{1}\left( s\right) ,x_{1}\left( s\right) \right\rangle +\Lambda _{8}\left(
s\right) \Bigg ]\mathrm{d}W\left( s\right) ,
\end{eqnarray*}%
where

\begin{eqnarray*}
\Lambda _{1}\left( t\right) &=&b_{x}^{\top }\left( t\right) p\left( t\right)
+\sigma _{x}^{\top }\left( t\right) q\left( t\right) -\Gamma \left( t\right)
, \\
\Lambda _{2}\left( t\right) &=&b_{xu}^{\top }\left( t\right) v\left(
t\right) p\left( t\right) +\sigma _{xu}^{\top }\left( t\right) v\left(
t\right) q\left( t\right) +P\left( t\right) b_{u}\left( t\right) v\left(
t\right) \\
&&+\sigma _{x}^{\top }\left( t\right) P\left( t\right) \sigma _{u}\left(
t\right) v\left( t\right) +Q\left( t\right) \sigma _{u}\left( t\right)
v\left( t\right) , \\
\Lambda _{3}\left( t\right) &=&p\left( t\right) b_{xx}\left( t\right)
+q\left( t\right) \sigma _{xx}\left( t\right) +P\left( t\right) b_{x}\left(
t\right) +b_{x}^{\top }\left( t\right) P\left( t\right) \\
&&+\sigma _{x}^{\top }\left( t\right) P\left( t\right) \sigma _{x}\left(
t\right) +Q\left( t\right) \sigma _{x}\left( t\right) +\sigma _{x}^{\top
}\left( t\right) Q\left( t\right) -\Pi \left( t\right) , \\
\Lambda _{4}\left( t\right) &=&\varepsilon \left\langle p\left( t\right)
b_{u}\left( t\right) ,v\left( t\right) \right\rangle +\frac{1}{2}\varepsilon
^{2}\left\langle p\left( t\right) b_{uu}\left( t\right) v\left( t\right)
,v\left( t\right) \right\rangle \\
&&+\varepsilon \left\langle q\left( t\right) \sigma _{u}\left( t\right)
,v\left( t\right) \right\rangle +\frac{1}{2}\varepsilon ^{2}\left\langle
q\left( t\right) \sigma _{uu}\left( t\right) v\left( t\right) ,v\left(
t\right) \right\rangle \\
&&+\frac{1}{2}\varepsilon ^{2}\left\langle \sigma _{u}\left( t\right)
P\left( t\right) \sigma _{u}^{\top }v\left( t\right) ,v\left( t\right)
\right\rangle , \\
\Lambda _{5}\left( t\right) &=&p\left( t\right) \sigma _{x}\left( t\right)
+q\left( t\right) , \\
\Lambda _{6}\left( t\right) &=&p\left( t\right) v\left( t\right) \sigma
_{xu}\left( t\right) +P\left( t\right) \sigma _{u}\left( t\right) v\left(
t\right) , \\
\Lambda _{7}\left( t\right) &=&p\left( t\right) \sigma _{xx}\left( t\right)
+P\left( t\right) \sigma _{x}\left( t\right) +\sigma _{x}^{\top }\left(
t\right) P\left( t\right) +Q\left( t\right) , \\
\Lambda _{8}\left( t\right) &=&\varepsilon \left\langle \sigma _{u}\left(
t\right) p\left( t\right) ,v\left( t\right) \right\rangle +\frac{1}{2}%
\varepsilon ^{2}\left\langle p\left( t\right) \sigma _{uu}\left( t\right)
v\left( t\right) ,v\left( t\right) \right\rangle .
\end{eqnarray*}

\begin{remark}
Note that $\Gamma \left( t\right) $ and $\Pi \left( t\right) $ do not appear
in the $\mathrm{d}W\left( s\right) $-term.
\end{remark}

Define
\begin{equation*}
\left\{
\begin{array}{rcl}
\mathrm{d}y^{\varepsilon }(t) & = & -f\left( (t,x^{\varepsilon }\left(
t\right) ,y^{\varepsilon }\left( t\right) ,z^{\varepsilon }\left( t\right)
,u^{\varepsilon }\left( t\right) \right) \mathrm{d}t+z^{\varepsilon }\left(
t\right) \mathrm{d}W\left( t\right) , \\
y^{\varepsilon }\left( T\right) & = & \Phi \left( x^{\varepsilon }\left(
T\right) \right) .%
\end{array}%
\right.
\end{equation*}%
Let
\begin{eqnarray}
\bar{y}^{\varepsilon }\left( t\right) &=&y^{\varepsilon }(t)-\left[
\left\langle p\left( t\right) ,\left( \varepsilon x_{1}\left( t\right) +%
\frac{1}{2}\varepsilon ^{2}x_{2}\left( t\right) \right) \right\rangle +\frac{%
1}{2}\varepsilon ^{2}\left\langle P\left( t\right) x_{1}\left( t\right)
,x_{1}\left( t\right) \right\rangle \right]  \label{exp1} \\
\bar{z}^{\varepsilon }\left( t\right) &=&z^{\varepsilon }(t)-\Bigg [%
\left\langle \left( \varepsilon x_{1}\left( s\right) +\frac{1}{2}\varepsilon
^{2}x_{2}\left( s\right) \right) ,\Lambda _{5}\left( s\right) \right\rangle
+\varepsilon ^{2}\left\langle x_{1}\left( s\right) ,\Lambda _{6}\left(
s\right) \right\rangle  \notag \\
&&+\frac{1}{2}\varepsilon ^{2}\left\langle \Lambda _{7}\left( s\right)
x_{1}\left( s\right) ,x_{1}\left( s\right) \right\rangle +\Lambda _{8}\left(
s\right) \Bigg ].  \label{exp2}
\end{eqnarray}%
Clearly, from Lemma \ref{l3}, we have
\begin{eqnarray*}
\Phi \left( x^{\varepsilon }\left( T\right) \right) &=&\Phi \left( \bar{x}%
\left( T\right) \right) +\left\langle p\left( T\right) ,\left( \varepsilon
x_{1}\left( T\right) +\frac{1}{2}\varepsilon ^{2}x_{2}\left( T\right)
\right) \right\rangle \\
&&+\frac{1}{2}\varepsilon ^{2}\left\langle P\left( t\right) x_{1}\left(
t\right) ,x_{1}\left( t\right) \right\rangle +o\left( \varepsilon
^{2}\right) .
\end{eqnarray*}%
After some tedious computations, we have%
\begin{eqnarray}
\bar{y}^{\varepsilon }\left( t\right) &=&\Phi \left( \bar{x}\left( T\right)
\right) +\int_{t}^{T}\Bigg [f\left( s,x^{\varepsilon }\left( s\right)
,y^{\varepsilon }\left( s\right) ,z^{\varepsilon }\left( s\right)
,u^{\varepsilon }\left( s\right) \right)  \notag \\
&&+\left\langle \left( \varepsilon x_{1}\left( s\right) +\frac{1}{2}%
\varepsilon ^{2}x_{2}\left( s\right) \right) ,\Lambda _{1}\left( s\right)
\right\rangle +\varepsilon ^{2}\left\langle x_{1}\left( s\right) ,\Lambda
_{2}\left( s\right) \right\rangle  \notag \\
&&+\frac{1}{2}\varepsilon ^{2}\left\langle \Lambda _{3}\left( s\right)
x_{1}\left( s\right) ,x_{1}\left( s\right) \right\rangle +\Lambda _{4}\left(
s\right) \Bigg ]\mathrm{d}s-\int_{t}^{T}\bar{z}^{\varepsilon }\left(
s\right) \mathrm{d}W\left( s\right) +o\left( \varepsilon ^{2}\right) .
\label{var1}
\end{eqnarray}%
Put
\begin{equation*}
\hat{y}^{\varepsilon }\left( t\right) =\bar{y}^{\varepsilon }(t)-\bar{y}%
\left( t\right) ,\text{ }\hat{z}^{\varepsilon }\left( t\right) =\bar{z}%
^{\varepsilon }(t)-\bar{z}\left( t\right) ,
\end{equation*}%
then we attain
\begin{eqnarray}
\hat{y}^{\varepsilon }\left( t\right) &=&o\left( \varepsilon ^{2}\right)
+\int_{t}^{T}\Bigg [f\left( s,x^{\varepsilon }\left( s\right)
,y^{\varepsilon }\left( s\right) ,z^{\varepsilon }\left( s\right)
,u^{\varepsilon }\left( s\right) \right)  \notag \\
&&-f\left( s,\bar{x}\left( s\right) ,\bar{y}\left( s\right) ,\bar{z}\left(
s\right) ,\bar{u}\left( s\right) \right) \\
&&+\left\langle \left( \varepsilon x_{1}\left( s\right) +\frac{1}{2}%
\varepsilon ^{2}x_{2}\left( s\right) \right) ,\Lambda _{1}\left( s\right)
\right\rangle +\varepsilon ^{2}\left\langle x_{1}\left( s\right) ,\Lambda
_{2}\left( s\right) \right\rangle  \notag \\
&&+\frac{1}{2}\varepsilon ^{2}\left\langle \Lambda _{3}\left( s\right)
x_{1}\left( s\right) ,x_{1}\left( s\right) \right\rangle +\Lambda _{4}\left(
s\right) \Bigg ]\mathrm{d}s-\int_{t}^{T}\hat{z}^{\varepsilon }\left(
s\right) \mathrm{d}W\left( s\right) ,  \label{var2}
\end{eqnarray}%
where%
\begin{eqnarray*}
&&f\left( s,x^{\varepsilon }\left( s\right) ,y^{\varepsilon }\left( s\right)
,z^{\varepsilon }\left( s\right) ,u^{\varepsilon }\left( s\right) \right)
-f\left( s,\bar{x}\left( s\right) ,\bar{y}\left( s\right) ,\bar{z}\left(
s\right) ,\bar{u}\left( s\right) \right) \\
&=&f\left( s,\bar{x}\left( s\right) ,\bar{y}\left( s\right) ,\bar{z}\left(
s\right) +\varepsilon p\left( t\right) \sigma _{u}\left( t\right) v\left(
t\right) +\frac{\varepsilon ^{2}}{2}p\left( t\right) \sigma _{uu}\left(
t\right) v^{2}\left( t\right) ,\bar{u}\left( s\right) \right) \\
&&-f\left( s,\bar{x}\left( s\right) ,\bar{y}\left( s\right) ,\bar{z}\left(
s\right) ,\bar{u}\left( s\right) \right) +f\left( s,x^{\varepsilon }\left(
s\right) ,y^{\varepsilon }\left( s\right) ,z^{\varepsilon }\left( s\right)
,u^{\varepsilon }\left( s\right) \right) \\
&&-f\left( s,\bar{x}\left( s\right) ,\bar{y}\left( s\right) ,\bar{z}\left(
s\right) +\varepsilon p\left( t\right) \sigma _{u}\left( t\right) v\left(
t\right) +\frac{\varepsilon ^{2}}{2}p\left( t\right) \sigma _{uu}\left(
t\right) v^{2}\left( t\right) ,\bar{u}\left( s\right) \right) \\
&=&f\left( s,\bar{x}\left( s\right) ,\bar{y}\left( s\right) ,\bar{z}\left(
s\right) +\varepsilon p\left( t\right) \sigma _{u}\left( t\right) v\left(
t\right) +\frac{\varepsilon ^{2}}{2}p\left( t\right) \sigma _{uu}\left(
t\right) v^{2}\left( t\right) ,\bar{u}\left( s\right) \right) \\
&&-f\left( s,\bar{x}\left( s\right) ,\bar{y}\left( s\right) ,\bar{z}\left(
s\right) ,\bar{u}\left( s\right) \right) \\
&&+f\Big (s,\bar{x}\left( s\right) +\varepsilon x_{1}\left( s\right) +\frac{%
\varepsilon ^{2}}{2}x_{2}\left( s\right) ,\bar{y}\left( s\right) +\hat{y}%
^{\varepsilon }\left( t\right) +\Lambda _{9}\left( s\right) , \\
&&\bar{z}\left( s\right) +\hat{z}^{\varepsilon }\left( t\right) +\Lambda
_{10}\left( s\right) ,u^{\varepsilon }\left( s\right) \Big )-f\left( s,\bar{x%
}\left( s\right) ,\bar{y}\left( s\right) ,\bar{z}\left( s\right) ,\bar{u}%
\left( s\right) \right) +o\left( \varepsilon ^{2}\right) .
\end{eqnarray*}%
Next we are going to seek $\Gamma \left( \cdot \right) ,$ $\Pi \left( \cdot
\right) ,$ determined by the optimal quadruple $(\bar{x}\left( \cdot \right)
,\bar{y}\left( \cdot \right) ,\bar{z}\left( \cdot \right) ,\bar{u}\left(
\cdot \right) ),$ such that
\begin{eqnarray*}
&&f\left( s,\bar{x}\left( s\right) +\varepsilon x_{1}\left( s\right) +\frac{1%
}{2}\varepsilon ^{2}x_{2}\left( s\right) ,\bar{y}\left( s\right) +\Lambda
_{9}\left( s\right) ,\bar{z}\left( s\right) +\Lambda _{10}\left( s\right)
,u^{\varepsilon }\left( s\right) \right) \\
&&-f\left( s,\bar{x}\left( s\right) ,\bar{y}\left( s\right) ,\bar{z}\left(
s\right) ,\bar{u}\left( s\right) \right) +\left\langle \left( \varepsilon
x_{1}\left( s\right) +\frac{1}{2}\varepsilon ^{2}x_{2}\left( s\right)
\right) ,\Lambda _{1}\left( s\right) \right\rangle \\
&&+\frac{1}{2}\varepsilon ^{2}\left\langle \Lambda _{3}\left( s\right)
x_{1}\left( s\right) ,x_{1}\left( s\right) \right\rangle =o\left(
\varepsilon ^{2}\right) ,
\end{eqnarray*}%
where
\begin{eqnarray*}
\Lambda _{9}\left( s\right) &=&\left\langle \left( \varepsilon x_{1}\left(
s\right) +\frac{1}{2}\varepsilon ^{2}x_{2}\left( s\right) \right) ,p\left(
s\right) \right\rangle +\frac{1}{2}\varepsilon ^{2}\left\langle P\left(
s\right) x_{1}\left( s\right) ,x_{1}\left( s\right) \right\rangle , \\
\Lambda _{10}\left( s\right) &=&\Bigg [\left\langle \left( \varepsilon
x_{1}\left( s\right) +\frac{1}{2}\varepsilon ^{2}x_{2}\left( s\right)
\right) ,\Lambda _{5}\left( s\right) \right\rangle +\varepsilon
^{2}\left\langle x_{1}\left( s\right) ,\Lambda _{6}\left( s\right)
\right\rangle \\
&&+\frac{1}{2}\varepsilon ^{2}\left\langle \Lambda _{7}\left( s\right)
x_{1}\left( s\right) ,x_{1}\left( s\right) \right\rangle +\Lambda _{8}\left(
s\right) \Bigg ],
\end{eqnarray*}%
in which $o\left( \varepsilon ^{2}\right) $ does not involve the terms $%
x_{1}\left( \cdot \right) $ and $x_{2}\left( \cdot \right) .$ Note that in
BSDE (\ref{var1}), there appears the term $x_{1}^{\top }\left( s\right)
\Lambda _{3}\left( s\right) x_{1}\left( s\right) .$ Hence, we make use of
Taylor's expansion to

\begin{eqnarray*}
&&f\left( s,\bar{x}\left( s\right) +\varepsilon x_{1}\left( s\right) +\frac{1%
}{2}\varepsilon ^{2}x_{2}\left( s\right) ,\bar{y}\left( s\right) +\Lambda
_{9}\left( s\right) ,\bar{z}\left( s\right) +\Lambda _{10}\left( s\right)
,u^{\varepsilon }\left( s\right) \right) \\
&=&f\left( s,\bar{x}\left( s\right) ,\bar{y}\left( s\right) ,\bar{z}\left(
s\right) ,\bar{u}\left( s\right) \right) \\
&&+\left[ \varepsilon x_{1}\left( s\right) +\frac{1}{2}\varepsilon
^{2}x_{2}\left( s\right) ,\Lambda _{9}\left( s\right) ,\Lambda _{10}\left(
s\right) ,\varepsilon v\left( s\right) \right] \cdot \left[ f_{x}\left(
s\right) ,f_{y}\left( s\right) ,f_{z}\left( s\right) ,f_{u}\left( s\right) %
\right] ^{\top } \\
&&+\frac{1}{2}\left[ \varepsilon x_{1}\left( s\right) +\frac{1}{2}%
\varepsilon ^{2}x_{2}\left( s\right) ,\Lambda _{9}\left( s\right) ,\Lambda
_{10}\left( s\right) \right] \cdot \mathbf{H}_{1}f\left( s\right) \\
&&\cdot \left[ \varepsilon x_{1}\left( s\right) +\frac{\varepsilon ^{2}}{2}%
x_{2}\left( s\right) ,\Lambda _{9}\left( s\right) ,\Lambda _{10}\left(
s\right) \right] ^{\top }+\varepsilon v\left( s\right) f_{yu}\left( s\right)
\Lambda _{9}\left( s\right) \\
&&+\varepsilon v\left( s\right) f_{xu}\left( s\right) \left( \varepsilon
x_{1}\left( s\right) +\frac{1}{2}\varepsilon ^{2}x_{2}\left( s\right)
\right) +\varepsilon v\left( s\right) f_{zu}\Lambda _{10}\left( s\right) +%
\frac{1}{2}\varepsilon ^{2}v^{2}\left( s\right) f_{uu}\left( s\right) ,
\end{eqnarray*}%
where the Hessian matrix $\mathbf{H}_{1}$ is with respect to $\left(
x,y,z\right) .$

Then, we obtain
\begin{eqnarray*}
\Gamma \left( t\right) &=&b_{x}^{\top }\left( t\right) p\left( t\right)
+f_{x}\left( t\right) +f_{y}\left( t\right) p\left( t\right) \\
&&+\sigma _{x}\left( t\right) ^{\top }q\left( t\right) +f_{z}\left( t\right)
\left[ \sigma _{x}\left( t\right) ^{\top }+q\left( t\right) \right] , \\
\Pi \left( t\right) &=&P\left( t\right) b_{x}\left( t\right) +b_{x}^{\top
}\left( t\right) P\left( t\right) +f_{y}\left( t\right) P\left( t\right)
+b_{xx}^{\top }\left( t\right) p\left( t\right) \\
&&+[Q\left( t\right) \sigma _{x}\left( t\right) +\sigma _{x}\left( t\right)
^{\top }Q\left( t\right) +\sigma _{xx}\left( t\right) ^{\top }q\left(
t\right) +\sigma _{x}\left( t\right) ^{\top }P\left( t\right) \sigma
_{x}\left( t\right) ] \\
&&+[f_{z}\left( t\right) P\left( t\right) \sigma _{x}\left( t\right)
+f_{z}\left( t\right) \sigma _{x}\left( t\right) ^{\top }P\left( t\right)
+f_{z}\left( t\right) Q\left( t\right) +f_{z}\left( t\right) \sigma
_{xx}\left( t\right) ^{\top }p\left( t\right) ] \\
&&+\left[ I_{n\times n},p\left( t\right) ,\sigma _{x}\left( t\right) ^{\top
}p\left( t\right) +q\left( t\right) \right] \cdot \mathbf{H}_{1}f\left(
t\right) \\
&&\cdot \left[ I_{n\times n},p\left( t\right) ,\sigma _{x}\left( t\right)
^{\top }p\left( t\right) +q\left( t\right) \right] ^{\top },
\end{eqnarray*}%
where $I_{n\times n}$ denotes the identity matrix.%
%
%
%
%
%
%
%
%
%
%
%
%
%
%
%
%
%
%
%
%
%
%
%
%
%
%
%
%
%
%
%
%
%
%
%
%
%
%
%
%
%
%
%
%
%
%
%
%
%
%
%
%
%
%
%
%
%
%
%
%
%
%
%
%
%
%
%
%
%
%
%
%
%
%
%
%
%
%
%
%
%
%
%
%
%
%
%
%
%
%
%
%
%
%
%
%
%
%
%
%
%
%
%
%
%
%
%
%
%
%
%
%
%
%
%
%
%
%
%
%
%
%
%
%
%
%
%
%
%
%
%
%
%
%
%
%
%
%
%
%
%
%
%
%
%Now consider the following BSDE:%
%\begin{eqnarray}
%\hat{y}\left( t\right) &=&\int_{t}^{T}\Big [f_{y}\left( s\right) \hat{y}%
%\left( s\right) +f_{z}\left( s\right) \hat{z}\left( s\right) +p\left(
%s\right) b_{u}\left( s\right) v\left( s\right)  \notag \\
%&&+p\left( s\right) v^{2}\left( s\right) b_{uu}\left( s\right) +q\left(
%s\right) \sigma _{u}\left( s\right) v\left( s\right) +q\left( s\right)
%v^{2}\left( s\right) \sigma _{uu}\left( s\right)  \notag \\
%&&+\frac{1}{2}\sigma _{u}^{2}\left( s\right) v^{2}\left( s\right) P\left(
%s\right) +\big [2p\left( s\right) v\left( s\right) b_{xx}\left( s\right)
%+2q\left( s\right) v\left( s\right) \sigma _{xu}\left( s\right)  \notag \\
%&&+P\left( s\right) b_{u}\left( s\right) v\left( s\right) +P\left( s\right)
%\sigma _{x}\left( s\right) \sigma _{u}\left( s\right) +Q\left( s\right)
%\sigma _{u}\left( s\right) v\left( s\right) \big ]x_{1}\left( s\right) \Big ]%
%\mathrm{d}s  \notag \\
%&&-\int_{t}^{T}\hat{z}\left( s\right) \mathrm{d}W\left( s\right) .
%\label{var3}
%\end{eqnarray}

%\begin{remark}
%The adjoint equations derived here is the same as in Hu \cite{Hms}.
%\end{remark}

In the classical theory of optimal control for FBSDEs (cf \cite{Wu1998,
Wu2013}), there generally appear two groups of the first-order adjoint
equations, for instance$\left( \mathfrak{p}\left( \cdot \right) ,\mathfrak{q}%
\left( \cdot \right) \right) $ in Eqs. (\ref{firstadjoint}). The following
proposition will establish the relationship between them with $p\left( \cdot
\right) $ from (\ref{adj1}), which is very useful to study the connection
between maximum principle and dynamic programming (see Theorem \ref{MPDP}
below).

\begin{proposition}
\label{pro1}Suppose that Assumptions \emph{(A1)-(A2)} are in force. Then we
have%
\begin{equation*}
p\left( s\right) =-\frac{\mathfrak{p}^{\top }\left( s\right) }{\mathfrak{q}%
\left( s\right) },\text{ }s\in \left[ t,T\right] ,\text{ }P\text{-}a.s.,
\end{equation*}%
where $p\left( \cdot \right) $ and $\left( \mathfrak{p}\left( \cdot \right) ,%
\mathfrak{q}\left( \cdot \right) \right) $ are solutions to FBSDEs (\ref%
{adj1}) and (\ref{firstadjoint}), respectively.
\end{proposition}

The proof is just to apply the It\^{o}'s formula to $-\mathfrak{p}^{T}\left(
s\right) /\mathfrak{q}\left( s\right) ,$ so we omit it.

We define the \textit{classical} Hamiltonian function:

\begin{equation*}
H\left( t,x,y,z,u,p,q\right) =\left\langle p,b\left( t,x,u\right)
\right\rangle +\left\langle q,\sigma \left( t,x,u\right) \right\rangle
+f\left( t,x,y,z,u\right) ,
\end{equation*}%
where $\left( t,x,y,z,u,p,q\right) \in \left[ 0,T\right] \times \mathbb{R}%
^{n}\mathbb{\times R\times R\times R}^{k}\mathbb{\times R}^{n}\mathbb{\times
R}^{n}.$%

Then, we have%
\begin{eqnarray*}
\hat{y}^{\varepsilon }\left( t\right) &=&o\left( \varepsilon ^{2}\right)
+\int_{t}^{T}\Bigg [f\Big (s,\bar{x}\left( s\right) ,\bar{y}\left( s\right) ,%
\bar{z}\left( s\right) \\
&&+\varepsilon \left\langle p\left( s\right) \sigma _{u}\left( s\right)
,v\left( s\right) \right\rangle +\frac{\varepsilon ^{2}}{2}\left\langle
p\left( s\right) \sigma _{uu}\left( s\right) v\left( s\right) ,v\left(
s\right) \right\rangle ,\bar{u}\left( s\right) \Big ) \\
&&-f\left( s,\bar{x}\left( s\right) ,\bar{y}\left( s\right) ,\bar{z}\left(
s\right) ,\bar{u}\left( s\right) \right) +f_{y}\left( s\right) \hat{y}%
^{\varepsilon }\left( s\right) +f_{z}\left( s\right) \hat{z}^{\varepsilon
}\left( s\right) \\
&&+\varepsilon ^{2}\big <\big [p\left( s\right) b_{xu}\left( s\right)
+q\left( s\right) \sigma _{xu}\left( s\right) +P\left( s\right) b_{u}\left(
s\right) \\
&&+P\left( s\right) \sigma _{x}\left( s\right) \sigma _{u}\left( s\right)
+Q\left( s\right) \sigma _{u}\left( s\right) \big ]x_{1}\left( s\right)
,v\left( s\right) \big > \\
&&+\varepsilon \left\langle \left[ p\left( s\right) b_{u}\left( s\right)
+q\left( s\right) \sigma _{u}\left( s\right) +f_{u}\left( s\right) \right]
,v\left( s\right) \right\rangle \\
&&+\frac{\varepsilon ^{2}}{2}\left\langle \left[ q\left( s\right) \sigma
_{uu}\left( s\right) +p\left( s\right) b_{uu}\left( s\right) +\varepsilon
^{2}f_{uu}\left( s\right) +\sigma _{u}^{\top }\left( s\right) P\left(
s\right) \sigma _{u}\left( s\right) \right] v\left( s\right) ,v\left(
s\right) \right\rangle \\
&&+\varepsilon v^{\top }\left( s\right) f_{xu}\left( s\right) \left[
\varepsilon x_{1}\left( s\right) +\frac{\varepsilon ^{2}}{2}x_{2}\left(
s\right) \right] +\varepsilon v^{\top }\left( s\right) f_{yu}\left( s\right)
\Lambda _{9}\left( s\right) \\
&&+\varepsilon v^{\top }\left( s\right) f_{zu}\left( s\right) \Lambda
_{10}\left( s\right) \Bigg ]\mathrm{d}s-\int_{t}^{T}\hat{z}^{\varepsilon
}\left( s\right) \mathrm{d}W\left( s\right) .
\end{eqnarray*}%
Namely,
\begin{eqnarray}
\hat{y}^{\varepsilon }\left( t\right) &=&o\left( \varepsilon ^{2}\right)
+\int_{t}^{T}\Big \{f_{y}\left( s\right) \hat{y}^{\varepsilon }\left(
s\right) +f_{z}\left( s\right) \hat{z}^{\varepsilon }\left( s\right)
+\varepsilon ^{2}x_{1}^{\top }\left( s\right) \big [H_{xu}\left( s\right)
\notag \\
&&+P\left( s\right) \sigma _{x}\left( s\right) \sigma _{u}\left( s\right)
+Q\left( s\right) \sigma _{u}\left( s\right) +p\left( t\right) f_{yu}\left(
s\right) +P\left( s\right) b_{u}\left( s\right) \big ]v\left( s\right)
\notag \\
&&+\varepsilon \left[ H_{u}\left( s\right) +f_{z}\left( s\right) p\left(
s\right) \sigma _{u}\left( s\right) \right] v\left( s\right)  \notag \\
&&+\frac{\varepsilon ^{2}}{2}v^{\top }\left( s\right) \left[ H_{uu}\left(
s\right) +\sigma _{u}^{2}\left( s\right) P\left( s\right) +f_{zz}\left(
s\right) p^{2}\left( s\right) \sigma _{u}^{2}\left( s\right) +f_{z}\left(
s\right) p\sigma _{uu}\right] v\left( s\right)  \notag \\
&&+\varepsilon ^{2}v^{\top }\left( s\right) f_{zu}\left( s\right) \Big [%
\left( p\left( t\right) \sigma _{x}\left( t\right) +q\left( t\right) \right)
x_{1}\left( s\right) +p\left( t\right) \sigma _{u}\left( s\right) v\left(
s\right) \Big ]\Big \}\mathrm{d}s  \notag \\
&&-\int_{t}^{T}\hat{z}^{\varepsilon }\left( s\right) \mathrm{d}W\left(
s\right) .  \label{vbsde}
\end{eqnarray}

\begin{remark}
Note that FBSDEs (\ref{vbsde}) are somewhat different from (22) in Hu \cite%
{Hms}. Specifically, the term $A_{4}x_{1}\left( s\right) I_{E_{\varepsilon
}}\left( s\right) $ disappears in (22) since
\begin{equation*}
\mathbb{E}\left[ \left( \int_{0}^{T}\left\vert A_{4}\left( s\right)
x_{1}\left( s\right) I_{E_{\varepsilon }}\left( s\right) \right\vert \mathrm{%
d}s\right) ^{\beta }\right] =o\left( \varepsilon ^{\beta }\right) \text{ for
}\beta \geq 2
\end{equation*}
in \cite{Hms} by using spike variational approach. Nevertheless, the
corresponding term in our paper is just $\varepsilon ^{2}x_{1}^{\top }\left(
s\right) \Lambda _{2}\left( s\right) $. We will see a moment later that this
term is needed to define an extensive \textquotedblleft Hamiltonian
function\textquotedblright\ as follows.
\end{remark}

Define
\begin{eqnarray*}
\mathbb{H}\left( t,x,y,z,u,p,q,P,Q\right) &\triangleq &H_{xu}\left(
t,x,y,z,u,p,q\right) +b_{u}\left( t,x,u\right) P \\
&&+Q\sigma _{u}\left( t,x,u\right) +P\sigma _{x}\left( t,x,u\right) \sigma
_{u}\left( t,x,u\right) \\
&&+f_{yu}\left( t,x,y,z,u\right) p+f_{zu}\left( t,x,y,z,u\right) \left[
p\sigma _{x}\left( t,x,u\right) +q\right] ,
\end{eqnarray*}%
where
\begin{equation*}
\left( t,x,y,z,p,q,P,Q\right) \in \left[ 0,T\right] \times \mathbb{R}%
^{n}\times \mathbb{R}\times \mathbb{R}\times \mathbb{R}^{n}\times \mathbb{R}%
^{n}\times \mathbb{R}^{n\times n}\times \mathbb{R}^{n\times n}.
\end{equation*}

%\begin{remark}
%Observe that $\mathbb{H}$ is slightly different from $\mathbb{S}$ in Zhang
%et al. \cite{ZZconvex}, that is,
%\begin{eqnarray*}
%\mathbb{S}\left( t,x,u,y_{1},z_{1},y_{2},z_{2}\right) &=&H_{xu}\left(
%t,x,y,z,p,q\right) +b_{u}\left( t,x,u\right) y_{2} \\
%&&+z_{2}\sigma _{u}\left( t,x,u\right) +y_{2}\sigma _{x}\left( t,x,u\right)
%\sigma _{u}\left( t,x,u\right) ,
%\end{eqnarray*}%
%where $\left( y_{2},z_{2}\right) \in \mathbb{R}^{n\times n}\times \mathbb{R}%
%^{n\times n}.$
%\end{remark}

We now give the adjoint equation for BSDE (\ref{vbsde}) as follows:%
\begin{equation}
\left\{
\begin{array}{rcl}
\mathrm{d}\chi \left( t\right) & = & f_{y}\left( t\right) \chi \left(
t\right) \mathrm{d}t+f_{z}\left( t\right) \chi \left( t\right) \mathrm{d}%
W\left( t\right) , \\
\chi \left( 0\right) & = & 1.%
\end{array}%
\right.  \label{adjbsde}
\end{equation}

\begin{lemma}
\label{keylemma} Under the Assumptions \emph{(A1)-(A2)}, SDE (\ref{adjbsde})
admits a unique adapted strong solution $\chi \left( t\right) \in \mathcal{S}%
^{2}(0,T;\mathbb{R}).$ Moreover, we have
\begin{eqnarray}
\mathbb{E}\left[ \sup_{0\leq t\leq T}\left\vert \chi \left( t\right)
\right\vert ^{l}\right] &<&\infty ,\text{ }l\geq 1,  \label{keyest1} \\
\mathbb{E}\left[ \sup_{0\leq s\leq T}\left\vert \chi \left( s\right)
\right\vert ^{2}\right] ^{4} &<&\infty .  \label{keyest2}
\end{eqnarray}
\end{lemma}

\paragraph{Proof.}

The first inequality can be obtained from Theorem 6.16 of \cite{YZ}. We deal
with the second one. By It\^{o}'s formula, we have
\begin{equation*}
\sup_{0\leq t\leq T}\chi ^{2}\left( t\right) =1+\int_{0}^{T}\chi ^{2}\left(
t\right) \left( 2\left\vert f_{y}\left( t\right) \right\vert
+f_{z}^{2}\left( t\right) \right) \mathrm{d}t+\sup_{0\leq t\leq
T}\int_{0}^{t}2\chi ^{2}\left( s\right) f_{z}\left( s\right) \mathrm{d}%
W\left( s\right) .
\end{equation*}%
It follows that%
\begin{eqnarray*}
\left( \sup_{0\leq t\leq T}\chi ^{2}\left( t\right) \right) ^{4} &=&C\Bigg [%
1+\left( \int_{0}^{T}\chi ^{2}\left( t\right) \left( 2\left\vert f_{y}\left(
t\right) \right\vert +f_{z}^{2}\left( t\right) \right) \mathrm{d}t\right)
^{4} \\
&&+\left( \sup_{0\leq t\leq T}\int_{0}^{t}2\chi ^{2}\left( s\right)
f_{z}\left( s\right) \mathrm{d}W\left( s\right) \right) ^{4}\Bigg ].
\end{eqnarray*}%
But by the B-D-G inequality, we get%
\begin{equation*}
\mathbb{E}\left( \sup_{0\leq t\leq T}\int_{0}^{t}2\chi ^{2}\left( s\right)
f_{z}\left( s\right) \mathrm{d}W\left( s\right) \right) ^{4}\leq C\mathbb{E}%
\left[ \int_{0}^{t}\chi ^{8}\left( s\right) \left( f_{z}^{2}\left( t\right)
\right) ^{2}\mathrm{d}s\right] .
\end{equation*}%
Thus
\begin{eqnarray*}
\mathbb{E}\left[ \left( \sup_{0\leq t\leq T}\chi ^{2}\left( t\right) \right)
^{4}\right] &=&C\mathbb{E}\left[ 1+\left( \int_{0}^{T}\chi ^{2}\left(
t\right) \mathrm{d}t\right) ^{4}+\int_{0}^{t}\chi ^{8}\left( s\right)
\mathrm{d}s\right] \\
&\leq &C\mathbb{E}\left[ 1+\int_{0}^{t}\chi ^{8}\left( s\right) \mathrm{d}s%
\right] \\
&\leq &C\mathbb{E}\left[ 1+T\sup_{0\leq t\leq T}\left\vert \chi \left(
t\right) \right\vert ^{8}\right] \\
&<&\infty .
\end{eqnarray*}%
The second estimation comes from the H\"{o}lder inequality. We complete the
proof.\hfill $\Box $

Set%
\begin{eqnarray*}
y^{\varepsilon }(t) &=&\bar{y}\left( t\right) +\left[ p^{\top }\left(
t\right) \left( \varepsilon x_{1}\left( t\right) +\frac{\varepsilon ^{2}}{2}%
x_{2}\left( t\right) \right) +\frac{\varepsilon ^{2}}{2}x_{1}^{\top }\left(
t\right) P\left( t\right) x_{1}\left( t\right) \right] +\hat{y}^{\varepsilon
}\left( t\right) , \\
z^{\varepsilon }(t) &=&\bar{z}\left( t\right) +\Bigg [\Lambda _{5}\left(
s\right) \left( \varepsilon x_{1}\left( s\right) +\frac{\varepsilon ^{2}}{2}%
x_{2}\left( s\right) \right) +\varepsilon ^{2}\Lambda _{6}\left( s\right)
x_{1}\left( s\right) \\
&&+\frac{\varepsilon ^{2}}{2}x_{1}^{\top }\left( t\right) \Lambda _{7}\left(
s\right) x_{1}\left( s\right) +\Lambda _{8}\left( s\right) \Bigg ]+\hat{z}%
^{\varepsilon }\left( t\right) .
\end{eqnarray*}%
We are able to give the variational equations as follows:
\begin{eqnarray*}
y_{1}^{\varepsilon }\left( t\right) &=&\varepsilon p^{\top }\left( t\right)
x_{1}\left( t\right) , \\
z_{1}^{\varepsilon }\left( t\right) &=&\varepsilon \left[ x_{1}^{\top
}\left( s\right) \Lambda _{5}\left( s\right) +p^{\top }\left( t\right)
\sigma _{u}\left( t\right) v\left( t\right) \right]
\end{eqnarray*}%
and%
\begin{eqnarray*}
y_{2}^{\varepsilon }\left( t\right) &=&\frac{\varepsilon ^{2}}{2}\left[
p^{\top }\left( t\right) x_{2}\left( t\right) +x_{1}^{T}\left( t\right)
P\left( t\right) x_{1}\left( t\right) \right] +\hat{y}^{\varepsilon }\left(
t\right) , \\
z_{2}^{\varepsilon }\left( t\right) &=&\frac{\varepsilon ^{2}}{2}\Bigg [%
\Lambda _{5}\left( s\right) x_{2}\left( s\right) +2\Lambda _{6}\left(
s\right) x_{1}\left( s\right) \\
&&+x_{1}^{\top }\left( s\right) \Lambda _{7}\left( s\right) x_{1}\left(
s\right) +p^{\top }\left( t\right) v^{\top }\left( t\right) \sigma
_{uu}\left( t\right) v\left( t\right) \Bigg ]+\hat{z}^{\varepsilon }\left(
t\right) .
\end{eqnarray*}%
Obviously, we have%
\begin{equation*}
y^{\varepsilon }(0)-\bar{y}\left( 0\right) =\hat{y}^{\varepsilon }\left(
0\right) \geq 0,\text{ since the definition of value function.}
\end{equation*}%

\begin{lemma}
\label{est} Under the Assumptions \emph{(A1)-(A2)}, we have the following
estimation%
\begin{equation}
\mathbb{E}\left[ \sup_{0\leq t\leq T}\left\vert \hat{y}^{\varepsilon }\left(
t\right) \right\vert ^{2}+\int_{0}^{T}\left\vert \hat{z}^{\varepsilon
}\left( s\right) \right\vert ^{2}\mathrm{d}s\right] =O\left( \varepsilon
^{2}\right) .  \label{estv}
\end{equation}%
%
%
%
%
%
%
%
%
%
%
%
%
%
%
%
%
%
%
%
%
%
%
%
%
%
%
%
%
%
%
%
%
%
%
%
%
%
%
%
%
%
%
%
%
%
%
%
%
%
%
%
%
%
%
%
%
%
%
%
%
%
%
%
%
%
%
%
%
%
%
%
%
%
%
%
%
%
%
%
%
%
%
%
%
%
%
%
%
%
%
%
%
%
%
%
%
%
%
%
%
%
%
%
%
%
%
%
%
%
%
%
%
%
%
%
%
%
%
%
%
%
%
%
%
%\begin{equation}
%\mathbb{E}\left[ \sup_{0\leq t\leq T}\left\vert y^{\varepsilon }\left(
%t\right) -\hat{y}^{\varepsilon }\left( t\right) -\bar{y}\left( t\right)
%\right\vert ^{2}+\int_{0}^{T}\left\vert z^{\varepsilon }\left( s\right) -%
%\hat{z}^{\varepsilon }\left( s\right) -\bar{z}\left( s\right) \right\vert
%^{2}\mathrm{d}s\right] =o\left( \varepsilon ^{2}\right) ,  \label{estv2}
%\end{equation}
\end{lemma}

\paragraph{Proof.}

To prove (\ref{estv}), we consider (\ref{vbsde}) again. From assumptions
(A1)-(A2), one can check that the adjoint equations (\ref{adj1}) and (\ref%
{adj2}) have a unique adapted strong solution, respectively. Furthermore, by
classical approach, we are able to get the following estimates for $\beta
\geq 2,$%
\begin{equation}
\mathbb{E}\left[ \sup_{0\leq t\leq T}\left( \left\vert p\left( t\right)
\right\vert ^{\beta }+\left\vert P\left( t\right) \right\vert ^{\beta
}\right) +\left( \int_{0}^{T}\left( \left\vert q\left( t\right) \right\vert
^{2}+\left\vert Q\left( t\right) \right\vert ^{2}\right) \mathrm{d}t\right)
^{\frac{\beta }{2}}\right] <+\infty .  \label{estadj}
\end{equation}%
Applying Lemma \ref{estlemma} to (\ref{vbsde}), we get the desired result.
Indeed, since there appears a term
\begin{equation*}
\varepsilon \left[ H_{u}\left( s\right) +f_{z}\left( s\right) p\left(
s\right) \sigma _{u}\left( s\right) \right] v\left( s\right) ,
\end{equation*}
in BSDE (\ref{vbsde}), so we have the estimation with $O\left( \varepsilon
^{2}\right) .$ We complete the proof. %we reformulate (\ref%
\hfill $\Box $

We shall derive a variational inequality which is crucial to establish the
necessary condition for optimal control. Before this, we introduce the
following the other type of singular control using the Hamiltonian function:

\begin{definition}[Singular control in the classical sense]
\label{sc}We call a control $\breve{u}\left( \cdot \right) \in \mathcal{U}%
\left( 0,T\right) $ a singular control in the classical sense if $\breve{u}%
\left( \cdot \right) $ satisfies%
\begin{equation}
\left\{
\begin{array}{l}
\text{i) }H_{u}\left( t,\breve{x}\left( t\right) ,\breve{y}\left( t\right) ,%
\breve{z}\left( t\right) ,\breve{u}\left( t\right) ,\breve{p}\left( t\right)
,\breve{q}\left( t\right) \right) \\
\qquad +f_{z}\left( t,\breve{x}\left( t\right) ,\breve{y}\left( t\right) ,%
\breve{z}\left( t\right) \right) \cdot \sigma _{u}^{\top }\left( t,\breve{x}%
\left( t\right) ,\breve{u}\left( t\right) \right) \cdot \breve{p}\left(
t\right) =0,\text{ a.s., a.e., }t\in \left[ 0,T\right] ; \\
\text{ii) }H_{uu}\left( t,\breve{x}\left( t\right) ,\breve{y}\left( t\right)
,\breve{z}\left( t\right) ,\breve{u}\left( t\right) ,\breve{p}\left(
t\right) ,\breve{q}\left( t\right) \right) \\
\qquad +\sigma _{u}\left( t,\breve{x}\left( t\right) ,\breve{u}\left(
t\right) \right) \cdot \breve{P}\left( s\right) \cdot \sigma _{u}^{\top
}\left( t,\breve{x}\left( t\right) ,\breve{u}\left( t\right) \right) \\
\qquad +f_{z}\left( t,\breve{x}\left( t\right) ,\breve{y}\left( t\right) ,%
\breve{z}\left( t\right) \right) \cdot \sigma _{uu}^{\top }\left( t,\breve{x}%
\left( t\right) ,\breve{u}\left( t\right) \right) \cdot \breve{p}\left(
s\right) \\
\qquad +2f_{zu}\left( t,\breve{x}\left( t\right) ,\breve{y}\left( t\right) ,%
\breve{z}\left( t\right) \right) \cdot \breve{p}^{\top }\left( t\right)
\cdot \sigma _{u}\left( t,\breve{x}\left( t\right) ,\breve{u}\left( t\right)
\right) \\
\qquad +f_{zz}\left( t,\breve{x}\left( t\right) ,\breve{y}\left( t\right) ,%
\breve{z}\left( t\right) \right) \cdot \sigma _{u}^{\top }\left( t,\breve{x}%
\left( t\right) ,\breve{u}\left( t\right) \right) \breve{p}\left( t\right)
\cdot \breve{p}^{\top }\left( t\right) \sigma _{u}\left( t,\breve{x}\left(
t\right) ,\breve{u}\left( t\right) \right) =0, \\
\qquad \text{a.s., a.e., }t\in \left[ 0,T\right] ;%
\end{array}%
\right.  \label{scc}
\end{equation}%
where $\left( \breve{x}\left( \cdot \right) ,\breve{y}\left( \cdot \right) ,%
\breve{z}\left( \cdot \right) \right) $ denotes the state trajectories
driven by $\breve{u}\left( \cdot \right) .$ Moreover, $\left( \breve{p}%
\left( \cdot \right) ,\breve{q}\left( \cdot \right) \right) $ and $\left(
\breve{P}\left( \cdot \right) ,\breve{Q}\left( \cdot \right) \right) $
denote the adjoint processes given respectively by (\ref{adj1}) and (\ref%
{adj2}) with $\left( \bar{x}\left( t\right) ,\bar{y}\left( t\right) ,\bar{z}%
\left( t\right) ,\bar{u}\left( t\right) \right) $ replaced by \newline
$(\breve{x}\left( \cdot \right) ,\breve{y}\left( \cdot \right) ,\breve{z}%
\left( \cdot \right) ,\breve{u}\left( \cdot \right) )$. If this $\breve{u}%
\left( \cdot \right) $ is also optimal, then we call it a singular optimal
control in the classical sense.
\end{definition}

\begin{remark}
Hu \cite{Hms} first considers the forward-backward stochastic control
problem whenever the diffusion term $\sigma \left( t,x,u\right) $ depends on
the control variable $u$ with non-convex control domain. In order to to
establish the stochastic maximum principle, he introduces the $\mathcal{H}$%
-function of the following type:%
\begin{eqnarray*}
\mathcal{H}\left( t,x,y,z,u,p,q,P\right) &\triangleq &pb\left( t,x,u\right)
+q\sigma \left( t,x,u\right) \\
&&+\frac{1}{2}\left( \sigma \left( t,x,u\right) -\sigma \left( t,\bar{x},%
\bar{u}\right) \right) ^{\top }P\left( \sigma \left( t,x,u\right) -\sigma
\left( t,\bar{x},\bar{u}\right) \right) \\
&&+f\left( t,x,y,z+p\left( \sigma \left( t,x,u\right) -\sigma \left( t,\bar{x%
},\bar{u}\right) \right) ,u\right) .
\end{eqnarray*}%
Note that this Hamiltonian function is slightly different from Peng 1990
\cite{Peng1990}. The main difference of this variational equations with
those in (Peng 1990) \cite{Peng1990} appears in the term $p\left( t\right)
\delta \sigma \left( t\right) I_{E_{\varepsilon }}\left( t\right) $ (the
similar term $\varepsilon p\left( t\right) \sigma _{u}\left( t\right)
v\left( t\right) +\frac{\varepsilon ^{2}}{2}p\left( t\right) \sigma
_{uu}\left( t\right) v^{2}\left( t\right) $ in our paper) in variational
equation for BSDE and maximum principle for the definition of $p\left(
t\right) $ in the variation of $z$, which is $O(\varepsilon )$ for any order
expansion of $f$. So it is not helpful to use the second-order Taylor
expansion for treating this term. The stochastic maximum principle (see \cite%
{Hms}) says that if $\left( \breve{x}\left( t\right) ,\breve{y}\left(
t\right) ,\breve{z}\left( t\right) ,\breve{u}\left( t\right) \right) $ is an
optimal pair, then%
\begin{eqnarray}
&&\mathcal{H}\left( t,\breve{x}\left( t\right) ,\breve{y}\left( t\right) ,%
\breve{z}\left( t\right) ,\breve{u}\left( t\right) ,\breve{p}\left( t\right)
,\breve{q}\left( t\right) ,\breve{P}\left( t\right) \right)  \notag \\
&=&\min_{u\in U}\mathcal{H}\left( t,\breve{x}\left( t\right) ,\breve{y}%
\left( t\right) ,\breve{z}\left( t\right) ,u,\breve{p}\left( t\right) ,%
\breve{q}\left( t\right) ,\breve{P}\left( t\right) \right) .  \label{optimp}
\end{eqnarray}%
Apparently, Definition \ref{sc} says that a singular control in the
classical sense is the real one that fulfils trivially the first and
second-order necessary conditions in classical optimization theory dealing
with the maximization problem (\ref{optimp}), namely,%
\begin{equation}
\left\{
\begin{array}{c}
\mathcal{H}_{u}\left( t,\breve{x}\left( t\right) ,\breve{y}\left( t\right) ,%
\breve{z}\left( t\right) ,\breve{u}\left( t\right) ,\breve{p}\left( t\right)
,\breve{q}\left( t\right) ,\breve{P}\left( t\right) \right) =0,\text{ a.s.,
a.e., }t\in \left[ 0,T\right] ; \\
\mathcal{H}_{uu}\left( t,\breve{x}\left( t\right) ,\breve{y}\left( t\right) ,%
\breve{z}\left( t\right) ,\breve{u}\left( t\right) ,\breve{p}\left( t\right)
,\breve{q}\left( t\right) ,\breve{P}\left( t\right) \right) =0,\text{ a.s.,
a.e., }t\in \left[ 0,T\right] .%
\end{array}%
\right.  \label{equsc}
\end{equation}%
It is easy to verify that (\ref{scc}) is equivalent to (\ref{equsc}).
Certainly, one could investigate stochastic singular optimal controls for
forward-backward stochastic systems in other senses, say, in the sense of
process in Skorohod space, which can be seen in Zhang \cite{Zhangsingular}
via viscosity solution approach (Hamilton-Jacobi-Bellman inequality), or in
the sense of Pontryagin-type maximum principle (cf Tang \cite{Tang}). As
this complete remake of the various topics is much longer than the present
paper, it will be reported elsewhere.
\end{remark}

\begin{lemma}[Variational inequality]
\label{varineq} Under the Assumptions \emph{(A1)-(A2)}, it holds that%
\begin{eqnarray}
0 &\leq &\mathbb{E}\Bigg \{\chi \left( T\right) o\left( \varepsilon
^{2}\right) +\int_{0}^{T}\chi \left( s\right) \Big [\varepsilon \left\langle
f_{z}\left( s\right) \sigma _{u}^{\top }\left( s\right) p\left( s\right)
+H_{u}\left( s\right) ,v\left( s\right) \right\rangle  \notag \\
&&+\frac{\varepsilon ^{2}}{2}\Big (\left\langle f_{z}\left( s\right) \sigma
_{uu}^{\top }\left( s\right) p\left( s\right) v\left( s\right) ,v\left(
s\right) \right\rangle +2\left\langle \mathbb{H}\left( s\right) x_{1}\left(
s\right) ,v\left( s\right) \right\rangle  \notag \\
&&+2\left\langle f_{zu}\left( s\right) p\left( t\right) \sigma _{u}\left(
t\right) v\left( s\right) ,v\left( s\right) \right\rangle +\left\langle
\mathbb{\tilde{H}}\left( s\right) v\left( s\right) ,v\left( s\right)
\right\rangle \Big )\Big ]\mathrm{d}s\Bigg \},  \label{varin}
\end{eqnarray}%
where
\begin{eqnarray*}
\mathbb{H}\left( s\right) &=&\mathbb{H}\left( t,\bar{x}\left( s\right) ,\bar{%
y}\left( s\right) ,\bar{z}\left( s\right) ,\bar{u}\left( s\right) ,p\left(
s\right) ,q\left( s\right) ,P\left( s\right) ,Q\left( s\right) \right) , \\
\mathbb{\tilde{H}}\left( s\right) &=&H_{uu}\left( s\right) +\sigma
_{u}^{\top }\left( s\right) P\left( s\right) \sigma _{u}^{\top }\left(
s\right) +f_{zz}\left( s\right) \sigma _{u}^{\top }\left( s\right) p^{\top
}\left( s\right) p\left( s\right) \sigma _{u}\left( s\right) .
\end{eqnarray*}
\end{lemma}

\paragraph{Proof.}

Using It\^{o}'s formula to $\left\langle \chi \left( s\right) ,\hat{y}%
^{\varepsilon }\left( s\right) \right\rangle $ on $\left[ 0,T\right] ,$ we
get the desired result.\hfill $\Box $

\begin{theorem}
\label{th1}Assume that \emph{(A1)-(A2)} hold. If $\bar{u}\left( \cdot
\right) \in \mathcal{U}\left( 0,T\right) $ is a singular optimal control in
the classical sense, then%
\begin{equation}
0\leq \mathbb{E}\Bigg [\int_{0}^{T}\left\langle \chi \left( s\right) \mathbb{%
H}\left( s\right) x_{1}\left( s\right) ,v\left( s\right) \right\rangle
\mathrm{d}s\Bigg ],  \label{scm}
\end{equation}%
for any $v\left( \cdot \right) =u\left( \cdot \right) -\bar{u}\left( \cdot
\right) ,$ $u\left( \cdot \right) \in \mathcal{U}\left( 0,T\right) .$
\end{theorem}

\paragraph{Proof.}

According to the definition of value function, we have%
\begin{equation*}
\frac{J\left( u^{\varepsilon }\right) -J\left( \bar{u}\right) }{\varepsilon
^{2}}=\frac{y^{\varepsilon }(0)-\bar{y}\left( 0\right) }{\varepsilon ^{2}}=%
\frac{\hat{y}^{\varepsilon }\left( 0\right) }{\varepsilon ^{2}}\geq 0.
\end{equation*}%
Letting $\varepsilon \rightarrow 0+,$ we get the desired result from
Definition \ref{sc} and Lemma \ref{varineq}. \hfill $\Box $

\begin{remark}
Clearly, if $f$ does not depend on $\left( y,z\right) $, then $\chi \left(
\cdot \right) \equiv 1.$ Consequently, (\ref{scm}) reduces to
\begin{equation*}
\mathbb{E}\left[ \int_{0}^{T}\left\langle \mathbb{H}\left( s\right)
x_{1}\left( s\right) ,v\left( s\right) \right\rangle \mathrm{d}s\right] \geq
0,
\end{equation*}%
which is just the classical case studied in Zhang et al. \cite{ZZconvex} for
classical stochastic control problems. Meanwhile, our result actually
extends Peng \cite{Peng1993} to second order case.%\begin{remark}
%The Hamilton function (\ref{H}) is different from the one in Zhang and Shi
%[]. Specifically, in (\ref{H}), $g$ depends on $z+p\left( \sigma \left(
%t,x,u\right) -\sigma \left( t,\bar{x},\bar{u}\right) \right) $ instead of $z$
%only. If the control variable enters in $g,$ we will discuss in Section \ref%
%{sec5}.
%\end{remark}
\end{remark}

\begin{remark}
Recall that, for deterministic system, it is possible to derive pointwise
necessary conditions for optimal controls via the first suitable
integral-type necessary conditions and normally there is no obstacles to
establish the pointwise first-order necessary condition for optimal controls
whenever an integral type one is on the hand. Nevertheless, the classical
approach to handle the pointwise condition from the integral-type can not be
employed directly in the framework of the pointwise second-order condition
in the general stochastic setting because of certain feature the stochastic
systems owning. In order to derive the second order variational equations
for BSDE in Hu \cite{Hms}, the author there introduces two kinds of adjoint
equations and a new Hamiltonian function. The main difference of this
variational equations with those in (Peng 1990) \cite{Peng1990} lies in the
term $p\left( t\right) \delta \sigma \left( t\right) I_{E_{\varepsilon
}}\left( t\right) .$ Then, it is possible to get the maximum principle
basing one variational equation. Note that the order of the difference
between perturbed state, optimal state and first, second order state is $%
o\left( \varepsilon \right) .$
\end{remark}

%\begin{theorem}
%\label{th1}Assume that conditions \emph{(A1)--(A4}\textsl{)} hold. Let $%
%\left( p\left( t\right) ,q\left( t\right) \right) $ and $\left( P\left(
%t\right) ,Q\left( t\right) \right) $ be an admissible pair of the solutions
%to (\ref{adj1}) and (\ref{adj2}), respectively, corresponding to $\left(
%\bar{x}\left( t\right) ,\bar{y}\left( t\right) ,\bar{z}\left( t\right) ,\bar{%
%u}\left( t\right) \right) $. Then, the following maximum principle holds:%
%\begin{eqnarray*}
%&&\Big <\chi \left( t\right) ,H\left( t,\bar{x}\left( t\right) ,\bar{y}%
%\left( t\right) ,\bar{z}\left( t\right) ,u,p\left( t\right) ,q\left(
%t\right) ,k\left( t\right) ,P\left( t\right) \right) \\
%&&-H\left( t,\bar{x}\left( t\right) ,\bar{y}\left( t\right) ,\bar{z}\left(
%t\right) ,\bar{u}\left( t\right) ,p\left( t\right) ,q\left( t\right)
%,k\left( t\right) ,P\left( t\right) \right) \Big > \\
%&\geq &0,\text{ }\forall u\in U,\text{ a.s., a.e.,}
%\end{eqnarray*}%
%where $\left( \chi \left( \cdot \right) ,k\left( \cdot \right) \right) $ is
%the solution to (\ref{adj}) and $H\left( \cdot \right) $ is defined in (\ref%
%{H}).
%\end{theorem}

As observed in Theorem \ref{th1}, there appears a term $\mathbb{H}\left(
s\right) x_{1}\left( s\right) v\left( s\right) $. In order to deal with it,
we give the expression of $x_{1}\left( \cdot \right) ,$ mainly taken from
Theorem 1.6.14 in Yong and Zhou \cite{YZ}. To this end, consider the
following matrix-valued stochastic differential equation:%
\begin{equation}
\left\{
\begin{array}{rcl}
\mathrm{d}\Psi \left( t\right) & = & b_{x}\left( t\right) \Psi \left(
t\right) \mathrm{d}t+\sigma _{x}\left( t\right) \Psi \left( t\right) \mathrm{%
d}W\left( t\right) , \\
\Psi \left( 0\right) & = & I,\text{ }t\in \left[ 0,T\right] ,%
\end{array}%
\text{ }\right.  \label{matrix}
\end{equation}%
where $I$ denotes the identity matrix in $\mathbb{R}^{n\times n}$. Then,
\begin{equation}
x_{1}\left( t\right) =\Psi \left( t\right) \left[ \int_{0}^{t}\Psi
^{-1}\left( r\right) \left( b_{u}\left( r\right) -\sigma _{x}\left( r\right)
\sigma _{u}\left( r\right) \right) v\left( r\right) \mathrm{d}%
r+\int_{0}^{t}\Phi ^{-1}\left( r\right) \sigma _{u}\left( r\right) v\left(
r\right) \mathrm{d}W\left( r\right) \right] .  \label{x1}
\end{equation}%
Substituting the explicit representation (\ref{x1}) of $x_{1}$ into (\ref%
{scm}) yields%
\begin{eqnarray}
0 &\leq &\mathbb{E}\Bigg [\int_{0}^{T}\chi \left( s\right) \Big \{%
f_{z}\left( s\right) p\left( s\right) \sigma _{uu}\left( s\right)
v^{2}\left( s\right)  \notag \\
&&+\mathbb{H}\left( s\right) v\left( s\right) \Psi \left( s\right) \Big [%
\int_{0}^{s}\Psi ^{-1}\left( r\right) \left( b_{u}\left( r\right) -\sigma
_{x}\left( r\right) \sigma _{u}\left( r\right) \right) v\left( r\right)
\mathrm{d}r  \notag \\
&&+\int_{0}^{s}\Psi ^{-1}\left( r\right) \sigma _{u}\left( r\right) v\left(
r\right) \mathrm{d}W\left( r\right) \Big ]+v^{2}\left( s\right) f_{zu}\left(
s\right) p\left( t\right) \sigma _{u}\left( t\right) \Big \}\mathrm{d}s\Bigg
].  \label{que}
\end{eqnarray}%
Clearly, (\ref{que}) contains an It\^{o}'s integral. Next we shall borrow
the spike variation method from \cite{ZZconvex} to check its order with
perturbed control. More precisely, let $\varepsilon >0$ and $E_{\varepsilon
}\subset \left[ 0,T\right] $ be a Borel set with Borel measure $\left\vert
E_{\varepsilon }\right\vert =\varepsilon ,$ define
\begin{equation*}
u^{\varepsilon }\left( t\right) =\bar{u}\left( t\right) I_{E_{\varepsilon
}^{c}}\left( t\right) +u\left( t\right) I_{E_{\varepsilon }}\left( t\right) ,
\end{equation*}%
where $u\left( \cdot \right) \in \mathcal{U}\left( 0,T\right) .$ This $%
u^{\varepsilon }$ is called a spike variation of the optimal control $\bar{u}
$. For our aim, we only need to use $E_{\varepsilon }=\left[ l,l+\varepsilon %
\right] $ for $l\in \left[ 0,T-\varepsilon \right] $ and $\varepsilon >0$.
Let
\begin{equation*}
v\left( \cdot \right) =u^{\varepsilon }\left( \cdot \right) -\bar{u}\left(
\cdot \right) =\left( u\left( t\right) -\bar{u}\left( \cdot \right) \right)
I_{E_{\varepsilon }}\left( \cdot \right) .
\end{equation*}%
Then, inserting it into (\ref{que}), we have

\begin{equation*}
\int_{l}^{l+\varepsilon }\chi \left( s\right) \mathbb{H}\left( s\right)
\left( u^{\varepsilon }\left( s\right) -\bar{u}\left( s\right) \right) \Psi
\left( s\right) \int_{l}^{s}\Psi ^{-1}\left( r\right) \sigma _{u}\left(
r\right) \left( u^{\varepsilon }\left( r\right) -\bar{u}\left( r\right)
\right) \mathrm{d}W\left( r\right) \mathrm{d}s
\end{equation*}%
By H\"{o}lder inequality and Burkholder-Davis-Gundy inequality, we have%
\begin{eqnarray*}
&&\mathbb{E}\left[ \int_{l}^{l+\varepsilon }\chi \left( s\right) \mathbb{H}%
\left( s\right) v\left( s\right) \Psi \left( s\right) \int_{0}^{s}\Psi
^{-1}\left( r\right) \sigma _{u}\left( r\right) v\left( r\right) \mathrm{d}%
W\left( r\right) \mathrm{d}s\right] \\
&\leq &\left[ \mathbb{E}\int_{l}^{l+\varepsilon }\left\vert \chi \left(
s\right) \mathbb{H}\left( s\right) v\left( s\right) \Psi \left( s\right)
\right\vert ^{2}\mathrm{d}s\right] ^{\frac{1}{2}}\cdot \\
&&\left[ \mathbb{E}\int_{l}^{l+\varepsilon }\int_{l}^{s}\left\vert \Psi
^{-1}\left( r\right) \sigma _{u}\left( r\right) \left( u^{\varepsilon
}\left( r\right) -\bar{u}\left( r\right) \right) \right\vert ^{2}\mathrm{d}r%
\mathrm{d}s\right] ^{\frac{1}{2}} \\
&\leq &C\varepsilon ^{\frac{3}{2}},
\end{eqnarray*}%
since $\sup_{s\in \left[ 0,T\right] }\left\vert \chi \left( s\right)
\right\vert ^{2}<\infty $ from classical estimate for stochastic
differential equations.

\begin{lemma}[Martingale representation theorem]
\label{marting}Suppose that $\phi \in L_{\mathbb{F}}^{2}\left( \Omega
;L^{2}\left( \left[ 0,T\right] :\mathbb{R}^{n}\right) \right) .$ Then, there
exists a $\kappa \left( \cdot ,\cdot \right) \in L^{2}\left( \left[ 0,T%
\right] ;L_{\mathbb{F}}^{2}\left( \left[ 0,T\right] \times \Omega ;\mathbb{R}%
^{n}\right) \right) $ such that%
\begin{equation*}
\phi \left( t\right) =\mathbb{E}\left[ \phi \left( t\right) \right]
+\int_{0}^{t}\kappa \left( s,t\right) \mathrm{d}W\left( s\right) ,\text{
a.s., a.e., }t\in \left[ 0,T\right] .
\end{equation*}
\end{lemma}

The proof can seen in Zhang et al. \cite{ZZconvex}.

\begin{lemma}
Assume that $\emph{(A1)}$\emph{-}$\emph{(A2)}$ hold. Then,
\begin{equation*}
\chi \left( \cdot \right) \mathbb{H}\left( \cdot \right) \in L_{\mathbb{F}%
}^{4}\left( \Omega ;L^{2}\left( \left[ 0,T\right] :\mathbb{R}\right) \right)
.
\end{equation*}
\end{lemma}

\paragraph{Proof.}

We shall prove that%
\begin{equation*}
\mathbb{E}\left[ \int_{0}^{T}\left\vert \chi \left( s\right) \mathbb{H}%
\left( t\right) \right\vert ^{2}\mathrm{d}t\right] ^{2}<\infty .
\end{equation*}%
From (A1)-(A2), we have
\begin{equation*}
\left\vert \psi _{x}\right\vert \leq C\text{ and }\left\vert \psi
_{xu}\right\vert \leq C\text{ for }\psi =b,\text{ }\sigma ,\text{ }f.
\end{equation*}
Besides,
\begin{equation*}
\left\vert f_{yu}\right\vert \leq C,\left\vert f_{zu}\right\vert \leq C.
\end{equation*}%
Hence,
\begin{eqnarray*}
&&\mathbb{E}\left[ \int_{0}^{T}\left\vert \chi \left( t\right) \mathbb{H}%
\left( t\right) \right\vert ^{2}\mathrm{d}t\right] ^{2}\leq \mathbb{E}\left[
\sup_{0\leq t\leq T}\left\vert \chi \left( t\right) \right\vert
^{2}\int_{0}^{T}\left\vert \mathbb{H}\left( t\right) \right\vert ^{2}\mathrm{%
d}t\right] ^{2} \\
&\leq &\frac{1}{2}\mathbb{E}\left[ \left( \sup_{0\leq t\leq T}\left\vert
\chi \left( t\right) \right\vert ^{2}\right) ^{4}+\left(
\int_{0}^{T}\left\vert \mathbb{H}\left( t\right) \right\vert ^{2}\mathrm{d}%
t\right) ^{4}\right] \\
&\leq &\frac{1}{2}\mathbb{E}\left[ \sup_{0\leq t\leq T}\left\vert \chi
\left( t\right) \right\vert ^{2}\right] ^{4}+\frac{1}{2}\mathbb{E}\Bigg [%
\int_{0}^{T}\Big |H_{xu}\left( t,\bar{x}\left( t\right) ,\bar{y}\left(
t\right) ,\bar{u}\left( t\right) ,\bar{z}\left( t\right) ,\bar{u}\left(
t\right) ,\bar{p}\left( t\right) ,\bar{q}\left( t\right) \right) \\
&&+Q\sigma _{u}\left( t,\bar{x}\left( t\right) ,\bar{u}\left( t\right)
\right) +P\sigma _{x}\left( t,\bar{x}\left( t\right) ,\bar{u}\left( t\right)
\right) \sigma _{u}\left( t,\bar{x}\left( t\right) ,\bar{u}\left( t\right)
\right) \\
&&+b_{u}\left( t,\bar{x}\left( t\right) ,\bar{u}\left( t\right) \right)
P\left( t\right) +f_{yu}\left( t,\bar{x}\left( t\right) ,\bar{y}\left(
t\right) ,\bar{u}\left( t\right) ,\bar{z}\left( t\right) ,\bar{u}\left(
t\right) \right) p\left( t\right) \\
&&+f_{zu}\left( t,\bar{x}\left( t\right) ,\bar{y}\left( t\right) ,\bar{u}%
\left( t\right) ,\bar{z}\left( t\right) ,\bar{u}\left( t\right) \right) %
\left[ p\sigma _{x}\left( t,\bar{x}\left( t\right) ,\bar{u}\left( t\right)
\right) +q\left( t\right) \right] \Big |^{2}\mathrm{d}t\Bigg ]^{4} \\
&\leq &\frac{1}{2}\mathbb{E}\left[ \sup_{0\leq t\leq T}\left\vert \chi
\left( t\right) \right\vert ^{2}\right] ^{4}+C+C\mathbb{E}\Bigg [%
\int_{0}^{T}\left( \left\vert p\left( t\right) \right\vert ^{2}+\left\vert
q\left( t\right) \right\vert ^{2}+\left\vert P\left( t\right) \right\vert
^{2}+\left\vert Q\left( t\right) \right\vert ^{2}\right) \mathrm{d}t\Bigg ]%
^{4} \\
&<&\infty .
\end{eqnarray*}%
From Lemma \ref{keylemma} and the classical estimation in (\ref{estadj}), we
finish the proof. \hfill $\Box $

Therefore, by our assumption (A1)-(A2) and Lemma \ref{marting}, for any $%
u\in U$, there exists a
\begin{equation*}
\psi ^{u}\left( \cdot ,\cdot \right) \in L^{2}\left( \left[ 0,T\right] :L_{%
\mathbb{F}}^{2}\left( \Omega \times \left[ 0,T\right] :\mathbb{R}^{n}\right)
\right)
\end{equation*}
such that for a.e. $t\in \left[ 0,T\right] $%
\begin{equation}
\chi \left( t\right) \mathbb{H}^{\top }\left( t\right) \left( u-\bar{u}%
\left( t\right) \right) =\mathbb{E}\left[ \chi \left( t\right) \mathbb{H}%
^{\top }\left( t\right) \left( u-\bar{u}\left( t\right) \right) \right]
+\int_{0}^{t}\psi ^{u}\left( s,t\right) dW\left( s\right) .  \label{p1}
\end{equation}%
Using (\ref{p1}), we are able to assert the following:

\begin{theorem}
\label{m1}Suppose that \emph{(A1)-(A2)} are in force. Let $\bar{u}\left(
\cdot \right) $ be a singular optimal control in the classical sense, then
we have%
\begin{equation*}
\begin{array}{l}
\mathbb{E}\left\langle \chi \left( r\right) \mathbb{H}\left( r\right)
b_{u}\left( r\right) \left( u-\bar{u}\left( r\right) \right) ,u-\bar{u}%
\left( r\right) \right\rangle \\
+\partial _{r}^{+}\left( \chi \left( r\right) \mathbb{H}^{\top }\left(
r\right) \left( u-\bar{u}\left( r\right) \right) ,\sigma _{u}\left( r\right)
\left( u-\bar{u}\left( r\right) \right) \right) \geq 0,\text{ a.e., }r\in %
\left[ 0,T\right] ,%
\end{array}%
\end{equation*}%
where
\begin{eqnarray*}
&&\partial _{r}^{+}\left( \chi \left( r\right) \mathbb{H}^{\top }\left(
r\right) \left( u-\bar{u}\left( r\right) \right) ,\sigma _{u}\left( r\right)
\left( u-\bar{u}\left( r\right) \right) \right) \\
&=&2\underset{\alpha \rightarrow 0}{\limsup }\frac{1}{\alpha ^{2}}\mathbb{E}%
\int_{r}^{r+\alpha }\int_{r}^{t}\left\langle \psi ^{u}\left( s,t\right)
,\Psi \left( r\right) \Psi ^{-1}\left( s\right) \sigma _{u}\left( s\right)
\left( u-\bar{u}\left( s\right) \right) \right\rangle \mathrm{d}s\mathrm{d}t,
\end{eqnarray*}%
where $\psi ^{u}\left( s,t\right) $ is obtained by (\ref{p1}), and $\Psi $
is determined by (\ref{matrix}).
\end{theorem}

The proof is just to repeat the process in Theorem 3.10, \cite{ZZconvex}, so
we omit it.

Note that Theorem \ref{m1} is pointwise with respect to the time variable $t$
(but also the integral form). Now if each of $\chi \left( \cdot \right)
\mathbb{H}\left( \cdot \right) $ and $\bar{u}\left( \cdot \right) $ are
regular enough, then the function $\psi ^{u}\left( \cdot ,\cdot \right) $
admits an explicit representation.

Suppose the following:

\begin{description}
\item[(A3)] $\bar{u}\left( \cdot \right) \in \mathbb{L}_{2,\mathbb{F}%
}^{1,2}\left( \mathbb{R}^{k}\right) ,$ $\chi \left( \cdot \right) \mathbb{H}%
^{\top }\left( \cdot \right) \in \mathbb{L}_{2,\mathbb{F}}^{1,2}\left(
\mathbb{R}^{k\times n}\right) \cap L^{\infty }\left( \left[ 0,T\right]
\times \Omega ;\mathbb{R}^{k\times n}\right) .$
\end{description}

\begin{theorem}
\label{m2}Suppose that the Assumptions \emph{(A1)-(A3)} are in force. Let $%
\bar{u}\left( \cdot \right) $ be a singular optimal control in the classical
sense, then we have%
\begin{equation*}
\begin{array}{l}
\left\langle \chi \left( r\right) \mathbb{H}\left( r\right) b_{u}\left(
r\right) \left( u-\bar{u}\left( r\right) \right) ,u-\bar{u}\left( r\right)
\right\rangle +\left\langle \nabla \left( \chi \left( r\right) \mathbb{H}%
\left( r\right) \right) \sigma _{u}\left( r\right) \left( u-\bar{u}\left(
r\right) \right) ,u-\bar{u}\left( r\right) \right\rangle \\
-\left\langle \chi \left( r\right) \mathbb{H}\left( r\right) \sigma
_{u}\left( r\right) \left( u-\bar{u}\left( r\right) \right) ,\nabla \bar{u}%
\left( r\right) \right\rangle \geq 0,\text{ a.e., }r\in \left[ 0,T\right] ,%
\text{ }\forall u\in U,\text{ }P\text{-a.s..}%
\end{array}%
\end{equation*}
\end{theorem}

Observe that the expression (\ref{scm}) is similar to (3.17) in \cite%
{ZZconvex}. Therefore, the proof is repeated as in Theorem 3.13 in Zhang and
Zhang \cite{ZZconvex}.

\subsubsection{Example}

\label{example}

We provide a concrete example to illustrate our theoretical result (Theorem %
\ref{m2}) by looking at an example. If the FBSDEs considered in this paper
are linear, it is possible to implement our principles directly. For
convenience, we still adopt the notations introduced in Section \ref{sec3.2}.

\begin{example}
Consider the following FBSDEs with $n=1$ and $U=\left[ -1,1\right] .$%
\begin{equation}
\left\{
\begin{array}{rcl}
\mathrm{d}x\left( t\right) & = & \left[ u\left( t\right) \sin x\left(
t\right) +u^{3}\left( t\right) \right] \mathrm{d}t+x\left( t\right) u\left(
t\right) \mathrm{d}W\left( t\right) , \\
-\mathrm{d}y\left( t\right) & = & \left[ x^{3}\left( t\right) +x\left(
t\right) u\left( t\right) +y\left( t\right) +z\left( t\right) +u^{3}\left(
t\right) \right] \mathrm{d}t-z\left( t\right) \mathrm{d}W\left( t\right) ,
\\
x\left( 0\right) & = & 0,\text{ }y\left( T\right) =x\left( T\right) .%
\end{array}%
\right.  \label{exa}
\end{equation}%
One can easily get the solutions to (\ref{adjbsde}),
\begin{equation*}
\chi \left( t\right) =\exp \left\{ \frac{1}{2}t+W\left( t\right) \right\} >0,%
\text{ }0\leq t\leq T,\text{ }P\text{-a.s.}
\end{equation*}%
Set $\left( \bar{x}\left( t\right) ,\bar{y}\left( t\right) ,\bar{z}\left(
t\right) ,\bar{u}\left( t\right) \right) =\left( 0,0,0,0\right) .$ The
corresponding adjoint equations are (\ref{adj1}) and (\ref{adj2}), namely,
\begin{equation}
\left\{
\begin{array}{rcl}
-\mathrm{d}p\left( t\right) & = & \left[ p\left( t\right) +q\left( t\right) %
\right] \mathrm{d}t-q\left( t\right) \mathrm{d}W\left( t\right) , \\
p\left( T\right) & = & 1,%
\end{array}%
\right.  \label{a1}
\end{equation}%
and%
\begin{equation}
\left\{
\begin{array}{rcl}
-\mathrm{d}P\left( t\right) & = & \left[ P\left( t\right) +Q\left( t\right) %
\right] \mathrm{d}t-Q\left( t\right) \mathrm{d}W\left( t\right) , \\
P\left( T\right) & = & 0.%
\end{array}%
\right.  \label{a2}
\end{equation}%
We get immediately, the solutions to (\ref{a1}) and (\ref{a2}) are%
\begin{equation*}
\left\{
\begin{array}{l}
\left( p\left( t\right) ,q\left( t\right) \right) =\left( e^{T-t},0\right) ,%
\text{ }t\in \left[ 0,T\right] , \\
\left( P\left( t\right) ,Q\left( t\right) \right) =\left( 0,0\right) ,\text{
}t\in \left[ 0,T\right] ,%
\end{array}%
\right.
\end{equation*}%
respectively$.$ Hence, we have
\begin{eqnarray*}
\mathcal{H}_{u}\left( t,\bar{x}\left( t\right) ,\bar{y}\left( t\right) ,\bar{%
z}\left( t\right) ,\bar{u}\left( t\right) ,p\left( t\right) ,q\left(
t\right) ,P\left( t\right) \right) &\equiv &0, \\
\mathcal{H}_{uu}\left( t,\bar{x}\left( t\right) ,\bar{y}\left( t\right) ,%
\bar{z}\left( t\right) ,\bar{u}\left( t\right) ,p\left( t\right) ,q\left(
t\right) ,P\left( t\right) \right) &\equiv &0.
\end{eqnarray*}%
Therefore, $\bar{u}\left( t\right) =0$ is a singular control in the
classical sense. Moreover, we compute
\begin{equation*}
\nabla \bar{u}\left( t\right) =0,\text{ }\mathbb{H}\left( t\right) =p\left(
t\right) \cos \bar{x}\left( t\right) +1=e^{T-t}+1.
\end{equation*}%
Consequently, we get%
\begin{eqnarray*}
&&\left\langle \chi \left( r\right) \mathbb{H}\left( r\right) b_{u}\left(
r\right) \left( u-\bar{u}\left( r\right) \right) ,u-\bar{u}\left( r\right)
\right\rangle \\
&&+\left\langle \nabla \left( \chi \left( r\right) \mathbb{H}\left( r\right)
\right) \sigma _{u}\left( r\right) \left( u-\bar{u}\left( r\right) \right)
,u-\bar{u}\left( r\right) \right\rangle \\
&&-\left\langle \chi \left( r\right) \mathbb{H}\left( r\right) \sigma
_{u}\left( r\right) \left( u-\bar{u}\left( r\right) \right) ,\nabla \bar{u}%
\left( r\right) \right\rangle , \\
&=&0,\text{ }\forall u\in \left[ -1,1\right] ,\text{ a.e. }r\in \left[ 0,T%
\right] ,\text{ }P\text{-a.s.,}
\end{eqnarray*}%
which indicates that Theorem \ref{m2} always holds and $\bar{u}\left(
r\right) =0$ is a singular optimal control.
\end{example}

\section{Singular Optimal Controls via Dynamic Programming Principle}

\label{secverif}In this section, we proceed our control problem from the
view point of DPP. From now on, we focus on the following
\begin{equation}
\left\{
\begin{array}{rcl}
\mathrm{d}X_{s}^{t,x;v,\xi } & = & b\left( s,X_{s}^{t,x;v,\xi },v_{s}\right)
\mathrm{d}s+\sigma \left( s,X_{s}^{t,x;v,\xi },v_{s}\right) \mathrm{d}W_{s}+G%
\mathrm{d}\xi _{s}, \\
\mathrm{d}Y_{s}^{t,x;v,\xi } & = & -f\left( s,X_{s}^{t,x;v,\xi
},Y_{s}^{t,x;v,\xi },Z_{s}^{t,x;v,\xi },v_{s}\right) \mathrm{d}%
s+Z_{s}^{t,x;v,\xi }\mathrm{d}W_{s}-K\mathrm{d}\xi _{s}, \\
X_{t}^{t,x;v,\xi } & = & x,\text{ }Y_{T}^{t,x;v,\xi }=\Phi \left(
X_{T}^{t,x;v,\xi }\right) ,\qquad 0\leq t\leq s\leq T.%
\end{array}%
\right.  \label{dppfbsde}
\end{equation}%
Since the value function defined by the solution of controlled BSDE (\ref%
{dppfbsde}), so from the existence and uniqueness, $u$ defined in (\ref%
{value}) is well-defined.

\begin{remark}
\label{GK} We assume that $G_{n\times m}$ and $K_{1\times m}$ are
deterministic matrices. On the one hand, from the derivations in Theorem 5.1
of \cite{HS2}, it is convenient to show the \textquotedblleft
inaction\textquotedblright\ region for singular control; On the other hand,
we may regard $Y_{s}^{t,x;v,\xi }+K\xi _{s}$ together as a solution, in this
way, we are able to apply the classical It\^{o}'s formula, avoiding the
appearance of jump. We believe these assumptions can be removed properly,
but at present, we consider constant only in our paper. Whilst in order to
get the uniqueness of the solution to H-J-B inequality (\ref{pde1}), we add
the assumption $K^{i}>k_{0}>0,$ $1\leq i\leq m.$ More details, see Theorem
2.2 in \cite{Zhangsingular}.
\end{remark}

Set
\begin{eqnarray*}
\mathcal{L}\left( t,x,v\right) \varphi &=&\frac{1}{2}\text{\textrm{Tr}}%
\left( \sigma \sigma ^{\ast }\left( t,x,v\right) \varphi _{xx}\right)
+\left\langle \varphi _{x},b\left( t,x,v\right) \right\rangle , \\
\left( t,x,v\right) &\in &\left[ 0,T\right] \times \mathbb{R}^{n}\times U,%
\text{ }\varphi \in C^{1,2}\left( \left[ 0,T\right] \times \mathbb{R}%
^{n}\right) .
\end{eqnarray*}

\subsection{Verification Theorem via Viscosity Solutions}

\label{secveri}

Zhang \cite{Zhangsingular} has given a verification theorem for smooth
solution of the following H-J-B inequality:
\begin{equation}
\left\{
\begin{array}{l}
\min \Big (V_{x}^{\top }\left( t,x\right) G+K,\frac{\partial }{\partial t}%
V\left( t,x\right) \\
+\min_{v\in U}\mathcal{L}\left( t,x,v\right) V\left( t,x\right) +f\left(
t,x,V\left( t,x\right) ,V_{x}^{\top }\left( t,x\right) \sigma \left(
t,x,v\right) ,v\right) \Big )=0, \\
u\left( T,x\right) =\Phi \left( x\right) ,\text{ }0\leq t\leq T.%
\end{array}%
\right.  \label{pde1}
\end{equation}

\begin{lemma}
\label{P2}Define
\begin{equation*}
\mathcal{D}_{t}\left( V\right) :=\left\{ x\in \mathbb{R}^{n}:V\left(
t,x\right) <V\left( t,x+Gh\right) +Kh,\text{ }h\in \mathbb{R}_{+}^{m},\text{
}h\neq 0\right\} .
\end{equation*}%
Then the optimal state process $X^{t,x;\hat{v},\hat{\xi}}$ is continuous
whenever $\left( r,X_{r}^{t,x;\hat{v},\hat{\xi}}\right) \in \mathcal{D}%
_{r}\left( u\right) $. To be precise, we have%
\begin{equation*}
P\left( \Delta X_{r}^{t,x;\hat{v},\hat{\xi}}\neq 0,\text{ }X_{r}^{t,x;\hat{v}%
,\hat{\xi}}\in \mathcal{D}_{r}\left( V\right) \right) =0,\text{ }t\leq r\leq
T.
\end{equation*}
\end{lemma}

The proof can be seen in Zhang \cite{Zhangsingular}.

\begin{proposition}
\label{smooth} Suppose that $V$ is a classical solution of the H-J-B
inequality \emph{(\ref{pde1})} such that for some $l>1,$
\begin{equation*}
\left\vert V\left( t,x\right) \right\vert \leq C\left( 1+\left\vert
x\right\vert ^{l}\right) .
\end{equation*}
Then for any $\left[ 0,T\right] \times \mathbb{R}^{n}$, $\left( v,\xi
\right) \in \mathcal{U}:$%
\begin{equation*}
V\left( t,x\right) \leq J\left( t,x,v,\xi \right) .
\end{equation*}%
Furthermore, if there exists $\left( \hat{v},\hat{\xi}\right) \in \mathcal{U}
$ such that
\begin{eqnarray}
1 &=&P\left\{ \left( r,X_{r}^{t,x;\hat{v},\hat{\xi}}\right) \in \mathcal{D}%
_{r}\left( V\right) ,\text{ }0\leq r\leq T\right\} ,  \label{v11} \\
1 &=&P\left\{ \int_{\left[ t,T\right] }\left[ V_{x}^{\top }\left( r,x\right)
G+K\right] \mathrm{d}r=0\right\} ,  \label{v22} \\
1 &=&P\Bigg \{\left( s,X_{s+}^{t,x;\hat{v},\hat{\xi}}\right) \in \mathcal{D}%
_{r}\left( V\right) ,\text{ }t\leq s\leq T:  \notag \\
\hat{v}_{s} &\in &\min_{v\in U}\Big [V_{t}\left( s,X_{s}^{t,x;v,\xi }\right)
+\mathcal{L}\left( s,X_{s+}^{t,x;\hat{v},\hat{\xi}},v\right) V\left(
s,X_{s+}^{t,x;\hat{v},\hat{\xi}}\right)  \notag \\
&&+f\Big (s,X_{s+}^{t,x;\hat{v},\hat{\xi}},V\left( t,X_{s+}^{t,x;\hat{v},%
\hat{\xi}}\right) ,  \notag \\
&&\nabla V\left( s,X_{s+}^{t,x;\hat{v},\hat{\xi}}\right) \sigma \left(
s,X_{s+}^{t,x;\hat{v},\hat{\xi}},v\right) ,v\Big )\Big ]\Bigg \}  \label{v33}
\end{eqnarray}%
and%
\begin{equation}
P\left\{ V\left( s,X_{s}^{t,x;\hat{v},\hat{\xi}}\right) =V\left(
s,X_{s+}^{t,x;\hat{v},\hat{\xi}}\right) +K\Delta \hat{\xi}_{s},\text{ }t\leq
s\leq T\right\} =1.  \label{v44}
\end{equation}%
Then%
\begin{equation}
V\left( t,x\right) =J\left( t,x;\hat{v}\left( \cdot \right) ,\hat{\xi}\left(
\cdot \right) \right) .  \label{v55}
\end{equation}
\end{proposition}

In this section, we remove the \emph{unreal} condition, smooth on value
function, by means of viscosity solutions\footnote{%
In the classical optimal stochastic control theory, the value function is a
solution to the corresponding H-J-B equation whenever it has sufficient
regularity (Fleming and Rishel \cite{FR}, Krylov \cite{Kry}). Nevertheless,
when it is only known that the value function is continuous, then, the value
function is a solution to the H-J-B equation in the viscosity sense (see
Lions \cite{CIL}).}. We will recall the definition of a viscosity solution
for H-J-B variational inequality (\ref{pde1}) from \cite{CIL}. %
Below, $\mathbb{S}^{n}$ will denote the set of $n\times n$ symmetric
matrices.

Let us begin at introducing the following parabolic superjet:

%We first recall the definition of a viscosity solution
%for H-J-B variational inequality (\ref{pde1}) from \cite{CIL}.

%\begin{definition}
%\label{dvis1} Let $u\left( t,x\right) \in C\left( \left[ 0,T\right] \times
%\mathbb{R}^{n}\right) $ and $\left( t,x\right) \in \left[ 0,T\right] \times
%\mathbb{R}^{n}.$ For every $\varphi \in C^{1,2}\left( \left[ 0,T\right]
%\times \mathbb{R}^{n}\right) $
%
%\noindent (1) for each local maximum point $\left( t_{0},x_{0}\right) $ of $%
%u-\varphi $ in the interior of $\left[ 0,T\right] \times \mathbb{R}^{n},$ we
%have%
%\begin{equation}
%\min \left( \varphi _{x}^{\top }G+K,\frac{\partial \varphi }{\partial t}%
%+\min_{v\in U}\left\{ \mathcal{L}\varphi +f\left( t_{0},x_{0},\varphi
%,\nabla \varphi \sigma ,v\right) \right\} \right) \geq 0  \label{vis1}
%\end{equation}
%
%at $\left( t_{0},x_{0}\right) ,$ i.e., $u$ is a subsolution.
%
%\noindent (2) for each local mimimum point $\left( t_{0},x_{0}\right) $ of $%
%u-\varphi $ in the interior of $\left[ 0,T\right] \times \mathbb{R}^{n},$ we
%have%
%\begin{equation}
%\min \left( \varphi _{x}^{\top }G+K,\frac{\partial \varphi }{\partial t}%
%+\min_{v\in U}\left\{ \mathcal{L}\varphi +f\left( t_{0},x_{0},\varphi
%,\nabla \varphi \sigma ,v\right) \right\} \right) \leq 0  \label{vis2}
%\end{equation}
%
%at $\left( t_{0},x_{0}\right) ,$ i.e., $u$ is a supersolution.
%
%\noindent (3) $u\left( t,x\right) \in C\left( \left[ 0,T\right] \times
%\mathbb{R}^{n}\right) $ is said to be a viscosity solution of (\ref{pde1})
%if it is both a viscosity sub and supersolution.
%\end{definition}

\begin{definition}
\label{d1}Let $V\left( t,x\right) \in C\left( \left[ 0,T\right] \times
\mathbb{R}^{n}\right) $ and $\left( t,x\right) \in \left[ 0,T\right] \times
\mathbb{R}^{n}$. We denote by $\mathcal{P}^{2,+}V\left( t,x\right) $, the
\textquotedblleft parabolic superjet\textquotedblright\ of $V$ at $\left(
t,x\right) $ the set of triples $\left( p,q,X\right) \in \mathbb{R}\times
\mathbb{R}^{n}\times \mathbb{S}^{n}$ which are such that%
\begin{eqnarray*}
V\left( s,y\right) &\leq &V\left( t,x\right) +p\left( s-t\right)
+\left\langle q,x-y\right\rangle \\
&&+\frac{1}{2}\left\langle X\left( y-x\right) ,y-x\right\rangle +o\left(
\left\vert s-t\right\vert +\left\vert y-x\right\vert ^{2}\right) .
\end{eqnarray*}%
Similarly, we denote by $\mathcal{P}^{2,-}V\left( t,x\right) ,$ the
\textquotedblleft parabolic subjet\textquotedblright\ of $V$ at $\left(
t,x\right) $ the set of triples $\left( p,q,X\right) \in \mathbb{R}\times
\mathbb{R}^{n}\times \mathbb{S}^{n}$ which are such that%
\begin{eqnarray*}
V\left( s,y\right) &\geq &V\left( t,x\right) +p\left( s-t\right)
+\left\langle q,x-y\right\rangle \\
&&+\frac{1}{2}\left\langle X\left( y-x\right) ,y-x\right\rangle +o\left(
\left\vert s-t\right\vert +\left\vert y-x\right\vert ^{2}\right) .
\end{eqnarray*}
\end{definition}

\begin{lemma}
\label{text}Let $V\in C\left( \left[ 0,T\right] \times \mathbb{R}^{n}\right)
$ and $\left( t,x\right) \in \left[ 0,T\right] \times \mathbb{R}^{n}$ be
given. Then:

1) $\left( p,q,X\right) \in \mathcal{P}^{2,+}V\left( t,x\right) $ if and
only if there exists a function $\varphi \in C^{1,2}\left( \left[ 0,T\right]
\times \mathbb{R}^{n}\right) $ such that $V-\varphi $ attains a strict
maximum at $\left( t,x\right) $ and
\begin{equation*}
\left( \varphi \left( t,x\right) ,\varphi _{t}\left( t,x\right) ,\varphi
_{x}\left( t,x\right) ,\varphi _{xx}\left( t,x\right) \right) =\left(
V\left( t,x\right) ,p,q,X\right) .
\end{equation*}%
2) $\left( p,q,X\right) \in \mathcal{P}^{2,-}V\left( t,x\right) $ if and
only if there exists a function $\varphi \in C^{1,2}\left( \left[ 0,T\right]
\times \mathbb{R}^{n}\right) $ such that $V-\varphi $ attains a strict
minimum at $\left( t,x\right) $ and
\begin{equation*}
\left( \varphi \left( t,x\right) ,\varphi _{t}\left( t,x\right) ,\varphi
_{x}\left( t,x\right) ,\varphi _{xx}\left( t,x\right) \right) =\left(
V\left( t,x\right) ,p,q,X\right) .
\end{equation*}
\end{lemma}

More details can be seen in Lemma 5.4 and 5.5 in Yong and Zhou \cite{YZ}.

Define%
\begin{equation}
\mathcal{G}\left( t,x,q,X\right) =\left\{ \frac{1}{2}\text{\textrm{Tr}}%
\left( \sigma \sigma ^{\ast }\left( t,x,v\right) X\right) +\left\langle
q,b\left( t,x,v\right) \right\rangle +f\left( t,x,V\left( t,x\right)
,q^{\top }\sigma \left( t,x,v\right) \right) \right\} .  \label{Hamvis}
\end{equation}

\begin{definition}
\label{dvis2} (i) It can be said $V\left( t,x\right) \in C\left( \left[ 0,T%
\right] \times \mathbb{R}^{n}\right) $ is a viscosity subsolution of \emph{(%
\ref{pde1})} if $V\left( T,x\right) \geq \Phi \left( x\right) ,$ $x\in
\mathbb{R}^{n}$, and at any point $\left( t,x\right) \in \left[ 0,T\right]
\times \mathbb{R}^{n}$, for any $\left( p,q,X\right) \in \mathcal{P}%
^{2,+}V\left( t,x\right) $,%
\begin{equation}
\min \left( qG+K,p+\mathcal{G}\left( t,x,q,X\right) \right) \geq 0.
\label{opde1}
\end{equation}%
In other words, at any point $\left( t,x\right) ,$ we have both $qG+K\geq 0$
and
\begin{equation*}
p+\mathcal{G}\left( t,x,q,X\right) \geq 0.
\end{equation*}%
\noindent (ii) It can be said $V\left( t,x\right) \in C\left( \left[ 0,T%
\right] \times \mathbb{R}^{n}\right) $ is a viscosity supersolution of \emph{%
(\ref{pde1})} if $V\left( T,x\right) \leq \Phi \left( x\right) ,$ $x\in
\mathbb{R}^{n}$, and at any point $\left( t,x\right) \in \left[ 0,T\right]
\times \mathbb{R}^{n}$, for any $\left( p,q,X\right) \in \mathcal{P}%
^{2,-}V\left( t,x\right) $,%
\begin{equation}
\min \left( qG+K,p+\mathcal{G}\left( t,x,q,X\right) \right) \leq 0.
\end{equation}%
In other words, at any point where $qG+K\geq 0$, we have
\begin{equation*}
p+\mathcal{G}\left( t,x,q,X\right) \leq 0.
\end{equation*}%
\noindent (iii) It can be said $V\left( t,x\right) \in C\left( \left[ 0,T%
\right] \times \mathbb{R}^{n}\right) $ is a viscosity solution of \emph{(\ref%
{pde1})} if it is both a viscosity sub and super solution.
\end{definition}

We have the following result:

\begin{proposition}
\label{t4}Assume that \emph{(A1)-(A2)} are in force. Then there exists at
most one viscosity solution of H-J-B inequality \emph{(\ref{pde1})} in the
class of bounded and continuous functions.
\end{proposition}

%\begin{remark}
%We have put somewhat strong assumptions, namely, $b,$ $\sigma ,$ $f$ are
%bounded. These conditions may be removed by modifying the idea by Ishii \cite%
%{Is}.
%\end{remark}

We need a generalized It\^{o}'s formula. Define
\begin{eqnarray*}
\mathcal{L}\left( t,x,v\right) \Psi &=&\frac{1}{2}\text{\textrm{Tr}}\left(
\sigma \sigma ^{\ast }\left( t,x,v\right) D^{2}\Psi \right) +\left\langle
D\Psi ,b\left( t,x,v\right) \right\rangle , \\
\left( t,x,v\right) &\in &\left[ 0,T\right] \times \mathbb{R}^{n}\times U,%
\text{ }\Psi \in C^{1,2}\left( \left[ 0,T\right] \times \mathbb{R}%
^{n}\right) .
\end{eqnarray*}%
For any $\Psi \in C^{1,2}\left( \left[ 0,T\right] \times \mathbb{R}^{n};%
\mathbb{R}\right) $, by virtue of Dol\'{e}ans-Dade-Meyer formula (see \cite%
{HS2,CH}), we have%
\begin{eqnarray}
\Psi \left( s,X_{s}\right) &=&\Psi \left( t,x\right) +\int_{t}^{s}\Psi
_{t}\left( r,X_{r}\right) +\mathcal{L}\left( r,X_{r},v\right) \Psi \left(
r,X_{r}\right) \mathrm{d}r  \notag \\
&&+\int_{t}^{s}\Psi _{x}\left( r,X_{r}\right) \sigma \left(
r,X_{r},v_{r}\right) \mathrm{d}W_{r}+\int_{t}^{s}\Psi _{x}\left(
r,X_{r}\right) G\mathrm{d}\xi _{r}  \notag \\
&&+\sum_{t\leq r\leq s}\left\{ \Psi \left( r,X_{r+}\right) -\Psi \left(
r,X_{r}\right) -\Psi _{x}\left( r,X_{r}\right) \Delta X_{r}\right\} .
\label{ito1}
\end{eqnarray}

We begin to introduce a useful lemma.

\begin{lemma}
\label{l4}Assume that \emph{(A1)-(A2)} are in force. Let $\left( t,x\right)
\in \left[ 0,T\right) \times \mathbb{R}^{n}$ be fixed and let $\left(
X^{t,x;u}\left( \cdot \right) ,u\left( \cdot \right) \right) $ be an
admissible pair. Define processes
\begin{equation*}
\left\{
\begin{array}{l}
z_{1}\left( r\right) \doteq b\left( r,X^{t,x;u}\left( r\right) ,u\left(
r\right) \right) , \\
z_{2}\left( r\right) \doteq \sigma \left( r,X^{t,x;u}\left( r\right)
,u\left( r\right) \right) \sigma ^{\ast }\left( r,X^{t,x;u}\left( r\right)
,u\left( r\right) \right) , \\
z_{3}\left( r\right) \doteq f\left( r,X^{t,x;u}\left( r\right)
,Y^{t,x;u}\left( r\right) ,Z^{t,x;u}\left( r\right) ,u\left( r\right)
\right) .%
\end{array}%
\right.
\end{equation*}%
Then
\begin{equation}
\lim\limits_{h\rightarrow 0+}\frac{1}{h}\int_{t}^{t+h}\left\vert z_{i}\left(
r\right) -z_{i}\left( t\right) \right\vert \mbox{\rm
d}r=0,\quad \text{a.e. }t\in \left[ 0,T\right] ,\text{ }i=1,2,3.  \label{3.1}
\end{equation}
\end{lemma}

The proof can be found in \cite{YZ}.

\begin{lemma}
\label{l5}Let $g\in C\left( \left[ 0,T\right] \right) .$ Extend $g$ to $%
\left( -\infty ,+\infty \right) $ with $g\left( t\right) =g\left( T\right) $
for $t>T,$ and $g\left( t\right) =g\left( 0\right) ,$ for $t<0.$ Suppose
that there is a integrable function $\rho \in L^{1}\left( \left[ 0,T\right] ;%
\mathbb{R}\right) $ and some $h_{0}>0,$ such that
\begin{equation*}
\frac{g\left( t+h\right) -g\left( t\right) }{h}\leq \rho \left( t\right)
,\quad \text{a.e. }t\in \left[ 0,T\right] ,\qquad h\leq h_{0}.
\end{equation*}%
Then
\begin{equation*}
g\left( \beta \right) -g\left( \alpha \right) \leq \int_{\alpha }^{\beta
}\limsup_{h\rightarrow 0+}\frac{g\left( t+h\right) -g\left( t\right) }{h}%
\mbox{\rm
d}r,\text{ }\forall 0\leq \alpha \leq \beta \leq T.
\end{equation*}
\end{lemma}

The proof can be seen in Zhang \cite{zlqverif}.

The main result in this section is the following.%\paragraph{Proof}
%
%Applying Fatou's Lemma, we have
%\begin{eqnarray*}
%\int_{\alpha }^{\beta }\rho \left( r\right) \text{d}r &\geq &\int_{\alpha
%}^{\beta }\limsup_{h\rightarrow 0+}\frac{g\left( r+h\right) -g\left(
%r\right) }{h}\mbox{\rm
%d}r \\
%&\geq &\limsup_{h\rightarrow 0+}\int_{\alpha }^{\beta }\frac{g\left(
%r+h\right) -g\left( r\right) }{h}\mbox{\rm
%d}r \\
%&=&\limsup_{h\rightarrow 0+}\frac{\int_{\alpha +h}^{\beta +h}g\left(
%r\right) \text{d}r-\int_{\alpha }^{\beta }g\left( r\right) \mbox{\rm
%d}r}{h} \\
%&=&\limsup_{h\rightarrow 0+}\frac{\int_{\beta }^{\beta +h}g\left( r\right)
%\text{d}r-\int_{\alpha }^{\alpha +h}g\left( r\right) \mbox{\rm
%d}r}{h} \\
%&=&g\left( \beta \right) -g\left( \alpha \right) .
%\end{eqnarray*}%
%We complete the proof. \hfill $\Box $

\begin{theorem}[Verification Theorem]
\label{verif}Suppose that the Assumptions \emph{(A1)-(A2)} are in force. Let
$V\in C\left( \left[ 0,T\right] \times \mathbb{R}^{n}\right) ,$ be a
viscosity solution of the H-J-B equations (\ref{pde1}), satisfying the
following conditions:
\begin{equation}
\left\{
\begin{array}{l}
\text{i) }V\left( t+h,x\right) -V\left( t,x\right) \leq C\left( 1+\left\vert
x\right\vert ^{m}\right) h,\qquad m\geq 0, \\
\quad \text{for all }x\in \mathbb{R}^{n},0<t<t+h<T. \\
\text{ii) }V\text{ is semiconcave},\text{ uniformly in }t,\text{i.e}.\text{
there exists }C_{0}\geq 0\text{ } \\
\quad \text{such that for every }t\in \left[ 0,T\right] ,\text{ }V\left(
t,\cdot \right) -C_{0}\left\vert \cdot \right\vert ^{2}\text{is concave on }%
\mathbb{R}^{n}.%
\end{array}%
\right.  \label{3.2}
\end{equation}%
Then we have
\begin{equation}
V\left( t,x\right) \leq J\left( t,x;u\left( \cdot \right) ,\xi \left( \cdot
\right) \right) ,\text{ }  \label{3.3}
\end{equation}%
for any $\left( t,x\right) \in \left( 0,T\right] \times \mathbb{R}^{n}$ and
any $u\left( \cdot \right) \times \xi \left( \cdot \right) \in \mathcal{U}%
\left( t,T\right) .$

Furthermore, let $\left( t,x\right) \in \left( 0,T\right] \times \mathbb{R}%
^{n}$ be fixed and let
\begin{equation*}
\left( \bar{X}^{t,x;\bar{u},\bar{\xi}}\left( \cdot \right) ,\bar{Y}^{t,x;%
\bar{u},\bar{\xi}}\left( \cdot \right) ,\bar{Z}^{t,x;\bar{u},\bar{\xi}%
}\left( \cdot \right) \bar{u}\left( \cdot \right) ,\bar{\xi}\left( \cdot
\right) \right)
\end{equation*}%
be an admissible pair such that there exist a function $\varphi \in
C^{1,2}\left( \left[ 0,T\right] ;\mathbb{R}^{n}\right) $ and a triple
\begin{equation}
\left( \overline{p},\overline{q},\overline{\Theta }\right) \in \left( L_{%
\mathcal{F}_{t}}^{2}\left( \left[ t,T\right] ;\mathbb{R}\right) \times L_{%
\mathcal{F}_{t}}^{2}\left( \left[ t,T\right] ;\mathbb{R}^{n}\right) \times
L_{\mathcal{F}_{t}}^{2}\left( \left[ t,T\right] ;\mathbf{S}^{n}\right)
\right)  \label{3.4}
\end{equation}%
satisfying
\begin{equation}
\left\{
\begin{array}{l}
\left( \bar{p}\left( s\right) ,\bar{q}\left( s\right) ,\bar{\Theta}\left(
s\right) \right) \in \mathcal{P}^{2,+}V\left( s,\bar{X}^{t,x;\bar{u},\bar{\xi%
}}\left( s\right) \right) , \\
\left( \frac{\partial \varphi }{\partial t}\left( s,\bar{X}^{t,x;\bar{u},%
\bar{\xi}}\left( s\right) \right) ,D_{x}\varphi \left( s,\bar{X}^{t,x;\bar{u}%
,\bar{\xi}}\left( s\right) \right) ,D^{2}\varphi \left( s,\bar{X}^{t,x;\bar{u%
},\bar{\xi}}\left( s\right) \right) \right) =\left( \overline{p}\left(
s\right) ,\overline{q}\left( s\right) ,\overline{\Theta }\left( s\right)
\right) , \\
\bar{p}\left( t\right) G+K=0,\text{ a.e. }t\in \left[ 0,T\right] ,\text{ }P%
\text{-a.s.}%
\end{array}%
\right.  \label{3.5}
\end{equation}%
and
\begin{equation}
\mathbb{E}\left[ \int_{t}^{T}\left[ \overline{p}\left( s\right) +\mathcal{G}%
\left( s,\bar{X}^{t,x;\bar{u},\bar{\xi}}\left( s\right) ,\overline{\varphi }%
\left( s\right) ,\overline{p}\left( s\right) ,\overline{\Theta }\left(
s\right) ,\overline{u}\left( s\right) \right) \right] \mbox{\rm
d}s\right] \leq 0,  \label{3.6}
\end{equation}%
where $\overline{\varphi }\left( s\right) =\varphi \left( s,\bar{X}^{t,x;%
\bar{u},\bar{\xi}}\left( s\right) \right) $ and $\mathcal{G}$ is defined in (%
\ref{Hamvis}). Then $(\bar{X}^{t,x;\bar{u},\bar{\xi}}\left( \cdot \right) ,%
\bar{u}\left( \cdot \right) ,\bar{\xi}\left( \cdot \right) )$ is an optimal
pair$.$
\end{theorem}

In order to prove Theorem \ref{verif}, we need the following lemma:

\begin{lemma}
\label{l6}Let $v$ be a viscosity subsolution of the H-J-B equations (\ref%
{pde1}) satisfying (\ref{3.2}). Then we have%
\begin{equation}
\mathbb{E}\frac{1}{h}\left[ v\left( s+h,\bar{X}^{t,x;\bar{u},\bar{\xi}%
}\left( s+h\right) \right) -v\left( s,\bar{X}^{t,x;\bar{u},\bar{\xi}}\left(
s\right) \right) \right] \leq \rho \left( s\right) ,\text{ }0<h\leq T-t,
\label{keyverif}
\end{equation}%
where $\rho \left( s\right) \in L^{1}\left( \left[ t,T\right] :\mathbb{R}%
\right) .$
\end{lemma}

The proof can be seen in the Appendix.

\paragraph{Proof of Theorem \protect\ref{verif}.}

We have (\ref{3.3}) from the uniqueness of viscosity solutions of the H-J-B
equations (\ref{pde1}). It remains to show that $\left( \bar{X}^{t,x;\bar{u},%
\bar{\xi}}\left( \cdot \right) ,\bar{u}\left( \cdot \right) ,\bar{\xi}\left(
\cdot \right) \right) $ is an optimal, we now fix $t_{0}\in \left[ t,T\right]
$ such that (\ref{3.4}) and (\ref{3.5}) hold at $t_{0}.$ For $z_{1}\left(
\cdot \right) =\overline{b}\left( \cdot \right) ,$ $z_{2}\left( \cdot
\right) =\overline{\sigma }\left( \cdot \right) \overline{\sigma }\left(
\cdot \right) ^{\ast },$ $z_{3}\left( \cdot \right) =\overline{f}\left(
\cdot \right) .$ We claim that the set of such points is of full measure in $%
\left[ t,T\right] $ by Lemma 7 in \cite{zlqverif}. Now we fix $\omega
_{0}\in \Omega $ such that the regular conditional probability $P\left(
\left. \cdot \right\vert \mathcal{F}_{t_{0}}^{t}\right) \left( \omega
_{0}\right) $, given $\mathcal{F}_{t_{0}}^{t}$ is well defined. In this new
probability space, the random variables $\bar{X}^{t,x;\bar{u},\bar{\xi}%
}\left( t_{0}\right) ,\overline{p}\left( t_{0}\right) ,\overline{q}\left(
t_{0}\right) ,\overline{\Theta }\left( t_{0}\right) $ are almost surely
deterministic constants and equal to
\begin{equation*}
\bar{X}^{t,x;\bar{u},\bar{\xi}}\left( t_{0},\omega _{0}\right) ,\overline{p}%
\left( t_{0},\omega _{0}\right) ,\overline{q}\left( t_{0},\omega _{0}\right)
,\overline{\Theta }\left( t_{0},\omega _{0}\right) ,
\end{equation*}%
respectively. We remark that in this probability space the Brownian motion $%
W $ is still the a standard Brownian motion although now $W\left(
t_{0}\right) =W\left( t_{0},\omega _{0}\right) $ almost surely. The space is
now equipped with a new filtration $\left\{ \mathcal{F}_{r}^{t}\right\}
_{t\leq r\leq T}$ and the control process $\overline{u}\left( \cdot \right) $
is adapted to this new filtration. For $P$-a.s. $\omega _{0}$ the process $%
\bar{X}^{t,x;\bar{u},\bar{\xi}}\left( \cdot \right) $ is a solution of (1.1)
on $\left[ t_{0},T\right] $ in $\left( \Omega ,\mathcal{F},P\left( \left.
\cdot \right\vert \mathcal{F}_{t_{0}}^{t}\right) \left( \omega _{0}\right)
\right) $ with the inial condition $\bar{X}^{t,x;\bar{u},\bar{\xi}}\left(
t_{0}\right) =\bar{X}^{t,x;\bar{u},\bar{\xi}}\left( t_{0},\omega _{0}\right)
.$

Then on the probability space $\left( \Omega ,\mathcal{F},P\left( \left.
\cdot \right\vert \mathcal{F}_{t_{0}}^{t}\right) \left( \omega _{0}\right)
\right) $, we are going to apply It\^{o}'s formula to $\varphi $ on $\left[
t_{0},t_{0}+h\right] $ for any $h>0,$%
\begin{eqnarray*}
&&\ \varphi \left( t_{0}+h,\bar{X}^{t,x;\bar{u},\bar{\xi}}\left(
t_{0}+h\right) \right) -\varphi \left( t_{0},\bar{X}^{t,x;\bar{u},\bar{\xi}%
}\left( t_{0}\right) \right) \\
\ &=&\int_{t_{0}}^{t_{0}+h}\Big [\frac{\partial \varphi }{\partial t}\left(
r,\bar{X}^{t,x;\bar{u},\bar{\xi}}\left( r\right) \right) +\left\langle
D_{x}\varphi \left( r,\bar{X}^{t,x;\bar{u},\bar{\xi}}\left( r\right) \right)
,\overline{b}\left( r\right) \right\rangle \\
&&\ +\frac{1}{2}\text{tr}\left\{ \overline{\sigma }\left( r\right) ^{\ast
}D_{xx}\varphi \left( r,\bar{X}^{t,x;\bar{u},\bar{\xi}}\left( r\right)
\right) \overline{\sigma }\left( r\right) \right\} \Big ]\mbox{\rm
d}r \\
&&+\int_{t_{0}}^{t_{0}+h}\varphi _{x}\left( r,\bar{X}^{t,x;\bar{u},\bar{\xi}%
}\left( r\right) \right) G\mathrm{d}\bar{\xi}_{r}+\int_{t_{0}}^{t_{0}+h}%
\left\langle D_{x}\varphi \left( r,\bar{X}^{t,x;\bar{u},\bar{\xi}}\left(
r\right) \right) ,\overline{\sigma }\left( r\right) \right\rangle
\mbox{\rm
d}W_{r} \\
&&+\sum_{t_{0}\leq r\leq t_{0}+h}\Big \{\varphi \left( r,\bar{X}^{t,x;\bar{u}%
,\bar{\xi}}\left( r\right) \right) -\varphi \left( r,\bar{X}^{t,x;\bar{u},%
\bar{\xi}}\left( r\right) \right) \\
&&-\varphi _{x}\left( r,\bar{X}^{t,x;\bar{u},\bar{\xi}}\left( r\right)
\right) \Delta \bar{X}^{t,x;\bar{u},\bar{\xi}}\left( r\right) \Big \}.
\end{eqnarray*}%
Taking conditional expectation value $\mathbb{E}^{\mathcal{F}%
_{t_{0}}^{t}}\left( \cdot \right) \left( \omega _{0}\right) ,$ dividing both
sides by $h$, and using (\ref{3.5}), we have
\begin{eqnarray}
&&\mathbb{E}\left[ v\left( s+\theta ,,\bar{X}^{t,x;\bar{u},\bar{\xi}}\left(
s+\theta \right) \right) -v\left( s,,\bar{X}^{t,x;\bar{u},\bar{\xi}}\left(
s\right) \right) \right]  \notag \\
&\geq &\mathbb{E}\left[ \varphi \left( s+\theta ,,\bar{X}^{t,x;\bar{u},\bar{%
\xi}}\left( s+\theta \right) \right) -\varphi \left( s,,\bar{X}^{t,x;\bar{u},%
\bar{\xi}}\left( s\right) \right) \right]  \notag \\
&=&\mathbb{E}\Bigg [\int_{s}^{s+\theta }\Big [\frac{\partial \varphi }{%
\partial t}\left( r,\bar{X}^{t,x;\bar{u},\bar{\xi}}\left( r\right) \right)
+\left\langle D_{x}\varphi \left( r,\bar{X}^{t,x;\bar{u},\bar{\xi}}\left(
r\right) \right) ,\overline{b}\left( r\right) \right\rangle  \notag \\
&&+\frac{1}{2}\text{tr}\left\{ \overline{\sigma }\left( r\right) ^{\ast
}D_{xx}\varphi \left( r,\bar{X}^{t,x;\bar{u},\bar{\xi}}\left( r\right)
\right) \overline{\sigma }\left( r\right) \right\} \Big ]\mbox{\rm
d}r  \notag \\
&&+\int_{t_{0}}^{t_{0}+h}\varphi _{x}\left( r,\bar{X}^{t,x;\bar{u},\bar{\xi}%
}\left( r\right) \right) G\mathrm{d}\bar{\xi}_{r}+\sum_{s\leq r\leq s+\theta
}\Big \{\varphi \left( r,\bar{X}^{t,x;\bar{u},\bar{\xi}}\left( r\right)
\right)  \notag \\
&&-\varphi \left( r,\bar{X}^{t,x;\bar{u},\bar{\xi}}\left( r\right) \right)
-\varphi _{x}\left( r,\bar{X}^{t,x;\bar{u},\bar{\xi}}\left( r\right) \right)
\Delta \bar{X}^{t,x;\bar{u},\bar{\xi}}\left( r\right) \Big \}\Bigg ]
\label{nec1}
\end{eqnarray}%
We now handle the last two terms. Note that
\begin{equation*}
\Delta \bar{X}^{t,x;\bar{u},\bar{\xi}}\left( r\right) =G\Delta \xi _{r}
\end{equation*}%
and
\begin{equation*}
\bar{X}^{t,x;\bar{u},\bar{\xi}}\left( r+\right) =\bar{X}^{t,x;\bar{u},\bar{%
\xi}}\left( r\right) +\Delta \bar{X}^{t,x;\bar{u},\bar{\xi}}\left( r\right) =%
\bar{X}^{t,x;\bar{u},\bar{\xi}}\left( r\right) +G\Delta \xi _{r}.
\end{equation*}%
Thus%
\begin{eqnarray}
&&-\mathbb{E}\left[ \int_{t_{0}}^{t_{0}+h}\varphi _{x}\left( r,\bar{X}^{t,x;%
\bar{u},\bar{\xi}}\left( r\right) \right) G\mathrm{d}\bar{\xi}_{r}\right] +%
\mathbb{E}\left[ \varphi _{x}\left( r,\bar{X}^{t,x;\bar{u},\bar{\xi}}\left(
r\right) \right) \Delta \bar{X}^{t,x;\bar{u},\bar{\xi}}\left( r\right) %
\right]  \notag \\
&=&-\mathbb{E}\left[ \int_{t_{0}}^{t_{0}+h}\varphi _{x}\left( r,\bar{X}^{t,x;%
\bar{u},\bar{\xi}}\left( r\right) \right) G\mathrm{d}\bar{\xi}_{r}^{c}\right]
\notag \\
&=&\mathbb{E}\left[ \int_{t_{0}}^{t_{0}+h}K\mathrm{d}\bar{\xi}_{r}^{c}\right]
.  \label{v2}
\end{eqnarray}%
We now deal the term
\begin{eqnarray}
&&-\mathbb{E}\left[ \sum_{t_{0}\leq r\leq t_{0}+h}\left\{ \varphi \left( r,%
\bar{X}^{t,x;\bar{u},\bar{\xi}}\left( r+\right) \right) -\varphi \left( r,%
\bar{X}^{t,x;\bar{u},\bar{\xi}}\left( r\right) \right) \right\} \right]
\notag \\
&=&-\mathbb{E}\left[ \sum_{t_{0}\leq r\leq t_{0}+h}\left\{
\int_{0}^{1}\varphi _{x}\left( r,\bar{X}^{t,x;\bar{u},\bar{\xi}}\left(
r\right) +\theta \Delta \bar{X}^{t,x;\bar{u},\bar{\xi}}\left( r\right)
\right) G\Delta \bar{\xi}_{r}\mathrm{d}\theta \right\} \right]  \notag \\
&=&\mathbb{E}\left[ K\Delta \bar{\xi}_{r}\right] .  \label{v3}
\end{eqnarray}%
Combining (\ref{v2}) and (\ref{v3}), we have
\begin{eqnarray}
&&\frac{1}{h}\mathbb{E}^{\mathcal{F}_{t_{0}}^{t}\left( \omega _{0}\right) }%
\left[ v\left( t_{0}+h,\bar{X}^{t,x;\bar{u},\bar{\xi}}\left( t_{0}+h\right)
\right) -v\left( t_{0},\bar{X}^{t,x;\bar{u},\bar{\xi}}\left( t_{0}\right)
\right) \right]  \notag \\
&=&\frac{1}{h}\mathbb{E}^{\mathcal{F}_{t_{0}}^{t}\left( \omega _{0}\right) }%
\Bigg \{\int_{t_{0}}^{t_{0}+h}\Big [\frac{\partial \varphi }{\partial t}%
\left( r,\bar{X}^{t,x;\bar{u},\bar{\xi}}\left( r\right) \right)
+\left\langle D_{x}\varphi \left( r,\bar{X}^{t,x;\bar{u},\bar{\xi}}\left(
r\right) \right) ,\overline{b}\left( r\right) \right\rangle  \notag \\
&&+\frac{1}{2}\text{tr}\left\{ \overline{\sigma }\left( r\right) ^{\ast
}D_{xx}\varphi \left( r,\bar{X}^{t,x;\bar{u},\bar{\xi}}\left( r\right)
\right) \overline{\sigma }\left( r\right) \right\} \Big ]\mbox{\rm
d}r-\int_{t_{0}}^{t_{0}+h}K\mathrm{d}\bar{\xi}_{r}\Bigg \}.
\end{eqnarray}%
Letting $h\rightarrow 0,$ and employing the similar delicate method as in
the proof of Theorem 4.1 of Gozzi et al. \cite{GSZ}, we have
\begin{eqnarray*}
&&\frac{1}{h}\limsup_{h\rightarrow 0+}\mathbb{E}^{\mathcal{F}%
_{t_{0}}^{t}\left( \omega _{0}\right) }\left[ v\left( t_{0}+h,\bar{X}^{t,x;%
\bar{u},\bar{\xi}}\left( t_{0}+h\right) \right) -v\left( t_{0},\bar{X}^{t,x;%
\bar{u},\bar{\xi}}\left( t_{0}\right) \right) \right] \\
&\leq &\frac{\partial \varphi }{\partial t}\left( t_{0},\bar{X}^{t,x;\bar{u},%
\bar{\xi}}\left( t_{0},\omega _{0}\right) \right) +\left\langle D_{x}\varphi
\left( t_{0},\bar{X}^{t,x;\bar{u},\bar{\xi}}\left( t_{0},\omega _{0}\right)
\right) ,\overline{b}\left( t_{0}\right) \right\rangle \\
&&+\frac{1}{2}\text{tr}\left\{ \overline{\sigma }\left( t_{0}\right) ^{\ast
}D_{xx}\varphi \left( t_{0},\bar{X}^{t,x;\bar{u},\bar{\xi}}\left(
t_{0},\omega _{0}\right) \right) \overline{\sigma }\left( t_{0}\right)
\right\} -K\mathrm{d}\bar{\xi}_{t_{0}} \\
&=&\overline{p}\left( t_{0},\omega _{0}\right) +\left\langle \overline{q}%
\left( t_{0},\omega _{0}\right) ,\overline{b}\left( t_{0}\right)
\right\rangle +\frac{1}{2}\text{tr}\left\{ \overline{\sigma }\left(
t_{0}\right) ^{\ast }\overline{\Theta }\left( t_{0},\omega _{0}\right)
\overline{\sigma }\left( t_{0}\right) \right\} -K\mathrm{d}\bar{\xi}_{t_{0}}
\end{eqnarray*}%
From Lemma \ref{l6}, that there exist $\rho \in L^{1}\left( t_{0},T;\mathbb{R%
}\right) $ and $\rho _{1}\in L^{1}\left( \Omega ;\mathbb{R}\right) $ such
that
\begin{eqnarray}
&&\mathbb{E}\left[ \frac{1}{h}\left[ v\left( t+h,\bar{X}^{t,x;\bar{u},\bar{%
\xi}}\left( t+h\right) \right) -v\left( t,\bar{X}^{s,y;u,\xi }\left(
t\right) \right) \right] \right]  \notag \\
&\leq &\rho \left( t\right) ,\text{ for }h\leq h_{0}\text{, for some }h_{0}>0
\label{3.8}
\end{eqnarray}%
and
\begin{eqnarray}
&&\mathbb{E}^{\mathcal{F}_{t_{0}}^{t}\left( \omega _{0}\right) }\left[ \frac{%
1}{h}\left[ v\left( t+h,\bar{X}^{t,x;\bar{u},\bar{\xi}}\left( t+h\right)
\right) -v\left( t,\bar{X}^{s,y;u,\xi }\left( t\right) \right) \right] %
\right]  \notag \\
&\leq &\rho _{1}\left( \omega _{0}\right) ,\text{ }h\leq h_{0}\text{, for
some }h_{0}>0.  \label{3.9}
\end{eqnarray}%
holds, respectively.

By virtue of Fatou's Lemma, noting (\ref{3.9}), we obtain
\begin{eqnarray}
&&\limsup_{h\rightarrow 0+}\frac{1}{h}\mathbb{E}\left[ v\left( t_{0}+h,\bar{X%
}^{t,x;\bar{u},\bar{\xi}}\left( t_{0}+h\right) \right) -v\left( t_{0},\bar{X}%
^{t,x;\bar{u},\bar{\xi}}\left( t_{0}\right) \right) \right]  \notag \\
&=&\limsup_{h\rightarrow 0+}\frac{1}{h}\mathbb{E}\left[ \mathbb{E}^{\mathcal{%
F}_{t_{0}}^{t}\left( \omega _{0}\right) }\left\{ v\left( t_{0}+h,\bar{X}%
^{t,x;\bar{u},\bar{\xi}}\left( t_{0}+h\right) \right) -v\left( t_{0},\bar{X}%
^{t,x;\bar{u},\bar{\xi}}\left( t_{0}\right) \right) \right\} \right]  \notag
\\
&\leq &\mathbb{E}\left[ \limsup_{h\rightarrow 0+}\frac{1}{h}\mathbb{E}^{%
\mathcal{F}_{t_{0}}^{t}\left( \omega _{0}\right) }\left\{ v\left( t_{0}+h,%
\bar{X}^{t,x;\bar{u},\bar{\xi}}\left( t_{0}+h\right) \right) -v\left( t_{0},%
\bar{X}^{t,x;\bar{u},\bar{\xi}}\left( t_{0}\right) \right) \right\} \right]
\notag \\
&\leq &\mathbb{E}\left[ \overline{p}\left( t_{0}\right) +\left\langle
\overline{q}\left( t_{0}\right) ,\overline{b}\left( t_{0}\right)
\right\rangle +\frac{1}{2}\text{tr}\left\{ \overline{\sigma }\left(
t_{0}\right) ^{\ast }\overline{\Theta }\left( t_{0}\right) \overline{\sigma }%
\left( t_{0}\right) \right\} -K\mathrm{d}\bar{\xi}_{t_{0}}\right] ,
\label{3.10}
\end{eqnarray}%
for a.e. $t_{0}\in \left[ t,T\right] .$ Then the rest of the proof goes
exactly as in \cite{GSZ}.

We apply Lemma 8 in \cite{zlqverif} to $g\left( t\right) =\mathbb{E}\left[
v\left( t,\bar{X}^{t,x;\bar{u},\bar{\xi}}\left( t\right) \right) \right] ,$
using (\ref{3.8}), \ref{3.6}) and (\ref{3.10}) to get
\begin{eqnarray*}
&&\mathbb{E}\left[ v\left( T,\bar{X}^{t,x;\bar{u},\bar{\xi}}\left( T\right)
\right) -v\left( t,x\right) \right] \\
&\leq &\mathbb{E}\left[ \int_{t}^{T}\overline{p}\left( t\right)
+\left\langle \overline{q}\left( t\right) ,\overline{b}\left( t\right)
\right\rangle +\frac{1}{2}\text{tr}\left[ \overline{\sigma }\left( t\right)
^{\ast }\overline{\Theta }\left( t\right) \overline{\sigma }\left( t\right) %
\right] \mbox{\rm
d}t-K\mathrm{d}\bar{\xi}_{t_{0}}\right] \\
&\leq &-\mathbb{E}\left[ \int_{t}^{T}\overline{f}\left( t\right)
\mbox{\rm
d}t+K\mathrm{d}\bar{\xi}_{t_{0}}\right] .
\end{eqnarray*}%
From this we claim that
\begin{eqnarray*}
v\left( t,x\right) &\geq &\mathbb{E}\left[ v\left( T,\bar{X}^{t,x;\bar{u},%
\bar{\xi}}\left( T\right) \right) +\int_{t}^{T}\overline{f}\left( r\right)
\mbox{\rm
d}r+\int_{s}^{T}K\mathrm{d}\bar{\xi}_{r}\right] \\
&=&\mathbb{E}\left[ \Phi \left( \bar{X}^{t,x;\bar{u},\bar{\xi}}\left(
T\right) \right) +\int_{t}^{T}\overline{f}\left( r\right) \mbox{\rm
d}r+\int_{s}^{T}K\mathrm{d}\bar{\xi}_{r}\right] ,
\end{eqnarray*}%
where%
\begin{equation*}
\bar{f}\left( r\right) =f\left( r,\bar{X}^{t,x;\bar{u},\bar{\xi}}\left(
r\right) ,v\left( r,\bar{X}^{t,x;\bar{u},\bar{\xi}}\left( r\right) \right) ,%
\bar{q}\left( r\right) \sigma \left( r,\bar{X}^{t,x;\bar{u},\bar{\xi}}\left(
r\right) ,\bar{u}\right) \right) .
\end{equation*}%
Thus, combining the above with the first assertion (\ref{3.3}), we prove the
$\left( \bar{X}^{t,x;\bar{u},\bar{\xi}}\left( \cdot \right) ,\overline{u}%
\left( \cdot \right) \right) $ is an optimal pair. The proof is thus
completed. \hfill $\Box $

\begin{remark}
The condition (\ref{3.6}) is just equivalent to the following:
\begin{eqnarray}
\overline{p}\left( s\right) &=&\min\limits_{u\in U}\mathcal{G}\left( s,\bar{X%
}^{t,x;\bar{u},\bar{\xi}}\left( s\right) ,\overline{\varphi }\left( s\right)
,\overline{q}\left( s\right) ,\overline{\Theta }\left( s\right) ,u\right)
\notag \\
&=&\mathcal{G}\left( s,\bar{X}^{t,x;\bar{u},\bar{\xi}}\left( s\right) ,%
\overline{\varphi }\left( s\right) ,\overline{q}\left( s\right) ,\overline{%
\Theta }\left( s\right) ,\overline{u}\left( s\right) \right) ,  \notag \\
\text{a.e. }s &\in &\left[ t,T\right] ,\text{ }P\text{-a.s.,}  \label{3.11}
\end{eqnarray}%
where $\overline{\varphi }\left( t\right) $ is defined in Theorem \ref{verif}%
. This is easily seen by recalling the fact that $v$ is the viscosity
solution of (\ref{pde1}):
\begin{equation*}
\overline{p}\left( s\right) +\min\limits_{u\in U}\mathcal{G}\left( s,\bar{X}%
^{t,x;\bar{u},\bar{\xi}}\left( s\right) ,\overline{\varphi }\left( s\right) ,%
\overline{q}\left( s\right) ,\overline{\Theta }\left( s\right) ,u\right)
\geq 0,
\end{equation*}%
which yields (\ref{3.11}) under (\ref{3.6}).
\end{remark}

\begin{remark}
Clearly, Theorem \ref{verif} is expressed in terms of parabolic superjet.
One could naturally ask whether a similar result holds for parabolic subjet.
The answer was positive for the deterministic case (in terms of the
first-order parabolic subjet, see Theorem 3.9 in \cite{YZ}). Unfortunately,
as claimed in Yong and Zhou \cite{YZ}, the answer is that the statement of
Theorem \ref{verif} is no longer valid whenever the parabolic superjet in (%
\ref{3.5}) is replaced by the parabolic subjet.
\end{remark}

Now let us present a non-smooth version of the necessity part of Theorem \ref%
{verif}. However, we just have \textquotedblleft partial\textquotedblright\
result.

\begin{theorem}
\label{verifnec}Assume that \emph{(A1)-(A2)} hold. Let $v\in C\left( \left[
0,T\right] \times \mathbb{R}^{n}\right) $ be a viscosity solution of the
H-J-B equations (\ref{pde1}) and let $\left( \bar{u}\left( \cdot \right) ,%
\bar{\xi}\left( \cdot \right) \right) $ be an optimal singular controls. Let
$\left( \bar{X}^{t,x;\bar{u},\bar{\xi}}\left( \cdot \right) ,\bar{Y}^{t,x;%
\bar{u},\bar{\xi}}\left( \cdot \right) ,\bar{Z}^{t,x;\bar{u},\bar{\xi}%
}\left( \cdot \right) ,\bar{u}\left( \cdot \right) ,\bar{\xi}\left( \cdot
\right) \right) $ be an admissible pair such that there exist a function $%
\varphi \in C^{1,2}\left( \left[ 0,T\right] ;\mathbb{R}^{n}\right) $ and a
triple
\begin{equation*}
\left( \bar{p},\bar{q},\bar{\Theta}\right) \in \left( L_{\mathcal{F}%
_{t}}^{2}\left( \left[ t,T\right] ;\mathbb{R}\right) \times L_{\mathcal{F}%
_{t}}^{2}\left( \left[ t,T\right] ;\mathbb{R}^{n}\right) \times L_{\mathcal{F%
}_{t}}^{2}\left( \left[ t,T\right] ;\mathbf{S}^{n}\right) \right)
\end{equation*}%
satisfying
\begin{equation}
\left\{
\begin{array}{l}
\left( \bar{p}\left( s\right) ,\bar{q}\left( s\right) ,\bar{\Theta}\left(
s\right) \right) \in \mathcal{P}^{2,-}v\left( s,\bar{X}^{t,x;\bar{u},\bar{\xi%
}}\left( s\right) \right) , \\
\bar{p}\left( s\right) G+K\geq 0,\text{ a.e. }s\in \left[ t,T\right] ,\text{
}P\text{-a.s.}%
\end{array}%
\right.  \label{vernec}
\end{equation}%
Then, it holds that
\begin{equation*}
\mathbb{E}\bar{p}\left( s\right) \leq -\mathbb{E}\left[ \mathcal{G}\left( s,%
\bar{X}^{t,x;\bar{u},\bar{\xi}}\left( s\right) ,\bar{q}\left( s\right) ,\bar{%
\Theta}\left( t\right) \right) \right] ,\text{ a.e. }s\in \left[ t,T\right] .
\end{equation*}
\end{theorem}

\paragraph{Proof.}

On the one hand, let $s\in \left[ t,T\right] $ and $\omega \in \Omega $ such
that
\begin{equation*}
\left( \bar{p}\left( s\right) ,\bar{q}\left( s\right) ,\bar{\Theta}\left(
s\right) \right) \in \mathcal{P}^{2,-}v\left( s,\bar{X}^{t,x;\bar{u},\bar{\xi%
}}\left( s\right) \right) .
\end{equation*}
By Lemma \ref{text}, we have a test function $\varphi \in C^{1,2}\left( %
\left[ 0,T\right] \times \mathbb{R}^{n}\right) $ with $\left( s,x\right) \in %
\left[ 0,T\right] \times \mathbb{R}^{n}$ $\left( p,q,\Theta \right) \in
\mathbb{R}\times \mathbb{R}^{n}\times \mathbf{S}^{n}$ such that $v-\varphi $
achieves its minimum at $\left( s,\bar{X}^{t,x;\bar{u},\bar{\xi}}\left(
s\right) \right) $ and
\begin{equation*}
\left( \frac{\partial \varphi }{\partial t}\left( s,\bar{X}^{t,x;\bar{u},%
\bar{\xi}}\left( s\right) \right) ,D_{x}\varphi \left( s,\bar{X}^{t,x;\bar{u}%
,\bar{\xi}}\left( s\right) \right) ,D^{2}\varphi \left( s,\bar{X}^{t,x;\bar{u%
},\bar{\xi}}\left( s\right) \right) \right) =\left( \bar{p}\left( s\right) ,%
\bar{q}\left( s\right) ,\bar{\Theta}\left( s\right) \right)
\end{equation*}%
holds. Then for sufficiently small $\theta >0$, a.e. $s\in \left[ t,T\right]
$.%
\begin{eqnarray*}
&&\mathbb{E}\left[ v\left( s+\theta ,,\bar{X}^{t,x;\bar{u},\bar{\xi}}\left(
s+\theta \right) \right) -v\left( s,,\bar{X}^{t,x;\bar{u},\bar{\xi}}\left(
s\right) \right) \right] \\
&\geq &\mathbb{E}\left[ \varphi \left( s+\theta ,,\bar{X}^{t,x;\bar{u},\bar{%
\xi}}\left( s+\theta \right) \right) -\varphi \left( s,,\bar{X}^{t,x;\bar{u},%
\bar{\xi}}\left( s\right) \right) \right] \\
&=&\mathbb{E}\Bigg [\int_{s}^{s+\theta }\Big [\frac{\partial \varphi }{%
\partial t}\left( r,\bar{X}^{t,x;\bar{u},\bar{\xi}}\left( r\right) \right)
+\left\langle D_{x}\varphi \left( r,\bar{X}^{t,x;\bar{u},\bar{\xi}}\left(
r\right) \right) ,\overline{b}\left( r\right) \right\rangle \\
&&+\frac{1}{2}\text{tr}\left\{ \overline{\sigma }\left( r\right) ^{\ast
}D_{xx}\varphi \left( r,\bar{X}^{t,x;\bar{u},\bar{\xi}}\left( r\right)
\right) \overline{\sigma }\left( r\right) \right\} \Big ]\mbox{\rm
d}r \\
&&+\int_{t_{0}}^{t_{0}+h}\varphi _{x}\left( r,\bar{X}^{t,x;\bar{u},\bar{\xi}%
}\left( r\right) \right) G\mathrm{d}\bar{\xi}_{r}+\sum_{s\leq r\leq s+\theta
}\Big \{\varphi \left( r,\bar{X}^{t,x;\bar{u},\bar{\xi}}\left( r\right)
\right) \\
&&-\varphi \left( r,\bar{X}^{t,x;\bar{u},\bar{\xi}}\left( r\right) \right)
-\varphi _{x}\left( r,\bar{X}^{t,x;\bar{u},\bar{\xi}}\left( r\right) \right)
\Delta \bar{X}^{t,x;\bar{u},\bar{\xi}}\left( r\right) \Big \}\Bigg ] \\
&\geq &\mathbb{E}\Bigg [\int_{s}^{s+\theta }\Big [\frac{\partial \varphi }{%
\partial t}\left( r,\bar{X}^{t,x;\bar{u},\bar{\xi}}\left( r\right) \right)
+\left\langle D_{x}\varphi \left( r,\bar{X}^{t,x;\bar{u},\bar{\xi}}\left(
r\right) \right) ,\overline{b}\left( r\right) \right\rangle \\
&&+\frac{1}{2}\text{tr}\left\{ \overline{\sigma }\left( r\right) ^{\ast
}D_{xx}\varphi \left( r,\bar{X}^{t,x;\bar{u},\bar{\xi}}\left( r\right)
\right) \overline{\sigma }\left( r\right) \right\} \Big ]\mbox{\rm
d}r-\int_{s}^{s+\theta }K\mathrm{d}\bar{\xi}_{r}\Bigg ].
\end{eqnarray*}%
The last inequality comes from the derivation in Theorem \ref{verif} by
means of the condition (\ref{vernec}). On the other hand, since $\left( \bar{%
X}^{t,x;\bar{u},\bar{\xi}}\left( \cdot \right) ,\bar{Y}^{t,x;\bar{u},\bar{\xi%
}}\left( \cdot \right) ,\bar{Z}^{t,x;\bar{u},\bar{\xi}}\left( \cdot \right)
\bar{u}\left( \cdot \right) ,\bar{\xi}\left( \cdot \right) \right) $ is
optimal, by DPP of optimality, it yields%
\begin{eqnarray*}
&&v\left( \tau ,,\bar{X}^{t,x;\bar{u},\bar{\xi}}\left( \tau \right) \right)
\\
&=&\mathbb{E}^{\mathcal{F}_{s}^{t}\left( \omega \right) }\left[ \Phi \left(
\bar{X}^{t,x;\bar{u},\bar{\xi}}\left( T\right) \right) +\int_{\tau }^{T}\bar{%
f}\left( r\right) \mathrm{d}r+\int_{\tau }^{T}K\mathrm{d}\bar{\xi}_{r}\right]
,\text{ }\forall \tau \in \left[ t,T\right] ,\text{ }P\text{-a.s.,}
\end{eqnarray*}%
which implies that
\begin{equation}
\mathbb{E}\left[ v\left( s+\theta ,,\bar{X}^{t,x;\bar{u},\bar{\xi}}\left(
s+\theta \right) \right) -v\left( s,,\bar{X}^{t,x;\bar{u},\bar{\xi}}\left(
s\right) \right) \right] =-\left[ \int_{s}^{s+\theta }\bar{f}\left( r\right)
\mathrm{d}r+\int_{\tau }^{T}K\mathrm{d}\bar{\xi}_{r}\right] .  \label{nec2}
\end{equation}%
Therefore, it follows from (\ref{nec2}) that
\begin{equation*}
\begin{array}{c}
\mathbb{E}\left[ \bar{p}\left( s\right) \right] \leq -\mathbb{E}\Big [\frac{%
\partial \varphi }{\partial t}\left( s,\bar{X}^{t,x;\bar{u},\bar{\xi}}\left(
s\right) \right) +\left\langle D_{x}\varphi \left( s,\bar{X}^{t,x;\bar{u},%
\bar{\xi}}\left( s\right) \right) ,\overline{b}\left( s\right) \right\rangle
\\
+\frac{1}{2}\text{tr}\left\{ \overline{\sigma }\left( s\right) ^{\ast
}D_{xx}\varphi \left( s,\bar{X}^{t,x;\bar{u},\bar{\xi}}\left( s\right)
\right) \overline{\sigma }\left( s\right) \right\} -\bar{f}\left( s\right) %
\Big ],\text{ a.e. }s\in \left[ t,T\right] ,%
\end{array}%
\end{equation*}%
where
\begin{equation*}
\bar{f}\left( r\right) =f\left( r,\bar{X}^{t,x;\bar{u},\bar{\xi}}\left(
r\right) ,v\left( r,\bar{X}^{t,x;\bar{u},\bar{\xi}}\left( r\right) \right) ,%
\bar{q}\left( r\right) \sigma \left( r,\bar{X}^{t,x;\bar{u},\bar{\xi}}\left(
r\right) ,\bar{u}\right) \right) .
\end{equation*}
We thus complete the proof. \hfill $\Box $

\subsection{Optimal Feedback Controls}

In this subsection, we describe the method to construct optimal feedback
controls by the verification Theorem \ref{verif}. First, let us recall the
definition of admissible feedback controls.

\begin{definition}
A measurable function $\left( \mathbf{u,\xi }\right) $ from $\left[ 0,T%
\right] \times \mathbb{R}^{n}$ to $U\times \left[ 0,\infty \right) ^{m}$ is
called an admissible feedback control pair if for any $\left( t,x\right) \in %
\left[ 0,T\right) \times \mathbb{R}^{n}$ there is a weak solution $%
X^{t,x;u,\xi }\left( \cdot \right) $ of the following SDE:
\begin{equation}
\left\{
\begin{array}{rcl}
\mbox{\rm
d}X^{t,x;\mathbf{u},\mathbf{\xi }}\left( r\right) & = & b\left( r,X^{t,x;%
\mathbf{u},\mathbf{\xi }}\left( r\right) ,\mathbf{u}\left( r\right) \right)
\mbox{\rm
d}r+\sigma \left( r,X^{t,x;\mathbf{u},\mathbf{\xi }}\left( r\right) ,\mathbf{%
u}\left( r\right) \right) \mbox{\rm
d}W\left( r\right) +G\mathrm{d}\mathbf{\xi }_{r}, \\
\mbox{\rm
d}Y^{t,x;\mathbf{u},\mathbf{\xi }}\left( r\right) & = & -f\left( r,X^{t,x;%
\mathbf{u},\mathbf{\xi }}\left( r\right) ,Y^{t,x;\mathbf{u},\mathbf{\xi }%
}\left( r\right) ,\mathbf{u}\left( r\right) \right) \mbox{\rm
d}r+\mbox{\rm
d}M^{t,x;\mathbf{u},\mathbf{\xi }}\left( r\right) -K\mathrm{d}\mathbf{\xi }%
_{r}, \\
X^{t,x;\mathbf{u},\mathbf{\xi }}\left( t\right) & = & x,\quad Y^{t,x;\mathbf{%
u},\mathbf{\xi }}\left( T\right) =\Phi \left( X^{t,x;\mathbf{u},\mathbf{\xi }%
}\left( T\right) \right) ,\text{ }r\in \left[ t,T\right] ,%
\end{array}%
\right.  \label{4.1}
\end{equation}%
where $M^{t,x;\mathbf{u},\mathbf{\xi }}$ is an $\mathbb{R}$-valued $\mathbb{F%
}^{t,x;\mathbf{u},\mathbf{\xi }}$-adapted right continuous and left limit
martingale vanishing in $t=0$ which is orthogonal to the driving Brownian
motion $W.$ Here $\mathbb{F}^{t,x;\mathbf{u},\mathbf{\xi }}=\left( \mathcal{F%
}_{s}^{X^{t,x;\mathbf{u},\mathbf{\xi }}}\right) _{s\in \left[ t,T\right] }$
is the smallest filtration and generated by $X^{t,x;\mathbf{u},\mathbf{\xi }%
} $, which is such that $X^{t,x;\mathbf{u},\mathbf{\xi }}$ is $\mathbb{F}%
^{t,x;\mathbf{u},\mathbf{\xi }}$-adapted. Obviously, $M^{t,x;\mathbf{u},%
\mathbf{\xi }}$ is a part of the solution of BSDE of (\ref{4.1}).
Simultaneously, we suppose that $f$ satisfies the Lipschitz condition with
respect to $\left( x,y,z\right) $. An admissible feedback control pair $%
\left( \mathbf{u}^{\star },\mathbf{\xi }^{\star }\right) $ is called optimal
if
\begin{equation*}
(X^{\star }\left( \cdot ;t,x\right) ,Y^{\star }\left( \cdot ;t,x\right) ,%
\mathbf{u}^{\star }\left( \cdot ,X^{\star }\left( \cdot ;t,x\right) \right) ,%
\mathbf{\xi }^{\star }\left( \cdot ,X^{\star }\left( \cdot ;t,x\right)
\right) )
\end{equation*}
is optimal for each $\left( t,x\right) $ is a solution of (\ref{4.1})
corresponding to $\left( \mathbf{u}^{\star },\mathbf{\xi }^{\star }\right) .$
\end{definition}

\begin{theorem}
\label{t3}Let $\left( \mathbf{u}^{\star },\mathbf{\xi }^{\star }\right) $ be
an admissible feedback control and $p^{\star },q^{\star },$ and $\Theta
^{\star }$ be measurable functions satisfying $\left( p^{\star }\left(
t,x\right) ,q^{\star }\left( t,x\right) ,\Theta \left( t,x\right) \right)
\in \mathcal{P}^{2,+}v\left( t,x\right) $ for all $\left( t,x\right) \in %
\left[ 0,T\right] \times \mathbb{R}^{n}.$ If
\begin{eqnarray}
&&p^{\star }\left( t,x\right) +\mathcal{G}\left( t,x,V\left( t,x\right)
,q^{\star }\left( t,x\right) ,\Theta ^{\star }\left( t,x\right) ,\mathbf{u}%
^{\star }\left( t,x\right) \right)  \notag \\
&=&\inf\limits_{\left( p,q,\Theta ,u\right) \in \mathcal{P}^{2,+}v\left(
t,x\right) \times U}\left[ p+\mathcal{G}\left( t,x,V\left( t,x\right)
,q,\Theta ,u\right) \right] =0  \label{4.3}
\end{eqnarray}%
and $q^{\star }\left( t,x\right) G+K\geq 0$ for all $\left( t,x\right) \in %
\left[ 0,T\right] \times \mathbb{R}^{n},$ then $\left( \mathbf{u}^{\star },%
\mathbf{\xi }^{\star }\right) $ is singular optimal control pair.
\end{theorem}

\paragraph{Proof}

From Theorem \ref{verif}, we get the desired result. \hfill $\Box $

\begin{remark}
In FBSDEs (\ref{4.1}), $Y^{t,x;u}\left( \cdot \right) $ is actually
determined by $\left( X^{t,x;u}\left( \cdot \right) ,u\left( \cdot \right)
,\xi \left( \cdot \right) \right) .$ Hence, we need to investigate the
conditions imposed in Theorem \ref{verif} to ensure the existence and
uniqueness of $X^{t,x;u}\left( \cdot \right) $ in law and the measurability
of the multifunctions $\left( t,x\right) \rightarrow \mathcal{P}%
^{2,+}v\left( t,x\right) $ to obtain $\left( p^{\star }\left( t,x\right)
,q^{\star }\left( t,x\right) ,\Theta \left( t,x\right) \right) \in \mathcal{P%
}^{2,+}v\left( t,x\right) $ that minimizes (\ref{4.3}). This can be done by
virtue of the celebrated Filippov's Lemma (cf \cite{YZ}).
\end{remark}

\subsection{The Connection between DPP and MP}

\label{secconne}In Section \ref{sec3}, we have obtained the first and second
order adjoint equations. In this part, we shall investigate the connection
between the general DPP and the MP for such singular controls problem
without the assumption that the value is sufficient smooth. By associated
adjoint equations and delicate estimates, it is possible to establish the
set inclusions among the super- and sub-jets of the value function and the
first-order and second- order adjoint processes as well as the generalized
Hamiltonian function.

\begin{theorem}
\label{MPDP}Assume that \emph{(A1)-(A2)} are in force. Suppose that $\left(
\bar{u},\bar{\xi}\right) $ be a singular optimal controls, $v\left( \cdot
,\cdot \right) $ is a value function, and $\left( \bar{X}^{t,x;\bar{u},\bar{%
\xi}}\left( \cdot \right) ,\bar{Y}^{t,x;\bar{u},\bar{\xi}}\left( \cdot
\right) ,\bar{Z}^{t,x;\bar{u},\bar{\xi}}\left( \cdot \right) ,\bar{u}\left(
\cdot \right) ,\bar{\xi}\left( \cdot \right) \right) $ is optimal
trajectory. Let $\left( p,q\right) \in \mathcal{S}^{2}(0,T;\mathbb{R}%
^{n})\times \mathcal{M}^{2}(0,T;\mathbb{R}^{n})$ and $\left( P,Q\right) \in
\mathcal{S}^{2}(0,T;\mathbb{R}^{n\times n})\times \mathcal{M}^{2}(0,T;%
\mathbb{R}^{n\times n})$ be the adjoint equations (\ref{adj1}), (\ref{adj2}%
), respectively. Then, we have%
\begin{equation}
\begin{array}{c}
P\left\{ K_{\left( i\right) }+p\left( t\right) G_{\left( i\right) }\left(
t\right) \geq 0,\text{ }t\in \left[ 0,T\right] ,\text{ }\forall i\right\} =1,
\\
\left\{ p\left( s\right) \right\} \times \left[ P\left( s\right) ,\infty
\right) \subseteq \mathcal{P}^{2,+}v\left( t,\bar{X}^{t,x;\bar{u},\bar{\xi}%
}\left( s\right) \right) , \\
\mathcal{P}^{2,-}v\left( t,\bar{X}^{t,x;\bar{u},\bar{\xi}}\left( s\right)
\right) \subseteq \left\{ p\left( s\right) \right\} \times \left[ -\infty
,P\left( s\right) \right) ,\text{ a.e. }s\in \left[ t,T\right] ,\text{
P-a.s..}%
\end{array}
\label{MPDPP1}
\end{equation}
\end{theorem}

\paragraph{Proof.}

From Theorem \ref{mpsing} and Proposition \ref{pro1}, we get the first part
of (\ref{MPDPP1}). From Theorem 3.1 in Nie, Shi and Wu \cite{NSW}, we get
the second and third results of (\ref{MPDPP1}).\hfill $\Box $

\section{Concluding remarks}

\label{sec5}

In this paper, on the one hand, we have derived a second order pointwise
necessary condition for singular optimal control in classical sense of
FBSDEs with convex control domain by means of the variation equations and
two adjoint equations, which is separately extends the work by Zhang and
Zhang \cite{ZZconvex} to stochastic recursive case, and Hu \cite{Hms} to
pointwise case in the framework of Malliavin calculus. A new necessary
condition for singular control has been obtained. Moreover, we investigate
the verification theorem for optimal controls via viscosity solution and
establish the connection between the adjoint equations and value function
also in viscosity solution sense.

There are still several interesting topics should be scheduled as follows:

\begin{itemize}
\item As an important issue, the \emph{existence} of optimal singular
controls has never been exploited. Haussmannand and Suo \cite{HS1} apply the
compactification method to study the classical and singular control problem
of It\^{o}'s type of stochastic differential equation, where the problem is
reformulated as a martingale problem on an appropriate canonical space after
the relaxed form of the classical control is introduced. Under some mild
continuity assumptions on the data, they obtain the existence of optimal
control by purely probabilistic arguments. Note that, in the framework of
BSDE with singular control, the trajectory of $Y$ seems to be a c\`{a}dl\`{a}%
g process (from French, for right continuous with left hand limits). Hence,
we may consider $Y$ in some space with appropriate topologies, for instance,
Skorokhod $M_{1}$ topology or Meyer-Zheng topology (see \cite{FH}) to obtain
the convergence of probability measures deduced by $Y$ involving relaxed
control. Related work from the technique of PDEs can be seen in \cite{BLRT,
BGM} references therein. From Wang \cite{WB}, one may construct the optimal
control via the existence of diffusion with refections (see \cite{CMR}).
However, it is interesting to extend this result to FBSDEs.

\item The matrices $K,G$ are deterministic. It is also interesting to extend
this restriction to time varying matrices, even the generator $b,\sigma ,f$
involving the singular control. Whenever the coefficients are random, the
H-J-B inequality will become stochastic PDEs. No doubt, stochastic viscosity
solution will be applied. For this direction, reader can refer to Buckdahn,
Ma \cite{BM1, BM2}, Peng \cite{PengSHJB} and Qiu \cite{Qiu}.

\item As for the general cases, \textit{i.e.}, the control regions are
assumed to be non-convex and both the drift and diffusion terms depend on
the control variable. Indeed, such a mathematical model, from view point of
application, is more reasonable and urgent in many real-life problems (for
instance, some finance models in which the controls may impact the
uncertainty, \textit{etc}). In near future, we shall remove the condition of
convex control region, employing the idea developed by Zhang et al. \cite%
{ZZgen}. It is worth mentioning that the analysis in \cite{ZZgen} is much
more complicated. Some new and useful tools, such as the multilinear
function valued stochastic processes, the BSDE for these processes are
introduced. Hence, it will be interesting to borrow these tools to
investigate the singular optimal controls problems for FBSDEs, which will
definitely promote and enrich the theories of FBSDEs.
\end{itemize}

\appendix\label{APP}

\section{Proofs of Lemmas}

\paragraph{Proof of Lemma \protect\ref{estvs}.}

We first prove the continuity of solution depending on parameter.

Set%
\begin{eqnarray*}
\hat{X}^{\alpha }\left( s\right) &=&X^{0,x;\bar{u},\xi ^{\alpha }}\left(
s\right) -X^{0,x;\bar{u},\bar{\xi}}\left( s\right) , \\
\hat{Y}^{\alpha }\left( s\right) &=&Y^{0,x;\bar{u},\xi ^{\alpha }}\left(
s\right) -Y^{0,x;\bar{u},\bar{\xi}}\left( s\right) , \\
\hat{Z}^{\alpha }\left( s\right) &=&Z^{0,x;\bar{u},\xi ^{\alpha }}\left(
s\right) -Z^{0,x;\bar{u},\bar{\xi}}\left( s\right) ,\text{ }s\in \left[ 0,T%
\right] .
\end{eqnarray*}%
It can be shown that
\begin{equation}
\left( \hat{X}^{\alpha }\left( t\right) ,\hat{Y}^{\alpha }\left( t\right) ,%
\hat{Z}^{\alpha }\left( t\right) \right) \text{ converges to }0\text{ in }%
\mathcal{N}^{2}\left[ 0,T\right] \text{ as }\alpha \rightarrow 0  \label{e4}
\end{equation}%
by standard estimates and the Burkholder-Davis-Gundy inequality$,$ so we
omit it.

Next, set
\begin{equation*}
\Delta X^{\alpha }\left( s\right) =\frac{\hat{X}^{\alpha }\left( s\right) }{%
\alpha }-x^{1}\left( s\right) ,\text{ }\Delta Y^{\alpha }\left( s\right) =%
\frac{\hat{Y}^{\alpha }\left( s\right) }{\alpha }-y^{1}\left( s\right) ,%
\text{ }\Delta Z^{\alpha }\left( s\right) =\frac{\hat{Z}^{\alpha }\left(
s\right) }{\alpha }-z^{1}\left( s\right) .
\end{equation*}%
Note that (\ref{e1}) has be obtained in \cite{BDM1}. We will prove (\ref{e2}%
) and (\ref{e3}).

Then,%
\begin{equation}
\left\{
\begin{array}{rcl}
\mathrm{d}\Delta X^{\alpha }\left( s\right) & = & \Delta b^{\alpha }\left(
s\right) \mathrm{d}s+\Delta \sigma ^{\alpha }\left( s\right) \mathrm{d}%
W\left( s\right) , \\
-\mathrm{d}\Delta Y^{\alpha }\left( s\right) & = & \Delta f^{\alpha }\left(
s\right) \mathrm{d}s-\Delta Z^{\alpha }\left( s\right) \mathrm{d}W\left(
s\right) , \\
\Delta X^{\alpha }\left( 0\right) & = & 0,\text{ }\Delta Y^{\alpha }\left(
T\right) =\Delta \Phi ^{\alpha }\left( T\right) ,\qquad 0\leq t\leq s\leq T,%
\end{array}%
\right.
\end{equation}%
where
\begin{eqnarray*}
\Delta b^{\alpha }\left( s\right) &=&\frac{b\left( s,X_{s}^{0,x;\bar{u},\xi
^{\alpha }},\bar{u}_{s}\right) -b\left( s,X_{s}^{0,x;\bar{u},\bar{\xi}},\bar{%
u}_{s}\right) }{\alpha }-\bar{b}_{x}\left( t\right) x^{1}\left( t\right) , \\
\Delta \sigma ^{\alpha }\left( s\right) &=&\frac{\sigma \left( s,X_{s}^{0,x;%
\bar{u},\xi ^{\alpha }},\bar{u}_{s}\right) -\sigma \left( s,X_{s}^{0,x;\bar{u%
},\bar{\xi}},\bar{u}_{s}\right) }{\alpha }-\bar{\sigma}_{x}\left( t\right)
x^{1}\left( t\right) , \\
\Delta f^{\alpha }\left( s\right) &=&\frac{f\left( s,X_{s}^{0,x;\bar{u},\xi
^{\alpha }},Y_{s}^{0,x;\bar{u},\xi ^{\alpha }},Z_{s}^{0,x;\bar{u},\xi
^{\alpha }},\bar{u}_{s}\right) -f\left( s,X_{s}^{0,x;\bar{u},\bar{\xi}%
},Y_{s}^{0,x;\bar{u},\bar{\xi}},Z_{s}^{0,x;\bar{u},\bar{\xi}},\bar{u}%
_{s}\right) }{\alpha } \\
&&-\bar{f}_{x}\left( s\right) x^{1}\left( s\right) -\bar{f}_{y}\left(
s\right) y^{1}\left( s\right) -\bar{f}_{z}\left( s\right) z^{1}\left(
s\right) , \\
\Delta \Phi ^{\alpha }\left( T\right) &=&\frac{\Phi \left( X_{T}^{0,x;\bar{u}%
,\xi ^{\alpha }}\right) -\Phi \left( X_{T}^{0,x;\bar{u},\bar{\xi}}\right) }{%
\alpha }-\Phi _{x}\left( X_{T}^{0,x;\bar{u},\bar{\xi}}\right) x^{1}\left(
T\right) .
\end{eqnarray*}%
Simple calculation yields
\begin{eqnarray*}
\Delta f^{\alpha }\left( s\right) &=&\int_{0}^{1}f_{x}\left( s,\lambda
X_{s}^{0,x;\bar{u},\xi ^{\alpha }}+\left( 1-\lambda \right) X_{s}^{0,x;\bar{u%
},\bar{\xi}},Y_{s}^{0,x;\bar{u},\xi ^{\alpha }},Z_{s}^{0,x;\bar{u},\xi
^{\alpha }},\bar{u}_{s}\right) \Delta X^{\alpha }\left( s\right) \mathrm{d}%
\lambda \\
&&+\int_{0}^{1}f_{y}\left( s,X_{s}^{0,x;\bar{u},\bar{\xi}},\lambda
Y_{s}^{0,x;\bar{u},\xi ^{\alpha }}+\left( 1-\lambda \right) Y_{s}^{0,x;\bar{u%
},\bar{\xi}},Z_{s}^{0,x;\bar{u},\xi ^{\alpha }},\bar{u}_{s}\right) \Delta
Y^{\alpha }\left( s\right) \mathrm{d}\lambda \\
&&+\int_{0}^{1}f_{z}\left( s,X_{s}^{0,x;\bar{u},\bar{\xi}},Y_{s}^{0,x;\bar{u}%
,\bar{\xi}},\lambda Z_{s}^{0,x;\bar{u},\xi ^{\alpha }}+\left( 1-\lambda
\right) Z_{s}^{0,x;\bar{u},\bar{\xi}},\bar{u}_{s}\right) \Delta Z^{\alpha
}\left( s\right) \mathrm{d}\lambda \\
&&+\Delta \rho ^{\alpha }\left( s\right) ,
\end{eqnarray*}%
where
\begin{eqnarray*}
\Delta \rho ^{\alpha }\left( s\right) &=&\left[ \int_{0}^{1}f_{x}\left(
s,\lambda X_{s}^{0,x;\bar{u},\xi ^{\alpha }}+\left( 1-\lambda \right)
X_{s}^{0,x;\bar{u},\bar{\xi}},Y_{s}^{0,x;\bar{u},\xi ^{\alpha }},Z_{s}^{0,x;%
\bar{u},\xi ^{\alpha }},\bar{u}_{s}\right) \mathrm{d}\lambda -\bar{f}%
_{x}\left( s\right) \right] x^{1}\left( s\right) \mathrm{d}s \\
&&+\left[ \int_{0}^{1}f_{y}\left( s,X_{s}^{0,x;\bar{u},\bar{\xi}},\lambda
Y_{s}^{0,x;\bar{u},\xi ^{\alpha }}+\left( 1-\lambda \right) Y_{s}^{0,x;\bar{u%
},\bar{\xi}},Z_{s}^{0,x;\bar{u},\xi ^{\alpha }},\bar{u}_{s}\right) \mathrm{d}%
\lambda -\bar{f}_{y}\left( s\right) \right] y^{1}\left( s\right) \mathrm{d}s
\\
&&+\left[ \int_{0}^{1}f_{z}\left( s,X_{s}^{0,x;\bar{u},\bar{\xi}},Y_{s}^{0,x;%
\bar{u},\bar{\xi}},\lambda Z_{s}^{0,x;\bar{u},\xi ^{\alpha }}+\left(
1-\lambda \right) Z_{s}^{0,x;\bar{u},\bar{\xi}},\bar{u}_{s}\right) \mathrm{d}%
\lambda -\bar{f}_{z}\left( s\right) \right] z^{1}\left( s\right) \mathrm{d}s
\end{eqnarray*}%
and
\begin{eqnarray*}
\Delta \Phi ^{\alpha }\left( T\right) &=&\int_{0}^{1}\Phi _{x}\left( \lambda
X_{T}^{0,x;\bar{u},\xi ^{\alpha }}+\left( 1-\lambda \right) X_{T}^{0,x;\bar{u%
},\bar{\xi}}\right) \Delta X^{\alpha }\left( s\right) \mathrm{d}\lambda \\
&&+\left[ \int_{0}^{1}\Phi _{x}\left( \lambda X_{T}^{0,x;\bar{u},\xi
^{\alpha }}+\left( 1-\lambda \right) X_{T}^{0,x;\bar{u},\bar{\xi}}\right)
\mathrm{d}\lambda -\Phi _{x}\left( X_{T}^{0,x;\bar{u},\bar{\xi}}\right) %
\right] x^{1}\left( T\right) \mathrm{d}s.
\end{eqnarray*}
From classical theory of BSDE, one can show
\begin{equation}
\mathbb{E}\left[ \sup_{0\leq t\leq T}\left\vert x^{1}\left( t\right)
\right\vert ^{2}\right] <\infty ,\mathbb{E}\left[ \sup_{0\leq t\leq
T}\left\vert y^{1}\left( t\right) \right\vert ^{2}\mathrm{d}t\right] <\infty
,\mathbb{E}\left[ \int_{0}^{T}\left\vert z^{1}\left( t\right) \right\vert
^{2}\mathrm{d}t\right] <\infty .  \label{e5}
\end{equation}%
By using (\ref{e4}) and (\ref{e5}), the dominated convergence theorem, Lemma %
\ref{estlemma} and Gronwall's lemma, we get the desired result by letting $%
\alpha \rightarrow 0$.\hfill $\Box $%
%\noindent \textbf{Acknowledgments.}

\paragraph{Proof of Lemma \protect\ref{l6}.}

From (\ref{3.2}) and (6) in Gozzi et al. in \cite{GSZ}, we have that if $%
\left( p,q,P\right) \in \mathcal{P}^{2,+}v\left( t,x\right) ,$ then%
\begin{eqnarray*}
&&v\left( t+h,X^{t,x;\bar{u},\bar{\xi}}\left( t+h\right) \right) -v\left(
t,X^{t,x;\bar{u},\bar{\xi}}\left( t\right) \right) \leq C\left( 1+\left\vert
X^{t,x;\bar{u},\bar{\xi}}\left( t\right) \right\vert ^{m}\right) h \\
&&+\left\langle q\left( t\right) ,X^{t,x;\bar{u},\bar{\xi}}\left( t+h\right)
-X^{t,x;\bar{u},\bar{\xi}}\left( t\right) \right\rangle +C_{0}\left\vert
X^{t,x;\bar{u},\bar{\xi}}\left( t+h\right) -X^{t,x;\bar{u},\bar{\xi}}\left(
t\right) \right\vert ^{2} \\
&=&I_{1}+I_{2}+I_{3}.
\end{eqnarray*}%
We shall deal with $I_{1},$ $I_{2},$ $I_{3},$ separately. For $I_{1},$ we
have $\mathbb{E}\left( 1+\left\vert X^{t,x;\bar{u},\bar{\xi}}\left(
t+h\right) \right\vert ^{m}\right) h\leq C\left( 1+\left\vert x\right\vert
^{m}\right) h,$ by classical estimate and the assumption $\mathbb{E}\left[
\left\vert \xi _{T}\right\vert ^{2}\right] <\infty .$ For $I_{2},$ from (7)
in \cite{GSZ} and H\"{o}lder inequality, we have%
\begin{eqnarray*}
&&\mathbb{E}\left\langle q\left( t\right) ,X^{t,x;\bar{u},\bar{\xi}}\left(
t+h\right) -X^{t,x;\bar{u},\bar{\xi}}\left( t\right) \right\rangle \\
&\leq &C\left( \mathbb{E}\left[ \left( 1+\left\vert X^{t,x;\bar{u},\bar{\xi}%
}\left( t\right) \right\vert ^{m_{1}}\right) ^{2}\right] \right) ^{\frac{1}{2%
}}\Bigg \{\left( \mathbb{E}\left[ \left\vert \int_{t}^{t+h}b\left( r,X^{0,x;%
\bar{u},\xi ^{\alpha }}\left( r\right) ,\bar{u}\left( r\right) \right)
\right\vert ^{2}\right] \right) ^{\frac{1}{2}} \\
&&+\left( \mathbb{E}\left[ \left\vert \int_{t}^{t+h}G\left( r\right) \mathrm{%
d}\bar{\xi}\left( r\right) \right\vert ^{2}\right] \right) ^{\frac{1}{2}}%
\Bigg \} \\
&\leq &C\left( \mathbb{E}\left[ \left( 1+\left\vert x\right\vert
^{m_{1}}\right) ^{2}\right] \right) ^{\frac{1}{2}},
\end{eqnarray*}%
since $\mathbb{E}\left[ \left\vert \xi _{T}\right\vert ^{2}\right] <\infty $
and the fact $\left( 1+\left\vert x\right\vert ^{2}\right) ^{\frac{1}{2}%
}\leq 1+\left\vert x\right\vert .$

Finally,
\begin{eqnarray*}
&&C_{0}\mathbb{E}\left\vert X^{t,x;\bar{u},\bar{\xi}}\left( t+h\right)
-X^{t,x;\bar{u},\bar{\xi}}\left( t\right) \right\vert ^{2}\leq C_{0}\mathbb{E%
}\left[ \left\vert \int_{t}^{t+h}b\left( r,X^{0,x;\bar{u},\xi ^{\alpha
}}\left( r\right) ,\bar{u}\left( r\right) \right) \right\vert ^{2}\right] \\
&&+C_{0}\mathbb{E}\left[ \left\vert \int_{t}^{t+h}\sigma \left( r,X^{0,x;%
\bar{u},\xi ^{\alpha }}\left( r\right) ,\bar{u}\left( r\right) \right)
\right\vert ^{2}\right] +C_{0}\mathbb{E}\left[ \left\vert
\int_{t}^{t+h}G\left( r\right) \mathrm{d}\bar{\xi}\left( r\right)
\right\vert ^{2}\right] \\
&\leq &C\left( 1+\left\vert x\right\vert ^{2}\right) \left( 2h^{2}+h\right) .
\end{eqnarray*}%
By It\^{o} isometry and classical estimate on SDE, we complete the
proof.\hfill $\Box $

\end{document}